\documentclass[10pt]{amsart}    
\usepackage{amsmath}
\usepackage{amssymb}
\usepackage{bm}
\usepackage{amsfonts}
\usepackage{enumerate}
\usepackage{amsmath,amsthm,amssymb}
\usepackage{mathrsfs}
\usepackage{color}
\allowdisplaybreaks

\usepackage[colorlinks, linkcolor=blue, anchorcolor=blue, citecolor=blue, bookmarksnumbered, bookmarksopen=true]{hyperref} 

\usepackage{geometry}
\newcommand{\R}{\mathbb{R}}

\newcommand{\B}{{\bf B}}
\newcommand{\D}{{\bf D}}
\newcommand{\G}{{\bf G}}
\renewcommand{\S}{{\bf S}}
\newcommand{\X}{{\bf X}}
\newcommand{\Y}{{\bf Y}}

\newcommand{\T}{\T}

\renewcommand{\div}{{\rm div}}
\newcommand{\curl}{\nabla\!\times\!} 

\newcommand{\pd}{\partial}

\newcommand{\sdot}{\!\cdot\!}

\def\norm#1{\left\| {#1} \right\|}  
\def\md#1{\partial_t{#1} + u \!\cdot\! \nabla{#1}}  
\def\rmd#1{\rho\partial_t{#1} + \rho u \!\cdot\! \nabla{#1}}

\def\a{\alpha}
\def\eps{\varepsilon}

\def\s{\sigma}

\def\T{\mathbb{T}}

\theoremstyle{plain} 
\newtheorem{theorem}{Theorem}[section]

\newtheorem{lemma}{Lemma}[section]
\newtheorem{proposition}{Proposition}[section]

\theoremstyle{remark} 
\newtheorem{remark}{\bf Remark}[section]

\theoremstyle{definition} 

\numberwithin{equation}{section}

\begin{document}
\title[Stability of composite waves for Navier-Stokes equations]{Stability of composite Wave of Planar Viscous Shock and Rarefaction for 3D Barotropic Navier-Stokes Equations}

\author{Jiajin Shi}
\address{\newline School of Mathematical Sciences, Nanjing Normal University, Nanjing 210023, P. R. China 
\newline and Institute of Applied Mathematics, Academy of Mathematics and Systems Science, Chinese Academy of Sciences, Beijing 100190, P. R. China}
\email{shijiajin@amss.ac.cn}

\author{Yi Wang}
\address{\newline Institute of Applied Mathematics, Academy of Mathematics and Systems Science, Chinese Academy of Sciences, Beijing 100190, P. R. China
\newline and School of Mathematical Sciences, University of Chinese Academy of Sciences, Beijing 100049, P. R. China}
\email{wangyi@amss.ac.cn}

\date{}

\thanks{\textbf{Acknowledgment.} The work of Yi Wang is partially supported by the NSFC grants (Grant No.s 12171459, 12090014, 12288201) and CAS Project for Young Scientists in Basic Research, Grant No. YSBR-031.}

\maketitle
We prove the nonlinear time-asymptotic stability of the composite wave consisting of a planar rarefaction wave and a planar viscous shock for the three-dimensional (3D) compressible barotropic Navier-Stokes equations under generic perturbations, in particular, without zero-mass conditions. 
It is shown that if the composite wave strength and the initial perturbations are suitably small, then 3D Navier-Stokes system admits a unique global-in-time strong solution which time-asymptotically converges to the corresponding composite wave up to a time-dependent shift for planar viscous shock.
Our proof is based on the $a$-contraction method with time-dependent shift and suitable weight function. 

\section{Introduction}
We are concerned with the large-time behavior of the solution to 3D compressible barotropic Navier-Stokes equations
\begin{equation}\label{eq:3D-NS}
  \begin{cases}
    \pd_{t}\rho + \div (\rho u) = 0 \,, \qquad x \in \Omega,\; t \geq 0, \\
    \pd_{t}(\rho u) + \div (\rho u \otimes u) +\nabla p(\rho) = \mu \Delta u + (\mu+\lambda)\nabla\div u \,,
  \end{cases}
\end{equation}
where the unknowns $\rho = \rho(t,x)>0$ and $u=u(t,x)$  represent respectively the
mass density and the velocity of the fluid with $x=(x_1,x_2,x_3) \in \Omega := \R \times \T^2 $ and $ \T:= \R / \mathbb{Z}$. The pressure $p= p(\rho)$ is given by the well-known $\gamma$-law
\begin{equation*}
  p(\rho) = b \rho^\gamma \,, 
\end{equation*}
where $b>0$, $\gamma>1$ are both constants depending on the fluid. Without loss of generality, we will take $b=1$ hereafter. The constants $\mu$ and $\lambda$ are viscosity coefficients satisfying the following physical constraints
\begin{equation*}
  \mu>0, \quad 2\mu+3\lambda \geq 0 \,.
\end{equation*}
The system \eqref{eq:3D-NS} is endowed with initial data with prescribed far-field states in the $x_1$-direction:
\begin{equation}\label{eq:initial data}
  (\rho,u)(t,x)|_{t=0} = (\rho_0(x),u_0(x)) \to (\rho_\pm,u_\pm), \quad {\rm as} \;\; x_1 \rightarrow \pm\infty,
\end{equation}
where $ \rho_\pm>0,\, u_\pm = (u_{1\pm},0,0)^t $ are given planar constant states, and the periodic boundary conditions are imposed on $(x_2,x_3) \in \T^2$ for the solution $(\rho,u)$. 

When the density $\rho$ is away from vacuum ($\rho>0$), it is convenient to introduce the specific volume $v:=1 / \rho$, and then \eqref{eq:3D-NS} can be transformed into the equivalent system with respect to $(v,u)$:
\begin{equation}\label{eq:3D-NS-v,u}
  \begin{cases}
    \rho (\md{v}) = \div u, \\
    \rho (\md{u}) +\nabla p = (2\mu+\lambda) \nabla\div u -\mu \curl\curl u,
  \end{cases}
\end{equation}
where we invoke the identity
\begin{equation*}
  \curl\curl u = \nabla\div u - \Delta u \,.
\end{equation*}

The large-time behavior of solutions to \eqref{eq:3D-NS}-\eqref{eq:initial data} is expected to be determined by the planar Riemann problem of the corresponding 3D compressible Euler system
\begin{equation}\label{eq:3D-Euler}
  \left\{ 
    \begin{aligned}
      &\pd_{t}\rho + \div (\rho u) = 0 \,,  \\
      &\pd_{t}(\rho u) + \div (\rho u \otimes u) +\nabla p(\rho) = 0 \,,  \\
      &(\rho,u)(0,x)=
      \begin{cases}
        (\rho_{-}, u_{-}), \quad x_1<0, \\
        (\rho_{+}, u_{+}), \quad x_1>0,
      \end{cases}
    \end{aligned}
  \right.
\end{equation}
which contains in general two nonlinear planar waves, i.e. planar shock and planar rarefaction wave. 

Such expectation has been well substantiated in one-dimensional case. In 1960, Il'in-Ole\u{\i}nik \cite{Oleinik1960} first validate the stability of single viscous shock and rarefaction wave for scalar quasi-linear equation with convex flux. Later, under zero-mass conditions on initial perturbations, Goodman \cite{Goodman1986} and Matsumura-Nishihara \cite{Matsumura-shock1985} independently proved the stability of single viscous shock to uniformly viscous conservation laws and compressible Navier-Stokes system, respectively. Then, Liu \cite{Liu-1985} and Szepessy-Xin \cite{Szepessy-1993} removed the zero-mass conditions in \cite{Goodman1986} by introducing a constant shift to the shock profile together with extra diffusion waves in the transverse characteristic fields, and the physical viscosity case was tackled by Liu-Zeng \cite{Liu-Zeng-2015}. 
For the stability of rarefaction wave, we refer to Matsumura-Nishihara \cite{Matsumura-rarefaction1986, Matsumura-rarefaction1992} for the isentropic Navier-Stokes system, and Liu-Xin \cite{liu-xin-1997}, Nishihara-Yang-Zhao \cite{Yang-Zhao-2004} for the non-isentropic case. 
As to the stability of composite waves, we refer to Huang-Matsumura \cite{Huang-Matsumura-2009} for two viscous shocks, and Huang-Li-Matsumura \cite{huang-li-Matsumura-2010} for the combination of a viscous contact wave and a rarefaction wave. The case of the superposition of a viscous shock and a rarefaction wave is more challenging \cite{Matsumura2018}, for the incompatibility between the standard anti-derivative method, used to study the stability of viscous shocks, and the direct energy method for rarefactions. Recently, Kang-Vasseur-Wang \cite{Kang-Vasseur-Wang-2023} conquered such difficulty by the method of $a$-contraction with shifts and then \cite{K-V-W-NSF} showed the stability of generic Riemann solutions consisting of a viscous shock, a rarefaction and a viscous contact wave for Navier-Stokes-Fourier system.

Comparatively speaking, the stability of elementary wave patterns for multi-dimensional viscous conservation laws remains less explored. 
For scalar viscous conservation laws, we refer to Goodman \cite{Goodman-1989}, Hoff-Zumbrun \cite{Hoff-Zumbrun-2000, Hoff-Zumbrun-2002} and Kang-Vasseur-Wang \cite{KVW-2019} for the stability of planar viscous shocks, as well as Xin \cite{Xin-1990}, Ito \cite{Ito-1996} and Nishikawa-Nishihara \cite{Nishikawa-2000} for the stability of planar rarefaction waves. 
As to the multi-dimensional Navier-Stokes system, Humpherys-Lyng-Zumbrun \cite{Zumbrun-2017} proved the spectral stability of a planar viscous shock under some spectral conditions in $\R^3$ , while Wang-Wang \cite{Teng-Yi-JEMS} and Kang-Lee \cite{kang-lee-2024} proved the nonlinear stability for single and two planar viscous shocks under generic perturbations in $\R \times \T^2$ respectively.  
In addition, Li-Wang-Wang \cite{Li-teng-yi-2018} proved the stability of planar rarefaction waves for the 3D Navier-Stokes-Fourier system in $\R \times \T^2$. 
Very recently, Meng \cite{meng2025} proved the stability of composite waves for the general $n \times n$ multi-dimensional conservation laws system in $\R \times \T^{d-1}$ with uniformly positive viscosities, in particular, including compressible barotropic Brenner-Navier-Stokes equations. However, this general result can not be directly applicable to the compressible barotropic Navier-Stokes system \eqref{eq:3D-NS} with physical viscosities. Meanwhile, it should be remarked that  the nonlinear stability of a planar viscous shock and its composite waves for the multi-dimensional compressible Navier-Stokes-Fourier system remains largely open even with the aid of $a$-contraction method, due to the temperature effects in the energy equation and their couplings.

The purpose of this paper is to prove the nonlinear time-asymptotic stability of the composite wave consisting of a planar rarefaction wave and a planar viscous shock for 3D barotropic Navier-Stokes system \eqref{eq:3D-NS} under generic perturbations without zero-mass conditions, which can be viewed as a 3D version of \cite{Kang-Vasseur-Wang-2023}. Compared with the tractable Lagrangian structure in one-dimensional case, 3D Navier-Stokes system is expressed in Eulerian coordinates, whose diffusive terms are more difficult to be fully utilized. Motivated by \cite{Teng-Yi-JEMS}, effective velocity is introduced to make full use of the dissipative effect. Also, some underlying physical structures are used to overcome the difficulties arising from the wave propagation along transverse directions and the interactions between the planar rarefaction wave and the planar viscous shock.

The rest of the paper is organized as follows. 
In Section \ref{sec:2}, we first outline the properties of planar rarefaction and viscous shock, then we state our main result. 
In Section \ref{sec:3}, we begin by introducing some useful inequalities, followed by the construction of the weight function and the shift function. The proof of our main result is then provided, based on the existence of local solution in Proposition \ref{prop:local-exist} and uniform-in-time a priori estimates in Proposition \ref{prop:apriori-estimate}. 
In Section \ref{sec:4}, we first reformulate the problem in terms of $(v,h)$, and then validate the uniform-in-time a priori estimates as stated in Proposition \ref{prop:apriori-estimate}.

\section{Preliminaries and main result}\label{sec:2}

In this section, we first exhibit the construction and properties of planar rarefaction and viscous shock. 
Then we state our main result on the nonlinear time-asymptotic stability of composite waves to 3D compressible Navier-Stokes equations \eqref{eq:3D-NS} under generic $H^2$-perturbations without zero-mass conditions.

\subsection{Planar rarefaction wave}
Throughout this paper, we consider the planar wave that depends only on $x_1$ and has null component in the $(x_2,x_3)$ direction. Thus, in order to construct planar rarefaction wave, we only need to consider the 1D Euler system
\begin{equation}\label{eq:1D-Euler-E}
  \begin{cases}
    \rho_t+(\rho u_1)_{x_1} =0,\\
    (\rho u_1)_t+(\rho u_1^2+p(\rho))_{x_1} = 0 ,
  \end{cases}
\end{equation}
which is strictly hyperbolic with both genuinely nonlinear characteristic fields when $\rho>0$.
It can be expressed as the following equivalent form when dealing with rarefaction waves,
\begin{equation*}
  \left( \begin{array}{l} \rho \\ u_1 \end{array} \right)_t + 
  \left(
  \begin{array}{cc}
    u_1 &\quad \rho \\
    p'(\rho)/\rho &\quad u_1
  \end{array} \right)
  \left( \begin{array}{l} \rho \\ u_1 \end{array} \right)_{x_1} = 0.
\end{equation*}
The Jacobi matrix
\begin{equation*}
  \left(
  \begin{array}{cc}
    u_1 &\quad \rho\\
    p'(\rho)/\rho &\quad u_1
  \end{array} \right)
\end{equation*}
has two distinct eigenvalues
\begin{equation*}
  \lambda_1(\rho,u_1)= u_1 -\sqrt{p'(\rho)},\qquad \lambda_2(\rho,u_1)= u_1 +\sqrt{p'(\rho)}
\end{equation*}
with corresponding right eigenvectors
\begin{equation*}
  r_i(\rho,u_1)= \left( 1,\, (-1)^i\frac{\sqrt{p'(\rho)}}{\rho} \right)^t,\qquad i=1,2,
\end{equation*}
such that
\begin{equation*}
  r_i(\rho,u) \sdot  \nabla_{(\rho,u)} \lambda_i(\rho,u)=(-1)^i \,\frac{\rho p''(\rho)+2p'(\rho)}{2\rho\sqrt{p'(\rho)}} \neq 0,\; i=1,2,  \quad {\rm if}~\rho>0.
\end{equation*}
Define the $i$-Riemann invariant $(i=1,2)$ as
\begin{equation*}
  \Sigma_i(\rho,u_1) := u_1 + (-1)^{i+1} \int^\rho \frac{\sqrt{p'(s)}}{s} ds ,
\end{equation*}
which satisfies
\begin{equation*}
  \nabla_{(\rho,u_1)} \Sigma_i(\rho,u_1) \cdot r_i(\rho,u) \equiv 0,\qquad \forall\, \rho>0,u.
\end{equation*}
For a given right state $ (\rho_0,u_{10}) \in \R_+ \times \R $, the 1-rarefaction wave curve $ R_1(\rho_0,u_{10}) $ can be defined as the integral curve of the 1-characteristic field, that is
\begin{equation*}
  R_1(\rho_0,u_{10}):\; \frac{du_1}{d\rho} = -\frac{\sqrt{p'(\rho)}}{\rho} .
\end{equation*}
When $p(\rho)= \rho^\gamma, \gamma>1$, it has the following explicit expression:
\begin{equation*}
  R_1(\rho,u_{1}) = \left\{ (\rho,u_{1}) : u_1-u_{10} = \frac{2\sqrt{\gamma}}{\gamma-1} \left( \rho^{-\frac{\gamma-1}{2}} - \rho_0^{-\frac{\gamma-1}{2}}  \right)   \right\} .
\end{equation*}

Now, we utilize the $i$-Riemann invariants to construct rarefaction waves for \eqref{eq:1D-Euler-E}. For the sake of brevity, we only consider 1-rarefaction wave and 2-rarefaction wave can be obtained similarly. To do this, we first study the Riemann problem for the inviscid Burgers equation
\begin{equation}\label{eq:jump Burgers}
  \begin{cases}
    w_t + ww_{x_1} =0,\quad \,t>0, ~ x_1 \in\R, \\
    w(0,x_1) = w_0^r(x_1) =
    \begin{cases}
      w_-,\; x_1<0,\\
      w_m,\; x_1>0.
    \end{cases}
  \end{cases}
\end{equation}
If $w_- < w_m$, the Riemann problem \eqref{eq:jump Burgers} admits a rarefaction wave solution
$w^r(t,x_1) = w^r(x_1/t)$ which reads as 
\begin{equation}\label{eq:rare-Burgers}
  w^r \left( \frac{x_1}{t} \right) = 
  \begin{cases}
    w_-,     &x_1\leq w_- t,\\
    \displaystyle \frac{x_1}{t},   &w_- t\leq x_1 \leq w_m t,\\
    w_m,     &x_1 \geq w_m t.
  \end{cases}
\end{equation}
For the given states $ (\rho_m,u_{1m}) \in \R_+ \times \R $ and $ (\rho_-,u_{1-}) \in R_1(\rho_m,u_{1m}) $, let
\begin{equation*}
  w_- = \lambda_1(\rho_-,u_{1-}), \quad w_m = \lambda_1(\rho_m,u_{1m}),
\end{equation*}
and $ w^r(t,x_1) := w^r(x_1/t) $ be the solution defined by \eqref{eq:rare-Burgers}. 

Then we can construct the 1-rarefaction wave $ (\rho^r,u_1^r)(x_1/t) $ as follows
\begin{equation}\label{eq:rare-def}
  \begin{cases}
    \lambda_1(\rho^r,u_1^r) (x_1/t) = w^r (x_1/t), \\
    \Sigma_1(\rho^r,u_1^r)(x_1/t) = \Sigma_1(\rho_-,u_{1-}) = \Sigma_1(\rho_m,u_{1m}),
  \end{cases}
\end{equation}
which is the self-similar solution of the 1D Euler system \eqref{eq:1D-Euler-E} with initial data
\begin{equation*}
  (\rho,u_1)|_{t=0} = 
  \begin{cases}
    (\rho_-,u_{1-}), \quad x_1<0, \\
    (\rho_m,u_{1m}), \quad x_1>0.
  \end{cases}
\end{equation*}

\subsection{Approximate rarefaction wave}
Motivated by \cite{Matsumura-rarefaction1986}, we construct the smooth approximate 1-rarefaction wave with the aid of the smooth solution to the Burgers equation
\begin{equation}\label{eq:smooth Burgers}
  \begin{cases}
    w_t + ww_{x_1} =0,\quad \,t>0, ~ x_1 \in\R, \\
    w(0,x_1) = w_0^R(x_1) = \displaystyle \frac{w_m + w_-}{2} + \frac{w_m - w_-}{2} \tanh x_1 .
  \end{cases}
\end{equation}
By the method of characteristics, the solution $w^R(t,x_1)$ to \eqref{eq:smooth Burgers} can be defined as
\begin{equation}
  w^R(t,x_1) = w_0^R(x_0(t, x_1)),\quad\quad x_1=x_0(t, x_1)+w_0^R(x_0(t, x_1))t \,.
\end{equation}
Also, $w^R(t,x_1)$ has the following properties whose proofs can be found in \cite{Matsumura-rarefaction1986} or \cite{Matsumura-rarefaction1992}.

\begin{lemma}\label{lem:w^R property}
  Let $ w_m>w_- $ and $ \tilde{w}:= w_m-w_- $, then the Cauchy problem \eqref{eq:smooth Burgers} admits a unique global smooth solution $ w^R(t, x_1) $ satisfying

  (1)~ $w_- < w^R(t, x_1) < w_m$ and $ w^R_{x_1} >0 $, for $x_1 \in \R$ and $t\geq 0$.
   
  (2)~ For any $t > 0$ and $p \in [1,+\infty]$, there exists a positive constant $C=C_p$  such that
    \begin{align*}
      &\|w^R_{x_1}(t, \cdot)\|_{L^p} \leq C\min \left\{ |\tilde{w}|, |\tilde{w}|^{1/p}t^{-1+1/p} \right\} , \\
      &\|w^R_{x_1x_1}(t, \cdot)\|_{L^p} \leq C\min \left\{ |\tilde{w}|,t^{-1} \right\} ,\\
      &\|w^R_{x_1x_1x_1}(t, \cdot)\|_{L^p} \leq C\min \left\{ |\tilde{w}|,t^{-1} \right\} .
    \end{align*}

  (3)~ For $x_1 \geq w_m t,\; t\geq 0 $, it holds that
    \begin{align*}
      &|w^R(t, x_1)-w_m| \leq \tilde{w}e^{-2|x_1-w_mt|}, \\
      &|(w^R_{x_1}, w^R_{x_1x_1})(t, x_1)|\leq C\tilde{w} e^{-2|x_1-w_mt|}.
    \end{align*}

  (4)~ For $x_1 \leq w_-t,\; t\geq 0 $, it holds that
    \begin{align*}
      &|w^R(t, x_1)-w_-| \leq \tilde{w}e^{-2|x_1-w_-t|}, \\
      &|(w^R_{x_1}, w^R_{x_1x_1})(t, x_1)|\leq C\tilde{w} e^{-2|x_1-w_-t|}.
    \end{align*}

  (5)~ $ \displaystyle \lim_{t \to +\infty} \sup_{x_1 \in \R} \left| w^R(t, x_1) - w^r(\frac{x_1}{t}) \right| = 0 $.
\end{lemma}

Correspondingly, the approximate 1-rarefaction wave $(\rho^R,u_1^R)(t,x_1)$ of the 1-rarefaction wave $(\rho^r,u_1^r)(x_1/t)$ can be defined by
\begin{equation}\label{eq:appro-rare-def}
  \begin{cases}
    w_- = \lambda_1(\rho_-,u_{1-}),\; w_m = \lambda_1(\rho_m,u_{1m}),\\
    \lambda_1(\rho^R,u^R_1)(t,x_1) = w^R(1+t, x_1), \\
    \Sigma_1(\rho^R,u^R_1)(t,x_1) = \Sigma_1(\rho_-, u_{1-}) =\Sigma_1 (\rho_m,u_{1m}),
  \end{cases}
\end{equation}
where $ w^R(t,x_1) $ is the solution of the Burgers equation \eqref{eq:smooth Burgers}. By direct computations, it can be verified that $(\rho^R,u_1^R)$ satisfies the 1D Euler system
\begin{equation*}
  \begin{cases}
    \rho_t+(\rho u_1)_{x_1} =0,\\
    (\rho u_1)_t+(\rho u_1^2+p(\rho))_{x_1} = 0 .
  \end{cases}
\end{equation*}
Also, the following lemma can be derived from Lemma \ref{lem:w^R property} (cf. \cite{Matsumura-rarefaction1986}).

\begin{lemma}\label{lem:appro-rare-property} 
  For any $ (\rho_m,u_{1m}) \in \R_+ \times \R$ and $ (\rho_-,u_{1-}) \in R_1(\rho_m,u_{1m}) $, let $\lambda_{1-} := \lambda_1(\rho_-,u_1-)$, $\lambda_{1m} := \lambda_1(\rho_m,u_{1m})$ and $\delta_R := |v_m - v_-| \sim |u_{1m}-u_-|$ denoting the rarefaction wave strength. 
  Then the approximate 1-rarefaction wave $(\rho^R, u_1^R)(t,x)$ defined by \eqref{eq:appro-rare-def} satisfies the following properties, where $ v^R := 1/\rho^R $,
  
  (1)~ For any $x_1 \in \R$ and any $t \geq 0 $, it holds $ u^R_{1,x_1}>0 $, $ v^R_{x_1}>0 $ and $ u^R_{1,x_1} \sim v^R_{x_1} $.
  
  (2)~ For any $t \geq 0$ and any $p \in [1,+\infty]$, it holds
  \begin{align*}
  &\|( v^R_{x_1}, u^R_{1,x_1})\|_{L^p} \leq C \min\{\delta_R, \delta_R^{1/p}(1+t)^{-1+1/p}\}, \\
  &\|( v^R_{x_1x_1}, u^R_{1,x_1x_1})\|_{L^p} \leq C \min\{\delta_R, (1+t)^{-1}\}, \\
  &|u^R_{1,x_1x_1}| \leq C |u^R_{1,x_1}|,\, |u^R_{1,x_1x_1x_1}| \leq C |u^R_{1,x_1x_1}|, \; \forall x\in\R.
  \end{align*}
  
  (3)~ For $x_1 \geq \lambda_{1m}(1+t),\, t\geq 0 $, it holds
  \begin{align*}
    &|(v^R, u_1^R)(t,x_1)-(v_m,u_m)| \leq C\delta_R \, e^{-2|x_1-\lambda_{1m}(1+t)|}, \\
    &|(v^R_{x_1}, u^R_{1,x_1})(t,x_1)|\leq C\delta_R\, e^{-2|x_1-\lambda_{1m}(1+t)|}.
  \end{align*}
  
  (4)~ For $x_1 \leq \lambda_{1-}(1+t),\, t \geq 0$, it holds 
  \begin{align*}
    &|( v^R, u_1^R)(t,x_1)-(\rho_-,u_-)| \leq C\delta_R \, e^{-2|x_1-\lambda_{1-}(1+t)|}, \\
    &|( v^R_{x_1}, u^R_{1,x_1})(t,x_1)|\leq C\delta_R\, e^{-2|x_1-\lambda_{1-}(1+t)|}.
  \end{align*}
  
  (5)~ $ \displaystyle  \lim_{t \to +\infty} \sup_{x \in \R} \left| (v^R, u_1^R)(t,x_1) - (v^r, u_1^r)\left( \frac{x_1}{t} \right) \right| = 0$.
\end{lemma}

\begin{remark}
  We wish to emphasize the distinction between the superscript $\centerdot^r$  and the superscript $\centerdot^R$. Throughout this paper, the former denotes the standard rarefaction waves, while the latter signifies the approximated rarefaction waves.
\end{remark}

\subsection{Planar viscous shock}
Like before, we consider the 1D barotropic Navier-Stokes system
\begin{equation}\label{eq:1D-NS-E}
  \begin{cases}
    \pd_{t}\rho + \pd_{1}(\rho u_1) = 0 , \\
    \pd_{t}(\rho u_1) + \pd_{1}(\rho u_1^2 + p(\rho)) = (2\mu+\lambda) \pd_{11}u_1 \,,
  \end{cases}
\end{equation}
with initial data
\begin{equation}\label{eq:1D-NS-E-data}
  (\rho,u_1)(0,x_1) = (\rho_0,u_{10})(x_1) = 
  \begin{cases}
    (\rho_m,u_{1m}), \quad &x_1<0, \\
    (\rho_+,u_{1+}), \quad &x_1>0.
  \end{cases}
\end{equation}
We hope to seek for travelling wave solutions of the form $(\rho,u_1)(t,x_1)=(\rho^s,u^s_1)(\xi)$, where $ \xi=x_1-\s t$ and $\s$ is a constant to be determined. Inserting the Ansatz $(\rho,u_1)(t,x_1)=(\rho^s,u^s_1)(\xi)$ into \eqref{eq:1D-NS-E} leads to
\begin{equation}\label{eq:rho,u-ODE}
  \begin{cases}
    -\s(\rho^s)'+(\rho^s u^s_1)'=0, \qquad\qquad\qquad\qquad\qquad\quad ':=\frac{d}{d\xi}, \\
    -\s(\rho^s u_1^s)'+(\rho^s (u_1^s)^2)'+p(\rho^s)'=(2\mu+\lambda)(u^s_1)'', 
  \end{cases}
\end{equation}
which is an ODE for $(\rho^s,u^s_1)(\xi)$. The initial data \eqref{eq:1D-NS-E-data} are then converted to
\begin{equation}\label{eq:rho,u-ODE-data}
  (\rho^s, u^s_1)(-\infty)=(\rho_m, u_{1m}),\quad (\rho^s, u^s_1)(+\infty)=(\rho_+, u_{1+}).
\end{equation}
Integrating \eqref{eq:rho,u-ODE} with respect to $ \xi $ over $ \R $ results in the following Rankine-Hugoniot condition
\begin{equation}\label{eq:sigma's eqn}
  \begin{cases}
    -\s(\rho_+ - \rho_m)+(\rho_+ u_{1+}-\rho_m u_{1m})=0,\\
    -\s(\rho_+ u_{1+} -\rho_m u_{1m})+(\rho_+ u_{1+}^2-\rho_m u_{1m}^2)+({p}(\rho_+)-{p}(\rho_m))=0, \\
  \end{cases}
\end{equation}
which is necessary for \eqref{eq:rho,u-ODE}-\eqref{eq:rho,u-ODE-data} to be solvable. 
Furthermore, we say $ (\rho_m,u_{1m}) \in S_2(\rho_+,u_{1+}) $ if it also satisfies the Lax entropy conditions
\begin{equation}\label{eq:Lax-E condition}
  \begin{cases}
    \s > \lambda_1(\rho_m,u_{1m}), \\
    \lambda_2(\rho_+,u_{1+}) < \sigma < \lambda_2(\rho_m,u_{1m}),
  \end{cases}
\end{equation}
where $ \lambda_1(\rho,u_1) = u_1 - \sqrt{p'(\rho)}$, $\lambda_2(\rho,u_1) = u_1 + \sqrt{p'(\rho)} $.

When away from vacuum, let $ v^s := 1/\rho^s $, then the ODE \eqref{eq:rho,u-ODE} can be converted to an ODE for $ (v^s,u^s_1)(\xi) $:
\begin{equation}\label{eq:2.8}
  \begin{cases}
    \rho_s \big(-\s(v^s)' + u^s_1(v^s)' \big) = (u^s_1)',  \\
    \rho_s \big(-\s(u^s_1)' + u^s_1(u^s_1)' \big) + p(v^s)' = (2\mu+\lambda)(u^s_1)'', 
  \end{cases}
\end{equation}
where $ p(v^s) = (v^s)^{-\gamma} $. Integrating \eqref{eq:rho,u-ODE}$ _1 $ over $ (-\infty,\xi) $ leads to
\begin{equation}\label{eq:sigma*-def}
  -\s\rho^s +\rho^s u^s_1 = -\s\rho_- + \rho_- u_{1-} =: \s_* \,.
\end{equation}
By virtue of the $ \s_* $ defined above, the system \eqref{eq:2.8} and the far-field conditions \eqref{eq:rho,u-ODE-data} can be rewritten as 
\begin{equation}\label{eq:v,u-ODE}
  \begin{cases}
    -\s_*(v^s)' = (u^s_1)', \\
    -\s_*(u_1^s)' + p(v^s)' = (2\mu+\lambda)(u^s_1)'', 
  \end{cases}
\end{equation}
and
\begin{equation}\label{eq:v,u-ODE-data}
  (v^s, u^s_1)(-\infty) = (v_m, u_{1m}),\quad (v^s, u^s_1)(+\infty) = (v_+, u_{1+}),
\end{equation}
where $ v_m=1/\rho_m,\,v_+=1/\rho_+ $. 
In addition, integrating \eqref{eq:v,u-ODE} over $ \R $ and combining \eqref{eq:sigma*-def}, we have
\begin{equation*}
  \begin{cases}
    -\s_* (v_+ - v_m) = u_{1+}-u_{1m}, \\
    -\s_* (u_{1+}-u_{1m}) + p(v_+)-p(v_m) = 0, 
  \end{cases}
\end{equation*}
then we can obtain another useful expression for $ \s_* $:
\begin{equation}\label{eq:sigma*>0}
  \s_* = \sqrt{-\frac{p(v_+)-p(v_m)}{v_+ - v_m}} > 0.
\end{equation}

Then, we summarize the properties of the 2-viscous shock in the following lemma, whose proof can be found in \cite{Matsumura-shock1985}.

\begin{lemma}\label{lem:viscous shock}
  For any given state $(v_+,u_{1+})$ and $(v_m,u_{1m}) \in S_2(v_+,u_{1+})$, let
  \begin{equation*}
    \delta_S := |p(v_m)-p(v_+)| \sim |v_m-v_+|\sim |u_{1m}-u_{1+}|,
  \end{equation*}
  which denotes the shock wave strength.
  Then, up to a constant shift, the ODE \eqref{eq:v,u-ODE}-\eqref{eq:v,u-ODE-data} admits a unique solution $ (v^s,u^s_1)(\xi) $ satisfying
  \begin{equation*}
    (v^s)'>0, \quad (u_1^s)' = -\s_*(v^s)' <0.
  \end{equation*}
  In addition, there exists a constant $ C>0 $ such that the following holds
  \begin{align*}
    &\left| v^s(\xi)-v_m \right|  \leq C\delta_S e^{-C\delta_S|\xi|}, \quad \forall \xi<0,\\
    &\left| v^s(\xi)-v_+ \right|  \leq C\delta_S e^{-C\delta_S|\xi|}, \quad \forall \xi>0,\\
    &\left| \big( (v^s)',(u_1^s)' \big) \right|  \leq C\delta_S^2 e^{-C\delta_S|\xi|}, \quad \forall \xi \in \R, \\
    &\left| \big( (v^s)'',(u_1^s)'' \big) \right| \leq C\delta_S \left| \big( (v^s)',(u_1^s)' \big) \right|, \quad \forall \xi \in \R, \\
    &\left| \big( (v^s)''',(u_1^s)''' \big) \right|  \leq C\delta_S^2 \left| \big( (v^s)',(u_1^s)' \big) \right|, \quad \forall \xi \in \R.
  \end{align*}  
\end{lemma}

\subsection{Notations}\label{sec:notation}
For the convenience of the readers, we list some notational conventions that taken throughout this paper.

(1)~ Let $\Omega$ denote the domain $\R \times \T^2$. For $x=(x_1,x_2,x_3) \in \Omega$, we define 
\begin{equation*}
  x' := (x_2,x_3),\quad dx' := dx_2dx_3 .
\end{equation*}

(2)~ Let $1 \leq r \leq +\infty$, we denote $L^r := L^r(\Omega)$ and 
\begin{equation*}
  \norm{f} := \norm{f}_{L^2} = \left( \int |f|^2 \,dx \right)^{\frac{1}{2}} .
\end{equation*}
For $m \in \mathbb{N}$, $H^m := H^m(\Omega)$ denotes the usual Sobolev space based on $L^{2}$  with the norm 
\begin{align*}
  \norm{f}_{H^m} := \left( \sum_{k=0}^m \norm{\nabla^k f}^2 \right)^{\frac{1}{2}},\quad \norm{(f,g)}_{H^m} := \norm{f}_{H^m} +\norm{g}_{H^m}.
\end{align*}

(3)~ Let $(v^s,u_1^s)$ be the 2-viscous shock defined in Lemma \ref{lem:viscous shock} and $\X(t)$ be any time-dependent shift, we define the shifted shock with uppercase superscript to distinguish
\begin{equation}\label{eq:v^S-def}
  (v^S,u^S_1)(t,x) := (v^s,u^s_1)\big( x_1- \s t -\X(t) \big),
\end{equation}
where $\s$ is the constant defined by \eqref{eq:sigma's eqn}.

(4)~ $(v^r,u_1^r)$ and $(v^R,u^R_1)$ and given by \eqref{eq:rare-def} and \eqref{eq:appro-rare-def}, respectively. Let $u^R:=(u^R_1,0,0)^t$, $u^S:=(u^S_1,0,0)^t$ and
\begin{equation}
  \begin{cases}
    \bar{v}(t,x) := (v^S + v^R - v_m)(t,x), \\
    \bar{u}_1(t,x) := (u^S_1 + u^R_1 - u_{1m})(t,x), \\
    \bar{u}(t,x) := \big(\bar{u}_1(t,x),0,0\big)^t, 
  \end{cases}
\end{equation}
all of which depend only on $(t,x_1)$.  In addition, we define the error terms as follows
\begin{equation}\label{eq:phi,psi-def}
  \begin{cases}
    \phi := v-\bar{v},  \\
    \psi := u-\bar{u} = (u_1-\bar{u}_1,0,0)^t = (\psi_1,0,0)^t.
  \end{cases}
\end{equation}

(5)~ For the pressure function $p(v)=v^{-\gamma}$, we often take the following abbreviations 
  \begin{equation*}
    p^S:=p(v^S),\quad p^R:=p(v^R),\quad \bar{p}:=p(\bar{v}).
  \end{equation*}

(6)~ $\delta_R$, $\delta_S$ and $\delta$ denote the strength of the 1-rarefaction, 2-viscous shock and the composite wave, respectively, which are defined as
\begin{equation}\label{eq:delta_R,S-def}
  \begin{cases}
    \delta_R := |v_m - v_-| \sim |u_{1m}-u_-|, \\
    \delta_S := |p(v_m)-p(v_+)| \sim |v_m-v_+|\sim |u_{1m}-u_{1+}|, \\
    \delta := \delta_R  +\delta_S .
  \end{cases}
\end{equation}

\subsection{Main result}

Now we can state our main result as follows. Roughly speaking, as long as the strength of the composite wave is suitably small, then it is nonlinearly stable under generic $H^2$-perturbations without zero-mass conditions.

\begin{theorem}\label{thm:main-ch2}
  For any given state $(v_+,u_{1+}) \in \R_+ \times \R $ there exist constants $\delta_0, \eps_0>0$ such that for each $(v_m,u_{1m})\in S_2(v_+,u_{1+})$, $(v_-,u_{1-})\in R_1(v_m,u_{1m})$ with
  $$ \delta := |v_+-v_m|+|v_m-v_-| \leq \delta_0 $$ 
  and each initial data with
  \begin{equation}\label{eq:initially small}
    \norm{(v_0-\bar{v}(0,\cdot), u_0-\bar{u}(0,\cdot))}_{H^2(\R \times \T^2)} \leq \varepsilon_0 ,
  \end{equation}
  the 3D barotropic Navier-Stokes system \eqref{eq:3D-NS-v,u} admits a unique global-in-time solution $ (v,u)(t,x)$ with $v := \frac{1}{\rho}$. 
  In addition, there exists an absolutely continuous shift $\X(t)$ such that
  \begin{equation}\label{eq:regulartiy}
    \begin{aligned}
      &v(t,x)-\bar{v}(t,x) \in C \big( [0,+\infty); H^2(\R \times \T^2)\big), \\
      &u(t,x)-\bar{u}(t,x) \in C \big( [0,+\infty); H^2(\R \times \T^2)\big), \\
      &\nabla \big( v(t,x)-\bar{v}(t,x) \big) \in L^2\big( 0,+\infty; H^1(\R \times \T^2) \big), \\
      &\nabla \big( u(t,x)-\bar{u}(t,x) \big) \in L^2\big( 0,+\infty; H^2(\R \times \T^2) \big).
    \end{aligned}
  \end{equation}
  Here, 
  \begin{equation}\label{eq:vbar ubar-def}
    \begin{aligned}
      \bar{v}(t,x) &:= v^R(t,x_1) + v^s(x_1-\s t-\X(t)) - v_m, \\
      \bar{u}(t,x) &:= u^R(t,x_1) + u^s(x_1-\s t-\X(t)) - u_m,
    \end{aligned}
  \end{equation}
  where $ (v^R,u^R)(t,x_1) $ and $ (v^s,u^s)(x_1-\s t) $ represent the planar approximate 1-rarefaction and the planar 2-viscous shock, i.e.
  \begin{align*}
    (v^R,u^R)(t,x_1) &:= \big( v^R,(u^R_1,0,0)^t \big) (t,x_1), \\
    (v^s,u^s)(x_1-\s t) &:= \big( v^s,(u^s_1,0,0)^t \big)(x_1-\s t),\\
    u_m &:= (u_{1m},0,0)^t .
  \end{align*}

  Furthermore, we have the following time-asymptotically stability. As $t \to +\infty$, 
  \begin{equation}\label{eq:asym-stable}
    \begin{aligned}
      \sup_{x\in\R} \Big| (v,u)(t,x) - \Big( &v^r(\frac{x_1}{t}) + v^s(x_1-\s t-\X(t)) - v_m , \\
      &\; u^r(\frac{x_1}{t}) +  u^s(x_1-\s t-\X(t)) - u_m \Big) \Big|  \to 0,
    \end{aligned}
  \end{equation}
  and
  \begin{equation}\label{eq:Xdot-->0}
    \lim_{t \to +\infty} \left| \dot{\X}(t) \right| = 0, 
  \end{equation}
  where
  \begin{equation*}
    (v^r,u^r)\left( \frac{x_1}{t} \right)  := \big( v^r,(u^r_1,0,0)^t \big)\left( \frac{x_1}{t} \right),
  \end{equation*}
  and $(v^r,u^r_1)\left( \frac{x_1}{t} \right)$ is the 1-rarefaction defined by \eqref{eq:rare-def} with end states $(v_-,u_{1-})$ and $(v_m,u_{1m})$.
\end{theorem}

\begin{remark}
  Theorem \ref{thm:main-ch2} states that if the two far-field states $(\rho_\pm,u_\pm)$ in \eqref{eq:initial data} are connected by the superposition of a rarefaction wave and a shock, then the solution to the three-dimensional compressible Navier-Stokes system \eqref{eq:3D-NS} tends to the superposition wave of the planar inviscid rarefaction wave and the planar viscous shock with shift $\X(t)$ in long time.
\end{remark}

\begin{remark}
  Compared with the 1D result in \cite{Kang-Vasseur-Wang-2023}, the 3D Navier-Stokes system is expressed in Eulerian coordinates, whose diffusive terms are more difficult to be fully utilized. Therefore, motivated by \cite{Teng-Yi-JEMS}, we need to introduce an effective velocity $h:= u - (2\mu+\lambda)\nabla v $ 
  to make full use of the dissipative effect. In addition, some underlying physical structures are used to overcome the difficulties arising from the wave propagation along transverse directions and the interactions between the planar rarefaction wave and the planar viscous shock.
\end{remark}

\begin{remark}
  The shift function is proved to satisfy the time-asymptotic behavior \eqref{eq:Xdot-->0}, which implies that
  \begin{equation*}
     \lim_{t \to +\infty} \frac{\X(t)}{t} = 0.
  \end{equation*}
  Hence, the shift function grows at most sub-linearly with respect to the time $t$, which time-asymptotically preserves the traveling wave profile of the shifted viscous shock. 
\end{remark}

\section{Proof of the main result}\label{sec:3}

The main strategy of our proof is the use of of the method of $a$-contraction with shifts \cite{KV-2017, Leger-2011}, which is based on the relative entropy introduced by DiPerna \cite{DiPerna1979} and Dafermos \cite{Dafermos1996}. 
To achieve this, we first give some technical lemmas and then construct the weight function and the time-dependent shift. 
Finally, we present the proof of Theorem \ref{thm:main-ch2} based on a local existence result and uniform-in-time a priori estimates.

\subsection{Some inequalities}

Firstly, we introduce a 3D weighted sharp Poincar\'e type inequality that given in \cite{Teng-Yi-JEMS}, which can be viewed as a 3D version of the 1D case in \cite{weighted-poincare}.

\begin{lemma}\label{lem:Poincare}
  For every $f:[0,1]\times\T^2 \to \R$ satisfying
  \begin{equation*}
    \int_{\T^2}\int_0^1 y_1(1-y_1)|\pd_{y_1}f|^2 + \frac{|\nabla_{y'}f|^2}{y_1(1-y_1)} \,dy_1dy' < \infty,
  \end{equation*}
  it holds
  \begin{equation}\label{eq:weighted Poincare}
    \begin{aligned}
      &\int_{\T^2}\int_0^1 |f-\bar{f}|^2 \,dy_1dy' \\
      \leq &\; \frac{1}{2} \int_{\T^2}\int_0^1 y_1(1-y_1) |\pd_{y_1}f|^2 dy_1dy'  + \frac{1}{16\pi^2} \int_{\T^2}\int_0^1 \frac{|\nabla_{y'}f|^2}{y_1(1-y_1)} dy_1dy' \,,
    \end{aligned}
  \end{equation}
  where $\displaystyle \bar{f}:=\int_{\T^2}\int_0^1 f dy_1dy'$, $y':=(y_2,y_3)$.
\end{lemma}

The following 3D Gagliardo-Nirenberg type inequality in the domain $\Omega = \R \times \T^2$ is also useful, whose proof can be found in \cite{3D-Gargliardo2022}.
\begin{lemma}
  Let $\Omega := \R \times \T^2$. There exists constant $C>0$, such that for each $g \in H^2(\Omega)$, it holds
  \begin{equation}\label{eq:G-N ineq}
    \|g\|_{L^{\infty}(\Omega)} \leq \sqrt{2} \|g\|^{\frac{1}{2}}_{L^2(\Omega)} \|\pd_{1}g\|^{\frac{1}{2}}_{L^2(\Omega)}
    + C\|\nabla g\|^{\frac{1}{2}}_{L^2(\Omega)} \|\nabla^2 g\|^{\frac{1}{2}}_{L^2(\Omega)}.
  \end{equation}
\end{lemma}

\subsection{Relative quantities}

Then we present some estimates concerning the relative quantities. 
For any function $F$ defined on some interval $I$ of $\R$, we define the associated relative quantity with $v,u \in I$, 
\begin{equation*}
  F(v|w) := F(v)-F(w)-F'(w)(v-w).
\end{equation*}
We will state some useful inequalities about the relative quantities associated to the pressure $p(v)=v^{-\gamma}$ and the internal energy $Q(v):=\frac{v^{1-\gamma}}{\gamma-1}$. The proofs, which are based on Taylor expansions, can be accessed in \cite{weighted-poincare}. 

\begin{lemma}\label{lem:relative quntities}
  Let $p(v)=v^{-\gamma}$ and $Q(v):=\frac{v^{1-\gamma}}{\gamma-1}$ with $\gamma>1$. For given constants $v_->0$, there exist constants $C,\delta_*>0$ such that the following holds true.
  
  (1)~ For any $v,w \in \R_+$ such that $0<w<2v_-$ and $0<v\leq 3v_-$, it holds
  \begin{equation}\label{eq:relative-1}
    |v-w|^2 \leq C Q(v|w),\qquad |v-w|^2 \leq C p(v|w). 
  \end{equation}

  (2)~ For any $\displaystyle v,w>\frac{v_-}{2}$, it holds
  \begin{equation}
    |p(v)-p(w)|\leq C|v-w|.
  \end{equation}

  (3)~ For any $0<\delta<\delta_*$ and $(v,w) \in \R^2_+$ satisfying $|p(v)-p(w)|<\delta$ and $|p(w)-p(v_-)|<\delta$, the following holds true:
  \begin{equation}\label{eq:relative-2}
    p(v|w) \leq \Big(\frac{\gamma+1}{2\gamma}\frac{1}{p(w)}+ C\delta\Big) |p(v)-p(w)|^2,  
  \end{equation}
  \begin{equation}\label{eq:relative-3}
    Q(v|w) \geq \frac{|p(v)-p(w)|^2}{2\gamma p^{1+\frac{1}{\gamma}}(w)}
    - \frac{1+\gamma}{3\gamma^2}\frac{(p(v)-p(w))^3}{p^{2+\frac{1}{\gamma}}(w)}, 
  \end{equation}
  \begin{equation}\label{eq:relative-4}
    Q(v|w) \leq \Big(\frac{1}{2\gamma p^{1+\frac{1}{\gamma}}(w)}+C\delta\Big) |p(v)-p(w)|^2.  
  \end{equation}  
\end{lemma}

At last, we list an estimate about the inverse function of $p(v)=v^{-\gamma}$, whose proof can be found in \cite{weighted-poincare}. 
\begin{lemma}\label{lem:inverse-pressure}
  For given constant $v_- > 0$, there exist $\delta_0>0, C>0$ such that for any $v_+ >0$ with $0< \delta:= p(v_-)-p(v_+) \leq \delta_0$ and $v_- \leq v \leq v_+$, it holds 
  \begin{equation}
    \left| \frac{v-v_-}{p(v)-p(v_-)} + \frac{v-v_+}{p(v_+)-p(v)}+ \frac{1}{2}\frac{p''(v_-)}{p'(v_-)^2}(v_--v_+) \right| \leq C\delta^2.
  \end{equation}
\end{lemma}

\subsection{Construction of the weight function and the time-dependent shift}\label{sec:a-X-def}

We define the weight function  $a(t,x)$ by
\begin{equation}\label{eq:a-weight-def}
 \begin{aligned}
    a(t,x) &:= 1 + \frac{\nu}{\delta_S} \big(p(v_m)-p(v^S)\big) \\ &= 1 + \frac{\nu}{\delta_S} \big( p(v_-)-p(v^s(x_1-\s t-\X(t))) \big),
 \end{aligned}
\end{equation}
where $\nu$ is a constant to be determined with
\begin{equation}
  \delta \ll \nu \lesssim \sqrt{\delta_S}.
\end{equation}
For the sake of clarity, we just choose
\begin{equation}\label{eq:nu-def}
  \nu = \sqrt{\delta_S}.
\end{equation}
Notice that $a(t,x)$ depends only on $(t,x_1)$ and satisfies
\begin{equation}
  1< a(t,x) < 1+\nu \,.
\end{equation}
Also, we have
\begin{equation}
  \pd_{1}a = -\frac{\nu}{\delta_S} p'(v^S) \pd_{1}v^S > 0,
\end{equation}
and then
\begin{equation}
  \pd_{1}a \sim \frac{\nu}{\delta_S} \pd_{1}v^S = \frac{1}{\nu} \pd_{1}v^S.
\end{equation}

With the above weight function $a(t,x)$ given, we can define the shift $\X(t)$ as the solution of the following Cauchy problem:
\begin{equation}\label{eq:X(t)-def} 
  \begin{cases}
    \displaystyle \dot{\X}(t) = -\frac{M}{\delta_S } \left[ \int_{\Omega} \frac{a\rho}{\s_*}\sdot \pd_{1}h^S_1\sdot \big(p(v)-p(\bar{v}) \big) - a\rho\sdot \pd_{1}[p(v^S)]\sdot (v-\bar{v}) \,dx \right], \\
    \X(0) = 0,
  \end{cases}
\end{equation}
where $h^S_1:= u^S_1 - (2\mu+\lambda)\pd_{1}v^S$, and 
\begin{equation}\label{eq:M-def}
  M:= \frac{5(\gamma+1)}{8\gamma}\sdot \frac{\s_m^3 v_m^2}{p(v_m)},\quad \s_m := \sqrt{-p'(v_m)}.
\end{equation}
The shift $\X(t)$ is well-defined for the well-posedness of the ODE \eqref{eq:X(t)-def}, which is guaranteed by the following lemma (cf. \cite{Kang-ODE-2020}).

\begin{lemma}\label{lem:ODE-existence}
  Let $p>1$ and $T>0$. Suppose that a function $F:[0,T] \times \R \rightarrow \R$ satisfies 
  \begin{equation*}
    \sup_{x \in \R} |F(t,x)| \leq f(t)\; \text{~and~}
    \sup_{x,y \in \R, x\neq y } \left| \frac{F(t,x)-F(t,y)}{x-y} \right| \leq g(t) \text{~for~} t \in [0,T]
  \end{equation*}
  for some functions $f \in L^1(0,T)$ and $ g\in L^p(0,T)$. Then for any $x_0 \in \R$, there exists a unique absolutely continuous function $\X:[0,T] \rightarrow \R$ satisfying
  \begin{equation}\label{eq:X-ODE-exist}
    \begin{cases}
      \dot{\X}(t) = F(t,\X(t)) \quad \mbox{for \textit{a.e.} } t \in [0,T],\\
      X(0)= x_0 \,.
    \end{cases}
  \end{equation}
\end{lemma}

Let $F(t,\X)$ be the right-hand side of the ODE \eqref{eq:X(t)-def}. Due to the facts that
\begin{equation*}
  \norm{a}_{C^1} \leq 2, \quad \norm{v^S}_{C^2} \leq \max\{ 1,v_+ \}, \quad \norm{\pd_{1}v^S}_{L^1} \lesssim \delta_S ,
\end{equation*}
we find that for some constant $C>0$, 
\begin{equation}\label{eq:3.12}
  \sup_{\X \in \R} |F(t,\X)| \leq \frac{C}{\delta_S} \left( \norm{v}_{L^\infty} + \norm{v^S}_{L^\infty} \right) \int_{\Omega} \pd_{1}v^S \,dx \leq C,
\end{equation}
and
\begin{equation*}
  \sup_{\X \in \R} |\pd_{\X} F(t,\X)| \leq \frac{C}{\delta_S} \norm{a}_{C^1} \left( \norm{v}_{L^\infty} + \norm{v^S}_{L^\infty} \right) \int_{\Omega} \pd_{1}v^S \,dx \leq C .
\end{equation*}
Thus by Lemma \ref{lem:ODE-existence}, the shift $\X(t)$ is well-defined and absolutely continuous. In addition, we have $|\dot{\X}(t)| \leq C$ from \eqref{eq:3.12}, and consequently
\begin{equation}\label{eq:X < Ct}
  \left| \X(t) \right|  \leq Ct, \qquad \forall \, t \in [0,T].
\end{equation}

\subsection{Proof of Theorem \ref{thm:main-ch2}}

Roughly speaking, we prove Theorem \ref{thm:main-ch2} by continuation arguments in virtue of a local existence result together with a priori estimates. 
For that, we first state the following local existence result, which can be proved by standard methods (see \cite{Nash1962} or \cite{Solonnikov1976}).

\begin{proposition}[Local existence]\label{prop:local-exist}
  Let $(\tilde{v},\tilde{u})(t,x)$ be the composite wave without shift, i.e.
  \begin{align*}
    \tilde{v}(t,x) &:= v^R(t,x_1) + v^s(x_1-\s t) - v_m, \\
    \tilde{u}(t,x) &:= u^R(t,x_1) + u^s(x_1-\s t) - u_m.
  \end{align*}
  For any $\Xi>0$, there exists $T_0>0$ such that for each initial data $(v_0,u_0)$ with
  \begin{equation*}
    \norm{\big( v_0-\tilde{v}(0,\cdot) \,,\, u_0-\tilde{u}(0,\cdot) \big)}_{H^2(\R\times\T^2)} \leq \Xi \,,
  \end{equation*}
  the 3D Navier-Stokes system \eqref{eq:3D-NS-v,u} admits a unique solution $(v,u)$ on $[0,T_0]$ satisfying
  \begin{align*}
    &v-\tilde{v} \in C([0,T_0];H^2(\R \times \T^2)), \quad \nabla(v-\tilde{v}) \in L^2(0,T_0;H^1(\R \times \T^2)), \\
    &u-\tilde{u} \in C([0,T_0];H^2(\R \times \T^2)), \quad \nabla(u-\tilde{u}) \in L^2(0,T_0;H^2(\R \times \T^2)),
  \end{align*}
  and for $t \in [0,T_0]$, it holds that
  \begin{equation}\label{eq:local-est}
    \begin{aligned}
      &\sup_{\tau \in [0,t]} \norm{(v-\tilde{v}, u-\tilde{u})(\tau)}_{H^2}^2 + \int_0^t \norm{\nabla(v-\tilde{v})}_{H^1}^2 + \norm{\nabla(u-\tilde{u})}_{H^2}^2 d\tau \\
      & \leq 4 \norm{\big( v_0-\tilde{v}(0,\cdot) \,,\, u_0-\tilde{u}(0,\cdot) \big)}_{H^2}^2.
    \end{aligned}
  \end{equation}
\end{proposition}

The pivotal part in the proof of Theorem \ref{thm:main-ch2} is the following a priori estimates. 

\begin{proposition}[A priori estimates]\label{prop:apriori-estimate}
  Let $(v,u)$ be the solution of \eqref{eq:3D-NS-v,u} on $[0,T]$ for some $T>0$, $(\bar{v},\bar{u})$ be the composite wave given by \eqref{eq:vbar ubar-def}. Then there exist positive constants $ \delta_0 \leq 1$, $\chi_0 \leq 1$, $C_0$ independent of $T$, such that if the composite wave strength $\delta:=\delta_R +\delta_S < \delta_0$ and
  \begin{align*}
    &v-\bar{v} \in C([0,T];H^2(\R \times \T^2)), \quad \nabla(v-\bar{v}) \in L^2(0,T;H^1(\R \times \T^2)), \\
    &u-\bar{u} \in C([0,T];H^2(\R \times \T^2)), \quad \nabla(u-\bar{u}) \in L^2(0,T;H^2(\R \times \T^2)),
  \end{align*}
\end{proposition}
with
\begin{equation}\label{eq:apriori-assumption}
  \chi := \sup_{0\leq t\leq T} \norm{(v-\bar{v},u-\bar{u})(t,\cdot)}_{H^2} \leq \chi_0,
\end{equation}
then it holds that
\begin{equation}\label{eq:apriori-estimate}
  \begin{aligned}
    &\sup_{0\leq t\leq T} \norm{(v-\bar{v},u-\bar{u})(t,\cdot)}_{H^2}^2 + \delta_S \int_0^T |\dot{\X}(t)|^2 \,dt \\
    &+ \int_0^T \big(\pd_{1}v^S + \pd_{1}v^R \big)|v-\bar{v}|^2 + \norm{\nabla(v-\bar{v})}_{H^1}^2 + \norm{\nabla(u-\bar{u})}_{H^2}^2 \,dt \\
    &\leq C_0 \norm{v(0,\cdot)-\bar{v}(0,\cdot)}_{H^2}^2 + C_0 \norm{u(0,\cdot)-\bar{u}(0,\cdot)}_{H^2}^2 + C_0 \delta_R^\frac{1}{3} \,.
  \end{aligned}
\end{equation}
In addition, by \eqref{eq:X(t)-def}, we have
\begin{equation}\label{eq:xdot-property}
  |\dot{\X}(t)| \leq C_0 \norm{(v-\bar{v})(t,\cdot)}_{L^\infty}, \quad \forall\, t \leq T.
\end{equation}

Based on Propositions \ref{prop:local-exist} and \ref{prop:apriori-estimate}, we can employ continuation arguments to establish the global-in-time existence of perturbations as stated in \eqref{eq:regulartiy}. Furthermore, Proposition \ref{prop:apriori-estimate} can be used to validate \eqref{eq:asym-stable} for the large-time behavior. These proofs are typical and utilize the same arguments presented in the preceding paper \cite{Kang-Vasseur-Wang-2023, Teng-Yi-JEMS}. Consequently, we will omit these details and complete the proof of Theorem \ref{thm:main-ch2}.

Thus, the proof of Theorem \ref{thm:main-ch2} is reduced to that of Proposition \ref{prop:apriori-estimate}, which is the heart of the whole proof. The next section is then devoted to the proof of Proposition \ref{prop:apriori-estimate}.

\section{Proof of Proposition \ref{prop:apriori-estimate} }\label{sec:4}
The proof of Proposition \ref{prop:apriori-estimate} will be given in the following subsections. 
Throughout this section, $C$ denotes a positive generic constant independent of $\delta$ and $\nu$ (see \eqref{eq:delta_R,S-def} and \eqref{eq:nu-def} for their definitions). 
Also, we use $f \lesssim g$ to denote $f \leq Cg$ for some generic constant $C>0$.

\subsection{Reformulation of the problem}
Motivated by \cite{BD-2006}, we introduce the multi-dimensional effective velocity $h$ as follows:
\begin{equation}\label{eq:h-def}
  h := u -(2\mu+\lambda) \nabla v \,.
\end{equation}
By \eqref{eq:h-def}, the system \eqref{eq:3D-NS-v,u} can be deformed into
\begin{equation}\label{eq:3D-NS-v,h}
  \begin{cases}
    \rho (\md{v}) - \div h = (2\mu+\lambda) \Delta v, \\
    \rho (\md{h}) +\nabla p = R \,,
  \end{cases}
\end{equation}
where
\begin{equation}\label{eq:R-def}
  \begin{aligned}
    R :=&\; \frac{2\mu+\lambda}{v} (\nabla u^t \sdot \nabla v -\div u \nabla v) - \mu \curl\curl u   \\
    =&\; (2\mu+\lambda) (\nabla u^t - \div u \,\mathbb{I}) \nabla\ln v - \mu \curl\curl u\,, 
  \end{aligned}
\end{equation}
and $\mathbb{I}$ denote the $3 \times 3$ identity matrix. 
Similarly, we define
\begin{equation}
  h^S_1 :=  u^S_1 -(2\mu+\lambda) \pd_{1} v^S,\quad h^S := (h^S_1,0,0) .
\end{equation}
For $(v^S,h^S_1)$, it is easy to see that
\begin{equation}\label{eq:shock-v,h-eqn}
  \begin{cases}
    \rho^S (-\s \pd_{1}v^S + u^S_1\sdot \pd_{1}v^S) - \pd_{1}h^S_1 = (2\mu+\lambda) \pd_{11}v^S, \\
    \rho^S (-\s \pd_{1}h^s_1 + u^S_1\sdot \pd_{1}h^S_1) + \pd_{1}[p(v^S)] = 0.
  \end{cases}
\end{equation}
With the constant $\s_*$ defined in \eqref{eq:sigma*-def}, \eqref{eq:shock-v,h-eqn} can be rewritten as
\begin{equation}\label{eq:new-shock-2}
  \begin{cases}
    -\s_* \pd_{1}v^S - \pd_{1}h^S_1 = (2\mu+\lambda) \pd_{11}v^S, \\
    -\s_* \pd_{1}h^S_1 + \pd_{1}[p(v^S)] = 0.
  \end{cases}
\end{equation}
Let
\begin{equation}\label{eq:hbar-def}
  \bar{h}_1 := u^R_1 + h^S_1 -u_{1m}, \quad \bar{h} := (\bar{h}_1,0,0) .
\end{equation}
Thus,
\begin{align*}
  \bar{h} = \bar{u} - (2\mu+\lambda)\nabla v^S 
  = \bar{u} - (2\mu+\lambda)\nabla \bar{v} + (2\mu+\lambda)\nabla v^R .
\end{align*}
It can be directly verified that
\begin{equation}\label{eq:vbar,hbar-eqn}
  \begin{cases}
    \rho(\md{\bar{v}}) - \div \bar{h} + \rho\dot{\X}(t)\sdot \pd_{1}v^S = (2\mu+\lambda) \Delta (\bar{v} - v^R) + Q_1, \\
    \rho(\md{\bar{h}}) +\nabla p(\bar{v}) + \rho\dot{\X}(t)\sdot \pd_{1}h^S = Q_2\sdot e_1 \,,
  \end{cases}
\end{equation}
where $e_1 := (1,0,0)^t$,
\begin{equation}\label{eq:Q_1-def}
  Q_1 := \rho(u_1-u^R_1)\sdot \pd_{1}v^R - \rho(v-v^R)\sdot \pd_{1}u^R_1 + F\sdot \pd_{1}v^S,
\end{equation}
\begin{equation}\label{eq:Q_2-def}
  Q_2 := \rho(u_1-u^R_1)\sdot \pd_{1}u^R_1 + \rho(v-v^R)\sdot \pd_{1}p^R + F\sdot \pd_{1}h^S_1 + \left[ \pd_{1}\bar{p}-\pd_{1}p^R-\pd_{1}p^S \right] ,
\end{equation}
\begin{equation}\label{eq:F-def}
  F := -\s(\rho-\rho^S) + (\rho u_1 - \rho^Su^S_1) = \s_*\rho(v-v^S) + \rho(u_1-u^S_1) .
\end{equation}
Subtracting \eqref{eq:vbar,hbar-eqn} from \eqref{eq:3D-NS-v,h}, we obtain the perturbed system for $(\phi,h-\bar{h}) = (v-\bar{v}, h-\bar{h})$ 
\begin{equation}\label{eq:phi-h-eqn}
  \begin{cases}
    \rmd{\phi} - \div(h-\bar{h}) - \rho\dot{\X}(t)\sdot \pd_{1}v^S = (2\mu+\lambda) \Delta(\phi + v^R) - Q_1 \,, \\
    \rmd{(h-\bar{h})}  +\nabla(p-\bar{p}) - \rho\dot{\X}(t)\sdot \pd_{1}h^S = R - Q_2\sdot e_1 \,.
  \end{cases}
\end{equation}

\subsection{Estimates on $\norm{(v-\bar{v},h-\bar{h})}$ }

In this subsection, we focus on the following lemma.
\begin{lemma}\label{lem:phi,h-estimate}
  Under the assumptions of Proposition \ref{prop:apriori-estimate}, there exists constant $C>0$ (independent of $\nu,\delta,\chi$ and $T$)  such that for all $t \in [0,T]$, it holds
  \begin{equation}\label{eq:phi,h-L^2-est}
    \begin{aligned}
      &\sup_{t\in[0,T]} \norm{(v-\bar{v},h-\bar{h})(t,\cdot)}^2 + \delta_S \int_0^t |\dot{\X}(\tau)|^2 d\tau  \\
      &+ \int_0^t G_2(\tau) + G_3(\tau) + G^R(\tau) + G^S(\tau) + D(\tau) \,d\tau \\
      &\leq C \norm{(v-\bar{v},h-\bar{h})(0,\cdot)}^2 +C(\delta+\chi)\int_0^t \norm{\nabla(u-\bar{u})}_{H^1}^2 d\tau + C\delta_R^{\frac{1}{3}} .
    \end{aligned}
  \end{equation}
\end{lemma}
where
\begin{equation}\label{eq:new good-def}
  \begin{aligned}
    G_2(t) &:= \frac{\nu}{\delta_S} \int_{\Omega} \pd_{1}v^S \left| h_1-\bar{h}_1 - \frac{p-\bar{p}}{\s_*} \right|^2 \,dx , \\
    G_3(t) &:= \frac{\nu}{\delta_S} \int_{\Omega} \pd_{1}v^S (h_2^2+h_3^2) \,dx , \\
    G^R(t) &:= \int_{\Omega} \pd_{1}u^R_1 \sdot p(v|\bar{v}) \,dx ,\\
    G^S(t) &:= \int_{\Omega} \pd_{1}v^S \sdot |p-\bar{p}|^2 \,dx ,\\
    D(t) &:= \int_{\Omega} |\nabla(p-\bar{p})|^2 \,dx .
  \end{aligned}
\end{equation}

First of all, by direct calculations, we can obtain the following lemma.
\begin{lemma}\label{lem:4.1}
  Let $a(t,x)$ be the weighted function defined in \eqref{eq:a-weight-def}, then we have
  \begin{equation}\label{eq:0-v,h-exp}
    \frac{d}{dt} \int_{\Omega} a\rho \left( Q(v|\bar{v}) + \frac{1}{2} |h-\bar{h}|^2 \right) \,dx = \dot{\X}(t)\Y(t) + \B(t) +\S(t) -\G(t) - \D(t),
  \end{equation}
  where
  \begin{align*}
    \Y(t) &:= -\int_{\Omega} \pd_{1}a \sdot \rho \left( Q(v|\bar{v}) + \frac{1}{2} |h-\bar{h}|^2 \right)  \,dx - \int_{\Omega} a\rho \pd_{1}v^S\sdot p'(\bar{v})(v-\bar{v}) \,dx \\
    &\;\quad + \int_{\Omega} a\rho \pd_{1}h^S_1\sdot (h_1-\bar{h}_1) \,dx , \\
    \B(t) &:= \sum_{i=1}^8 \B_i(t), 
  \end{align*}
  with
  \begin{align*}
    \B_1(t) &:= \frac{1}{2\s_*} \int_{\Omega} \pd_{1}a\sdot |p-\bar{p}|^2 \,dx, \\
    \B_2(t) &:= \s_* \int_{\Omega} a \pd_{1}v^S \sdot p(v|\bar{v}) \,dx, \\
    \B_3(t) &:= \int_{\Omega} F \sdot \pd_{1}a \left( Q(v|\bar{v}) + \frac{1}{2} |h-\bar{h}|^2 \right) dx, \\
    \B_4(t) &:= (2\mu+\lambda) \int_{\Omega} \frac{\pd_{1}a}{p'(v)} (p-\bar{p})\sdot \pd_{1}(p-\bar{p}) \,dx, \\
    \B_5(t) &:= (2\mu+\lambda) \int_{\Omega} \pd_{1}a\sdot (p-\bar{p}) \left( \frac{1}{p'(v)} - \frac{1}{p'(\bar{v})} \right) \sdot \pd_{1}\bar{p} \,dx, \\
    \B_6(t) &:= (2\mu+\lambda) \int_{\Omega} a \pd_{1}(p-\bar{p}) \sdot \pd_{1}\bar{p} \left( \frac{1}{p'(v)} - \frac{1}{p'(\bar{v})} \right) dx, \\
    \B_7(t) &:= -\int_{\Omega} a (p-\bar{p}) \sdot \Delta v^R \,dx ,\\
    \B_8(t) &:= \int_{\Omega} a(h-\bar{h}) \sdot R \,dx,
  \end{align*}
  and
  \begin{align*}
    \S(t) &:= \int_{\Omega} a \left[ p'(\bar{v})(v-\bar{v})\sdot Q_1 - (h_1-\bar{h}_1)\sdot Q_2 \right]  \,dx, \\
    \D(t) &:= (2\mu+\lambda) \int_{\Omega} \frac{a}{|p'(v)|} |\nabla(p-\bar{p})|^2 \,dx, \\
    \G(t) &:= \sum_{i=1}^3 \G_i(t) + \G^R(t), \\
    \G_1(t) &:= \s_* \int_{\Omega} \pd_{1}a \sdot Q(v|\bar{v}) \,dx,  \\
    \G_2(t) &:= \frac{\s_*}{2} \int_{\Omega} \pd_{1}a \sdot \left| h_1-\bar{h}_1 - \frac{p-\bar{p}}{\s_*} \right|^2 \,dx,  \\
    \G_3(t) &:= \s_* \int_{\Omega} \pd_{1}a \sdot \frac{h_2^2+h_3^2}{2} \,dx, \\
    \G^R(t) &:= \int_{\Omega} a \pd_{1}u^R_1\sdot p(v|\bar{v}) \,dx .
  \end{align*}
\end{lemma}

\begin{proof}
  Multiplying \eqref{eq:phi-h-eqn}$_1$ by $-a(p-\bar{p})$, we have
  \begin{equation}
    \begin{aligned}
      &\pd_{t}\big[a\rho Q(v|\bar{v})\big] + \div \big[ a\rho Q(v|\bar{v}) u \big] +a(p-\bar{p})\div(h-\bar{h})  \\
      =&\; \pd_{1}a \big[u_1-\s-\dot{\X}(t)\big] \rho Q(v|\bar{v}) -a\rho \dot{\X}(t) \pd_{1}v^S\sdot (p-\bar{p}) \\
      &-a p(v|\bar{v}) \big[ \rmd{\bar{v}} \big] -a(2\mu+\lambda)(p-\bar{p})\Delta(\phi+v^R) +aQ_1(p-\bar{p}) \\
      =&\;(\rho u_1 - \rho\s)\pd_{1}a\sdot Q(v|\bar{v}) - \dot{\X}(t)\big[ \rho\pd_{1}a\sdot Q(v|\bar{v}) + a\rho\pd_{1}v^S\sdot p'(\bar{v})(v-\bar{v}) \big] \\
      &- a\div\bar{u}\sdot p(v|\bar{v}) - a(2\mu+\lambda)(p-\bar{p})\Delta(\phi+v^R) + aQ_1 p'(\bar{v})(v-\bar{v}),
    \end{aligned}
  \end{equation}
  Here, we use the fact that
  \begin{equation*}
    \rmd{\bar{v}} = -\rho \dot{\X}(t) \pd_{1}v^S +\div\bar{u} +Q_1.
  \end{equation*}
  Recall the definition of $F$ and $\s_*$ in \eqref{eq:F-def} and \eqref{eq:sigma*-def}, it holds
  \begin{equation*}
    \rho u_1 - \rho\s = -\s\rho^S + \rho^Su^S_1 -\s(\rho-\rho^S) + (\rho u_1 - \rho^Su^S_1) = -\s_* + F .
  \end{equation*}
  Also, with the fact $-\s_* \pd_{1}v^S = \pd_{1}u^S_1$ form \eqref{eq:v,u-ODE}, we have
  \begin{equation}\label{eq:part_1}
    \begin{aligned}
      &\pd_{t} \big[a\rho Q(v|\bar{v})\big] + \div \big[ a\rho Q(v|\bar{v}) u \big] +a(p-\bar{p})\div(h-\bar{h})  \\
      =&\;(F- \s_*)\pd_{1}a\sdot Q(v|\bar{v}) - \dot{\X}(t)\big[ \rho\pd_{1}a\sdot Q(v|\bar{v}) + a\rho\pd_{1}v^S\sdot p'(\bar{v})(v-\bar{v}) \big] \\
      &+ a(\s_*\pd_{1}v^S - \pd_{1}u^R_1) \sdot p(v|\bar{v}) - a(2\mu+\lambda)(p-\bar{p})\Delta v^R \\
      &-(2\mu+\lambda) \div[a(p-\bar{p})\nabla\phi] + (2\mu+\lambda) \nabla[a(p-\bar{p})] \sdot \nabla\phi + aQ_1 p'(\bar{v})(v-\bar{v}).
    \end{aligned}
  \end{equation}
  Since
  \begin{equation*}
    \nabla v = \frac{\nabla p(v)}{p'(v)} = \frac{\nabla p(v)}{-\gamma p^{1+\frac{1}{\gamma}}(v)},
  \end{equation*}
  we can make the following transformations
  \begin{equation}
    \begin{aligned}
      &(2\mu+\lambda) \nabla [a(p-\bar{p})] \sdot \nabla\phi \\
      =&\, (2\mu+\lambda) a \nabla(p-\bar{p}) \sdot \left( \frac{\nabla p(v)}{p'(v)} - \frac{\nabla p(\bar{v})}{p'(\bar{v})} \right) + (2\mu+\lambda) \pd_{1}a \sdot (p-\bar{p}) \left( \frac{\pd_{1} p(v)}{p'(v)} - \frac{\pd_{1} p(\bar{v})}{p'(\bar{v})} \right)  \\
      =&\, -(2\mu+\lambda) a \frac{|\nabla(p-\bar{p})|^2}{|p'(v)|} - (2\mu+\lambda) a \pd_{1}(p-\bar{p}) \pd_{1}\bar{p} \sdot \left( \frac{1}{p'(v)} - \frac{1}{p'(\bar{v})} \right) \\
      & + (2\mu+\lambda) \pd_{1}a \sdot (p-\bar{p}) \frac{\pd_{1}(p-\bar{p})}{p'(v)} + (2\mu+\lambda) \pd_{1}a \sdot (p-\bar{p}) \pd_{1}\bar{p} \sdot \left( \frac{1}{p'(v)} - \frac{1}{p'(\bar{v})} \right).
    \end{aligned} 
  \end{equation}
  Then, multiplying \eqref{eq:phi-h-eqn}$_2$ by $a(h-\bar{h})$, we have
  \begin{equation}\label{eq:part_2}
    \begin{aligned}
      &\pd_{t}\left[ a\rho \frac{|h-\bar{h}|^2}{2} \right] + \div \left[ a\rho \frac{|h-\bar{h}|^2}{2} u \right] + \div \big[ a(p-\bar{p})(h-\bar{h}) \big] - a(p-\bar{p})\div(h-\bar{h}) \\
      =&\; \pd_{1}a \big[u_1-\s-\dot{\X}(t)\big] \rho \frac{|h-\bar{h}|^2}{2} + a\rho \dot{\X}(t) \pd_{1}h^S\sdot (h-\bar{h}) + \pd_{1}a(p-\bar{p})(h_1-\bar{h}_1) \\
      &+ aR\sdot (h-\bar{h}) -Q_2(h_1-\bar{h}_1) \\
      =&\; \pd_{1}a \sdot F\frac{|h-\bar{h}|^2}{2} +  \dot{\X}(t) \left[ a\rho\pd_{1}h^S\sdot (h-\bar{h}) - \rho \pd_{1}a\sdot \frac{|h-\bar{h}|^2}{2} \right] -\s_*\pd_{1}a \sdot F\frac{|h-\bar{h}|^2}{2}\\
      &+ \pd_{1}a(p-\bar{p})(h_1-\bar{h}_1) + aR\sdot (h-\bar{h}) -Q_2(h_1-\bar{h}_1) .
    \end{aligned}
  \end{equation}
  By direct computation, we have
  \begin{equation}
    \begin{aligned}
      &-\s_*\pd_{1}a \sdot F\frac{|h-\bar{h}|^2}{2} + \pd_{1}a(p-\bar{p})(h_1-\bar{h}_1) \\
      =& -\frac{\s_*}{2} \pd_{1}a \left| h_1-\bar{h}_1 - \frac{p-\bar{p}}{\s_*} \right|^2 +\pd_{1}a \sdot \frac{|p-\bar{p}|^2}{2\s_*} -\s_*\pd_{1}a\sdot \frac{h_2^2+h_3^2}{2}, 
    \end{aligned}
  \end{equation}
  Adding \eqref{eq:part_1} and \eqref{eq:part_2} together, and integrating the resultant equation on $\Omega$, we can achieve \eqref{eq:0-v,h-exp}. 
\end{proof}

For notational simplicity, we decompose the $\Y(t)$ in Lemma \ref{lem:4.1} as follows:
\begin{align*}
  \Y(t) &= -\int_{\Omega}\pd_{1}a \sdot \rho \left( Q(v|\bar{v}) + \frac{1}{2} |h-\bar{h}|^2 \right)  \,dx - \int_{\Omega} a\rho \pd_{1}v^S\sdot p'(\bar{v})(v-\bar{v}) \,dx \\
  &\;\quad + \int_{\Omega} a\rho \pd_{1}h^S_1\sdot (h_1-\bar{h}_1) \,dx =: \sum_{i=1}^5 \Y_i(t),
\end{align*}
where
\begin{align*}
  \Y_1(t) &:=  -\int_{\Omega} a\rho \pd_{1}v^S\sdot p'(v^S)(v-\bar{v}) \,dx,  \\
  \Y_2(t) &:=  \int_{\Omega} \frac{a\rho}{\s_*}\sdot \pd_{1}h^S_1\sdot (p-\bar{p}) \,dx,  \\
  \Y_3(t) &:= \int_{\Omega} a\rho \pd_{1}v^S\sdot [p'(v^S)-p'(\bar{v})](v-\bar{v}) \,dx,  \\
  \Y_4(t) &:= \int_{\Omega} a\rho \pd_{1}h^S_1\sdot \left( h_1-\bar{h}_1 - \frac{p-\bar{p}}{\s_*} \right)  \,dx,  \\
  \Y_5(t) &:= -\frac{1}{2} \int_{\Omega} \pd_{1}a \sdot \rho \left( h_1-\bar{h}_1 - \frac{p-\bar{p}}{\s_*} \right) \left( h_1-\bar{h}_1 + \frac{p-\bar{p}}{\s_*} \right)  \,dx,  \\
  \Y_6(t) &:= -\int_{\Omega} \pd_{1}a \sdot \rho \left( Q(v|\bar{v}) + \frac{h_2^2+h_3^2}{2}  \right) \,dx - \int_{\Omega} \pd_{1}a \sdot \rho \frac{|p-\bar{p}|^2}{2\s_*^2} \,dx.         
\end{align*}
Then, according to \eqref{eq:X(t)-def}, we have
\begin{equation}\label{eq:x=y_1+y_2}
  \dot{\X}(t) = -\frac{M}{\delta_S} \big[ \Y_1(t)+\Y_2(t) \big], 
\end{equation}
and
\begin{equation}
  \dot{\X}(t) \Y(t) = -\frac{M}{\delta_S} |\dot{\X}(t)|^2 + \dot{\X}(t) \sum_{i=3}^5 \Y_i(t).
\end{equation}

Secondly, we give the following leading order estimates, which is crucial for the proof of Lemma \ref{lem:phi,h-estimate}.
\begin{lemma}\label{lem:key lemma}
  There exists uniform-in-time $C_0>0$ such that for all $t \in [0,T]$, it holds
  \begin{equation}\label{eq:key lem-ineq}
    \begin{aligned}
      &-\frac{\delta_S}{2M} |\dot{\X}(t)|^2 +\B_1(t) +\B_2(t) - \G_1(t) - \frac{3}{4}\D(t)  
      \leq -C_0 G^S(t) + C\int_{\Omega} \pd_{1}a\sdot |p-\bar{p}|^3 \,dx ,
    \end{aligned}
  \end{equation}
  where $\displaystyle G^S(t) := \int_{\Omega} \pd_{1}v^S \sdot |p-\bar{p}|^2 \,dx. $ 
\end{lemma}

{\noindent \bf Proof.~} First, we introduce the following time-dependent change of variables:
\begin{equation}\label{eq:change of variable}
  \begin{cases}
    y_1 = \frac{p(v_m)-p(v^S)}{\delta_S} 
    = \frac{p(v_m)-p(v^s(x_1-\s t-\X(t)) )}{\delta_S}, \\
    y_2 =x_2, \\
    y_3 =x_3,
  \end{cases}
\end{equation}
and set
\begin{equation}\label{eq:w-def}
  w(t,y) := (p-\bar{p})(t,x(y)), \qquad y':=(y_2,y_3).
\end{equation}
By direct computation, we know that $y=(y_1,y') \in (0,1) \times \T^2$. Also, it follows from \eqref{eq:a-weight-def} that $a=1+\nu y_1$ and
\begin{equation}
  \frac{dy_1}{dx_1} = -\frac{1}{\delta_S}\pd_{x_1}p^S, \quad \pd_{x_1}a = \nu \frac{dy_1}{dx_1}, \quad |a-1|\leq \nu = \sqrt{\delta_S}.
\end{equation}
In the above new variables, we have
\begin{equation}
  \int_{\T^2}\int_{\R} \pd_{1}a |p-\bar{p}|^2 \,dx_1dx' = \nu \int_{\T^2}\int_0^1 w^2 \,dy_1dy', 
\end{equation}
and so
\begin{equation}
  \int_{\T^2}\int_{\R} \pd_{1}v^S |p-\bar{p}|^2 \,dx_1dx' \sim  \delta_S \int_{\T^2}\int_0^1 w^2 \,dy_1dy'. 
\end{equation}
For a finer estimate, we introduce two constants 
\begin{equation*}
  \s_m := \sqrt{-p'(v_m)}, \quad \alpha_m := \frac{\gamma+1}{2\gamma\s_mp(v_m)}.
\end{equation*}
Then, from \eqref{eq:sigma*>0} and the fact that  $\s_m^2 = -p'(v_m) = \gamma p^{1+\frac{1}{\gamma}}(v_m)$, we obtain
\begin{equation}\label{eq:equi-sigma}
  |\s_*-\s_m| \leq C\delta_S,
\end{equation}
and
\begin{equation}\label{eq:equi-sigma^2}
  |\s_m^2 + p'(\bar{v})| \leq C\delta,\quad \left| \frac{1}{\s_m^2} - \frac{1}{\gamma p^{1+\frac{1}{\gamma}}(\bar{v})} \right| \leq C\delta.
\end{equation}


$\bullet$ {\bf Estimates on $-\frac{\delta_S}{2M} |\dot{\X}(t)|^2$ }:
According to \eqref{eq:x=y_1+y_2}, we only have to estimate $\Y_1(t)$ and $\Y_2(t)$.
For one thing, 
\begin{align*}
  \Y_1(t) = -\int_{\Omega} a\rho \pd_{1}p^S (v-\bar{v}) \,dx = \delta_S \int_{\T^2}\int_0^1 \frac{a}{v} (v-\bar{v}) \,dy_1dy' .
\end{align*}
Thus, by Taylor expansion and that $|v-\bar{v}|\sim |p-\bar{p}|$, it holds
\begin{equation*}
  |p(v)-p(\bar{v})-p'(\bar{v})(v-\bar{v})| \leq  C|p-\bar{p}|^2,
\end{equation*}
and so
\begin{equation}\label{eq:equi-phi}
  \left| v-\bar{v} + \frac{p-\bar{p}}{\s_m^2} \right| \leq C(\delta+\chi)|p-\bar{p}|,  \quad
  \left| v-\bar{v} + \frac{p-\bar{p}}{\s_*^2} \right| \leq C(\delta+\chi)|p-\bar{p}|.
\end{equation}
Together with $|a-1|\leq \nu$ form \eqref{eq:a-weight-def}, we get
\begin{equation}\label{eq:Y_1-est}
  \left| \Y_1(t) + \frac{\delta_S}{\s_m^2 v_m} \int_{\T^2}\int_0^1 w \,dy_1dy' \right| \leq C\delta_S (\nu+\delta+\chi) \int_{\T^2}\int_0^1 |w| \,dy_1dy' .
\end{equation}
For another thing, by \eqref{eq:new-shock-2}, we yield
\begin{align*}
  \Y_2(t) = \int_{\Omega} \frac{a}{v\s_*^2}\sdot \pd_{1}p^S\sdot (p-\bar{p}) \,dx = -\frac{\delta_S}{\s_*^2} \int_{\T^2}\int_0^1 \frac{a}{v}\sdot w \,dy_1dy' ,
\end{align*}
and thus
\begin{equation}\label{eq:Y_2-est}
  \left| \Y_2(t) + \frac{\delta_S}{\s_m^2 v_m} \int_{\T^2}\int_0^1 w \,dy_1dy' \right| \leq C\delta_S (\nu+\delta+\chi) \int_{\T^2}\int_0^1 |w| \,dy_1dy' .
\end{equation}
By \eqref{eq:x=y_1+y_2}, \eqref{eq:Y_1-est} and \eqref{eq:Y_2-est}, we have
\begin{align*}
  \left| \dot{\X}(t) - \frac{2M}{\s_m^2 v_m} \int_{\T^2}\int_0^1 w \,dy_1dy' \right| &= \left| \sum_{i=1}^2 \frac{M}{\delta_S} \left( \Y_i(t) + \frac{\delta_S}{\s_m^2 v_m} \int_{\T^2}\int_0^1 w \,dy_1dy' \right) \right| \\
  &\leq C(\nu+\delta+\chi) \int_{\T^2}\int_0^1 |w| \,dy_1dy' ,
\end{align*}
which leads to
\begin{align*}
  \left( \left| \frac{2M}{\s_m^2 v_m} \int_{\T^2}\int_0^1 w \,dy_1dy' \right| -|\dot{\X}(t)| \right)^2  &\leq C(\nu+\delta+\chi)^2 \left( \int_{\T^2}\int_0^1 |w| \,dy_1dy' \right)^2 \\
  &\leq C(\nu+\delta+\chi)^2 \int_{\T^2}\int_0^1 |w|^2 \,dy_1dy' .
\end{align*}
Finally, with the algebraic inequality $\frac{p^2}{2} - q^2 \leq (p-q)^2, \forall\, p,q \geq 0$, we arrive at
\begin{equation*}
  \frac{2M^2}{\s_m^4 v_m^2} \left( \int_{\T^2}\int_0^1 w \,dy_1dy' \right)^2 - |\dot{\X}(t)|^2 \leq C(\nu+\delta+\chi)^2 \int_{\T^2}\int_0^1 |w|^2 \,dy_1dy',
\end{equation*}
and consequently
\begin{equation}\label{eq:xdot-est}
   -\frac{\delta_S }{2M} |\dot{\X}(t)|^2 \leq -\frac{M\delta_S}{\s_m^4 v_m^2} \left( \int_{\T^2}\int_0^1 w \,dy_1dy' \right)^2 + C\delta_S (\nu+\delta+\chi)^2 \int_{\T^2}\int_0^1 |w|^2 \,dy_1dy'.
\end{equation}


$\bullet$ {\bf Estimates on $\B_1(t)-\G_1(t)$ and $\B_2(t)$  }:
First, by \eqref{eq:equi-sigma}, 
\begin{equation}\label{eq:b_1-est}
  \begin{aligned}
    \B_1(t) &= \frac{1}{2\s_*} \int_{\Omega} \pd_{1}a\sdot |p-\bar{p}|^2 \,dx = \frac{\nu}{2\s_*} \int_{\T^2}\int_0^1 w^2 \,dy_1dy' \\
    &\leq \frac{\nu}{2\s_m} \int_{\T^2}\int_0^1 w^2 \,dy_1dy' + C\nu\delta_S \int_{\T^2}\int_0^1 w^2 \,dy_1dy'.
  \end{aligned}
\end{equation}
For $\B_2(t)$, by \eqref{eq:relative-2}, we have
\begin{equation}\label{eq:b_2-est}
  \begin{aligned}
    \B_2(t) &= \s_* \int_{\Omega} a \pd_{1}v^S \sdot p(v|\bar{v}) \,dx = \s_* \delta_S \int_{\T^2}\int_0^1 \frac{a}{|p'(v^S)|} \sdot p(v|\bar{v}) \,dy_1dy' \\
    &\leq \s_* \delta_S(1+\nu) \int_{\T^2}\int_0^1 \left( \frac{\gamma+1}{2\gamma}\frac{1}{p(v^S)}+ C\chi \right) \frac{|p-\bar{p}|^2}{|p'(v^S)|} \,dy_1dy' \\
    &\leq \s_* \delta_S(1+\nu) \int_{\T^2}\int_0^1 \left( \alpha_m \frac{\s_m p(v_m)}{p(v^S)}+ C\chi \right) \frac{w^2}{|p'(v^S)|} \,dy_1dy' \\
    &\leq \delta_S \alpha_m \int_{\T^2}\int_0^1 w^2 \,dy_1dy' + C(\nu+\delta+\chi) \int_{\T^2}\int_0^1 w^2 \,dy_1dy'.
  \end{aligned}
\end{equation}
For $\G_1(t)$, by \eqref{eq:relative-3} and \eqref{eq:equi-sigma^2}, we get
\begin{equation}\label{eq:g_1-est}
  \begin{aligned}
    \G_1(t) &= \s_* \int_{\Omega} \pd_{1}a \sdot Q(v|\bar{v}) \,dx, \\
    &\geq \s_* \int_{\Omega} \pd_{1}a \sdot \frac{|p-\bar{p}|^2}{2\gamma p^{1+\frac{1}{\gamma}}(\bar{v})} \,dx - \s_* \int_{\Omega} \pd_{1}a \sdot \frac{1+\gamma}{3\gamma^2}\frac{(p-\bar{p})^3}{p^{2+\frac{1}{\gamma}}(\bar{v})} \,dx  \\
    &\geq \frac{\nu\s_*}{2\s_m^2} (1-C\delta)\int_{\T^2}\int_0^1 w^2 \,dy_1dy' - C\nu \int_{\T^2}\int_0^1 w^3 \,dy_1dy' \\
    &\geq \frac{\nu}{2\s_m} \int_{\T^2}\int_0^1 w^2 \,dy_1dy' - C\nu\delta_S \int_{\T^2}\int_0^1 w^2 \,dy_1dy' - C\nu \int_{\T^2}\int_0^1 w^3 \,dy_1dy',
  \end{aligned}
\end{equation}
Hence, together with \eqref{eq:b_1-est}, we acquire
\begin{equation}\label{eq:b_new-est}
  \B_{new}(t) := \B_1(t)-\G_1(t) \leq C\nu\delta_S \int_{\T^2}\int_0^1 w^2 \,dy_1dy' + C\nu \int_{\T^2}\int_0^1 w^3 \,dy_1dy'.
\end{equation}


$\bullet$ {\bf Estimate on $\D(t)$ }:
First, by \eqref{eq:a-weight-def} and the change of variables \eqref{eq:change of variable}, we get
\begin{align*}
    &\D(t) \geq (2\mu+\lambda) \int_{\T^2}\int_{\R}
    \frac{|\pd_{x_1} \big( p(v)-p(\bar{v}) \big)|^2}{\gamma p^{1+\frac{1}{\gamma}}(v)} + \frac{|\nabla_{x'} \big(p(v)-p(\bar{v})\big)|^2}{\gamma p^{1+\frac{1}{\gamma}}(v)}dx_1dx' \\
    & =\underbrace{(2\mu+\lambda)\!\!\int_{\T^2}\!\int_{0}^1 \frac{|\pd_{y_1}w|^2}{\gamma p^{1+\frac{1}{\gamma}}(v)} \Big(\frac{dy_1}{dx_1}\Big) dy_1dy'}_{\D_{\rm I}(t)}
    + \underbrace{(2\mu+\lambda)\!\!\int_{\T^2}\!\int_{0}^1\frac{|\nabla_{y'}w|^2}{\gamma p^{1+\frac{1}{\gamma}}(v)} \Big(\frac{dx_1}{dy_1}\Big)dy_1dy'}_{\D_{\rm II}(t)}.
\end{align*}
For one thing, integrating \eqref{eq:new-shock-2} over $(-\infty,x_1]$ leads to
\begin{equation*}
  (2\mu+\lambda)\pd_{x_1} v^S=-\s_* (v^S-v_m)-\frac{1}{\s_* }\big(p(v^S)-p(v_m)\big).
\end{equation*}
For another thing, 
\begin{equation*}
  \pd_{x_1} v^S=\frac{\pd_{x_1} p(v^S)}{p'(v^S)}=\frac{\delta_S}{\gamma p^{1+\frac{1}{\gamma}}(v^S)} \frac{dy_1}{dx_1} \;.
\end{equation*}
Therefore, we have
\begin{align*}
  (2\mu+\lambda) \frac{\delta_S}{\gamma p^{1+\frac{1}{\gamma}}(v^S)} \frac{dy_1}{dx_1}
  & = -\s_* (v^S-v_m)-\frac{1}{\s_* }\big(p(v^S)-p(v_m)\big)  \\
  & = -\frac{1}{\s_* }\Big(\s_* ^2(v^S-v_m)+\big(p(v^S)-p(v_m)\big)\Big),
\end{align*}
which together with $\displaystyle \s_* ^2=-\frac{p(v_m)-p(v_+)}{v_m-v_+}$ implies
\begin{align*}
  & (2\mu+\lambda) \frac{\delta_S}{\gamma p^{1+\frac{1}{\gamma}}(v^S)} \frac{dy_1}{dx_1} \\
  = &\frac{-1}{\s_* (v_+-v_m)}\Big(\big(p(v_m)-p(v_+)\big)(v^S-v_m) + \big(p(v^S)-p(v_m)\big)(v_+-v_m)\Big) \\
  = &\frac{-1}{\s_* (v_+-v_m)}\Big(\big(p(v^S)-p(v_+)\big)(v^S-v_m) - (v^S-v_+)\big(p(v^S)-p(v_m)\big)\Big).
\end{align*}
By \eqref{eq:change of variable} and \eqref{eq:delta_R,S-def}, it holds $y_1=\frac{p(v_m)-p(v^S)}{\delta_S}$ and $1-y_1=\frac{p(v^S)-p(v_+)}{\delta_S}$, then
\begin{equation*}
  \frac{1}{y_1(1-y_1)}\frac{2\mu+\lambda}{\gamma p^{1+\frac{1}{\gamma}}(v^S)} \frac{dy_1}{dx_1}
  = \frac{\delta_S}{\s_* (v_+-v_m)} \left( \frac{v^S-v_m}{p(v^S)-p(v_m)} - \frac{v^S-v_+}{p(v^S)-p(v_+)} \right).
\end{equation*}
Thus
\begin{align*}
  & \left| \frac{1}{y_1(1-y_1)} \frac{2\mu+\lambda}{\gamma p^{1+\frac{1}{\gamma}}(v^S)} \frac{dy_1}{dx_1} - \frac{\delta_S p''(v_m)}{2\s_m(p'(v_m))^2} \right|  \\
  \leq& \left|\frac{1}{y_1(1-y_1)}\frac{2\mu+\lambda}{\gamma p^{1+\frac{1}{\gamma}}(v^S)}\frac{dy_1}{dx_1} -\frac{\delta_S p''(v_m)}{2 \s_*(p'(v_m))^2}\right|
  +\frac{\delta_S p''(v_m)}{2(p'(v_m))^2}\left|\frac{1}{\s_* }-\frac{1}{\s_m}\right|  \\
  \leq& \left| \frac{\delta_S}{\s_* (v_+-v_m)} \left( \frac{v^S-v_m}{p(v^S)-p(v_m)} +\frac{v^S-v_+}{p(v_+)-p(v^S)} \right) \right| 
  +\frac{\delta_S p''(v_m)}{2(p'(v_m))^2}\left|\frac{1}{\s_* }- \frac{1}{\s_m}\right| \,.
\end{align*}
Using Lemma \ref{lem:inverse-pressure}, we have
\begin{align*}
  \left|\frac{\delta_S}{\s_* (v_+-v_m)}\left(\frac{v^S-v_m}{p(v^S)-p(v_m)} + \frac{v^S-v_+}{p(v_+)-p(v^S)} \right) -\frac{\delta_S p''(v_m)}{2\s_*(p'(v_m))^2} \right| \le C\delta_S^2.
\end{align*}
Also, it follows from \eqref{eq:equi-sigma} that
\begin{equation*}
  \frac{\delta_S p''(v_m)}{2(p'(v_m))^2}\left|\frac{1}{\s_* }-\frac{1}{\s_m}\right| \leq C\delta_S^2 \,,
\end{equation*}
Hence, we get
\begin{equation}\label{coeffi-1}
    \left|\frac{1}{y_1(1-y_1)}\frac{2\mu+\lambda}{\gamma p^{1+\frac{1}{\gamma}}(v^S)}\frac{dy_1}{dx_1}
    -\frac{\delta_S p''(v_m)}{2\s_m(p'(v_m))^2}\right|\leq C\delta_S^2.
\end{equation}
In addition, 
\begin{equation*}
  \left| \left(\frac{p(v^S)}{p(v)}\right)^{1+\frac{1}{\gamma}}-1 \right| \leq C|v-v^S| \leq C(\chi+\delta_R) \,.
\end{equation*}
Thus, we obtain
\begin{align*}
  \D_{\rm I}(t) 
  & =\int_{\T^2}\int_0^1 y_1(1-y_1) |\pd_{y_1}w|^2 \left(\frac{p(v^S)}{p(v)}\right)^{1+\frac{1}{\gamma}} \frac{1}{y_1(1-y_1)}\frac{2\mu+\lambda}{\gamma p^{1+\frac{1}{\gamma}}(v^S)}
  \Big(\frac{dy_1}{dx_1}\Big) dy_1dy'  \\
  & \geq \Big( 1-C(\chi+\delta_R) \Big) \left(\frac{\delta_S p''(v_m)}{2\s_m(p'(v_m))^2}-C\delta_S^2\right) \int_{\T^2}\int_0^1 y_1(1-y_1) |\pd_{y_1}w|^2 dy_1dy'.  
\end{align*}
A direct computation gives rise to
\begin{equation*}
  \frac{p''(v_m)}{2\s_m(p'(v_m))^2} = \frac{\gamma+1}{2\gamma\s_mp(v_m)} = \a_m,
\end{equation*}
which leads to
\begin{equation*}
  \D_{\rm I}(t) \geq \a_m \delta_S(1-C(\delta+\chi))\int_{\T^2}\int_0^1y_1(1-y_1)|\pd_{y_1}w|^2dy_1dy'.
\end{equation*}
It follows from \eqref{coeffi-1} that
\begin{equation*}
  y_1(1-y_1) \frac{dx_1}{dy_1}\geq\frac{2\mu+\lambda}{\gamma p^{1+\frac{1}{\gamma}}(v^S)}\frac{1}{\a_m\delta_S + C\delta_S^2} \geq \frac{2\mu+\lambda}{2\a_m\delta_S|p'(v_m)|}.
\end{equation*}
Hence, we get
\begin{align*}
  \D_{\rm II}(t) &=(2\mu+\lambda)\int_{\T^2}\int_{0}^1\frac{|\nabla_{y'}w|^2}{y_1(1-y_1)} \left(\frac{p(v_m)}{p(v)}\right)^{1+\frac{1}{\gamma}} \frac{y_1(1-y_1)}{\gamma p^{1+\frac{1}{\gamma}}(v_m)} \Big(\frac{dx_1}{dy_1}\Big)dy_1dy' \\
  &\geq (1-C(\delta+\chi))(2\mu+\lambda)\int_{\T^2}\int_{0}^1\frac{|\nabla_{y'}w|^2}{y_1(1-y_1)} \frac{2\mu+\lambda}{|p'(v_m)|^2}\frac{1}{2\a_m\delta_S}dy_1dy' \\
  &\geq (1-C(\delta+\chi))\frac{\s_m}{\delta_S}\frac{(2\mu+\lambda)^2}{p''(v_m)} \int_{\T^2}\int_0^1\frac{|\nabla_{y'}w|^2}{y_1(1-y_1)}dy_1dy'.
\end{align*}
Combining the estimates on $\D_{\rm I}(t)$ and $\D_{\rm II}(t)$, we arrive at
\begin{equation}\label{eq:D-est}
  \begin{aligned}
    \D(t)&\geq \a_m\delta_S(1-C(\delta+\chi))\int_{\T^2}\int_0^1y_1(1-y_1)|\pd_{y_1}w|^2dy_1dy' \\
    &+ (1-C(\delta+\chi))\frac{\s_m}{\delta_S}\frac{(2\mu+\lambda)^2}{p''(v_m)} \int_{\T^2}\int_0^1\frac{|\nabla_{y'}w|^2}{y_1(1-y_1)}dy_1dy'.
  \end{aligned}
\end{equation}


$\bullet$ {\bf Proof of Lemma \ref{lem:key lemma} }:
First of all, in view of \eqref{eq:b_2-est}, \eqref{eq:b_new-est} and \eqref{eq:D-est}, it follows that
\begin{align*}
  &\B_1(t)+\B_2(t)-\G_1(t)-\frac{3}{4}\D(t) \\
  \leq&\; \a_m\delta_S \int_{\T^2}\int_0^1 w^2 \,dy_1dy' -\frac{3}{4} \a_m\delta_S \int_{\T^2}\int_0^1 y_1(1-y_1)|\pd_{y_1}w|^2 dy_1dy' \\
  & - \frac{C}{\delta_S} \int_{\T^2}\int_0^1 \frac{|\nabla_{y'}w|^2}{y_1(1-y_1)} dy_1dy' + C(\delta+\chi) \int_{\T^2}\int_0^1 y_1(1-y_1)|\pd_{y_1}w|^2dy_1dy'  \\
  &+ C\nu\delta_S \int_{\T^2}\int_0^1 w^2 \,dy_1dy' + C\nu \int_{\T^2}\int_0^1 w^3 \,dy_1dy'. 
\end{align*}
We denote the average of $w$ by $\bar{w}$, i.e.
\begin{equation*}
  \bar{w} := \int_{\T^2}\int_0^1 w \,dy_1dy' ,
\end{equation*}
It is east to show that
\begin{equation*}
  \int_{\T^2}\int_0^1 |w-\bar{w}|^2 \,dy_1dy' = \int_{\T^2}\int_0^1 w^2 \,dy_1dy' - \bar{w}^2.
\end{equation*}
Taking $\nu,\delta,\chi$ suitably small and using the weighted Poincar\'e inequality \eqref{eq:weighted Poincare}, we acquire
\begin{equation}\label{eq:4.38}
  \begin{aligned}
    &\B_1(t)+\B_2(t)-\G_1(t)-\frac{3}{4}\D(t) \\
    \leq&\; \frac{6}{5}\a_m\delta_S \int_{\T^2}\int_0^1 w^2 \,dy_1dy' -\frac{5}{4} \a_m\delta_S \int_{\T^2}\int_0^1 |w-\bar{w}|^2 dy_1dy' \\
    & - \frac{C}{\delta_S} \int_{\T^2}\int_0^1 \frac{|\nabla_{y'}w|^2}{y_1(1-y_1)} dy_1dy'  + C\nu \int_{\T^2}\int_0^1 w^3 \,dy_1dy'. 
  \end{aligned}
\end{equation}
Finally, choosing $\displaystyle M=\frac{5}{4} \a_m\s_m^4v_m^2$ and combining \eqref{eq:xdot-est} and \eqref{eq:4.38}, we arrive at
\begin{align*}
  &-\frac{\delta_S }{2M} |\dot{\X}(t)|^2 + \B_1(t)+\B_2(t)-\G_1(t)-\frac{3}{4}\D(t) \\
  \leq&\; \left( -\frac{\a_m\delta_S}{20} + C\delta_S (\nu+\delta+\chi)^2 \right) \int_{\T^2}\int_0^1 w^2 \,dy_1dy' \\
  &- \frac{C}{\delta_S} \int_{\T^2}\int_0^1 \frac{|\nabla_{y'}w|^2}{y_1(1-y_1)} dy_1dy' + C\nu \int_{\T^2}\int_0^1 w^3 \,dy_1dy' \\
  \leq&\; -C_0 \delta_S \int_{\T^2}\int_0^1 w^2 \,dy_1dy' + C\nu \int_{\T^2}\int_0^1 w^3 \,dy_1dy'.
\end{align*}
Notice that
\begin{equation*}
  G^S(t) \sim \delta_S \int_{\T^2}\int_0^1 w^2 \,dy_1dy',
\end{equation*}
so we obtain the desired inequality \eqref{eq:key lem-ineq} and complete the proof of Lemma \ref{lem:key lemma}.
\hfill  $\blacksquare$


At last, we give the proof of Lemma \ref{lem:phi,h-estimate}. To do so, we only need to estimate the remaining terms in \eqref{eq:phi,h-L^2-est} individually.

\vspace{0.5em}

{\noindent \bf Proof of Lemma \ref{lem:phi,h-estimate}:} First, by \eqref{eq:0-v,h-exp} and \eqref{eq:x=y_1+y_2}, it holds
\begin{align*}
  &\frac{d}{dt} \int_{\Omega} a\rho \left( Q(v|\bar{v}) + \frac{1}{2} |h-\bar{h}|^2 \right) \,dx = \dot{\X}(t)\Y(t) + \B(t) +\S(t) -\G(t) - \D(t) \\
  =& - \frac{\delta_S}{2M} |\dot{\X}(t)|^2 +\B_1(t) +\B_2(t) - \G_1(t) - \frac{3}{4}\D(t)  \\
  & - \frac{\delta_S}{2M} |\dot{\X}(t)|^2 + \dot{\X}(t)\sum_{i=3}^6 \Y_i(t) + \sum_{i=3}^8 \B_i(t) + \S(t) -\sum_{i=2}^3 \G_i(t) -\G^R(t) -\frac{1}{4}\D(t).
\end{align*}
Using Lemma \ref{lem:key lemma} and Cauchy's inequality, we get
\begin{align*}
  &\frac{d}{dt} \int_{\Omega} a\rho \left( Q(v|\bar{v}) + \frac{1}{2} |h-\bar{h}|^2 \right) \,dx  \\
  \leq& -C_0 G^S(t) + C\int_{\Omega} \pd_{1}a\sdot |p-\bar{p}|^3 \,dx - \frac{\delta_S}{4M} |\dot{\X}(t)|^2 + \frac{M}{\delta_S} \sum_{i=3}^6 |\Y_i(t)|^2 \\
  &+ \sum_{i=3}^8 \B_i(t) + \S(t) -\sum_{i=2}^3 \G_i(t) -\frac{1}{4}\D(t),
\end{align*}
Before estimating the above bad terms individually, we notice that, by \eqref{eq:apriori-assumption} and the Gagliardo-Nirenberg inequality \eqref{eq:G-N ineq}, we have
\begin{equation}\label{eq:4.41}
  \begin{aligned}
    &C\int_{\Omega} \pd_{1}a\sdot |p-\bar{p}|^3 \,dx \leq \frac{C}{\nu} \norm{p-\bar{p}}_{L^\infty}^2 \int_{\Omega} \pd_{1}v^S\sdot |p-\bar{p}|^3 \,dx  \\
    \leq&\; \frac{C}{\nu} \norm{\nabla(p-\bar{p})} \Big( \norm{p-\bar{p}} + \norm{\nabla^2(p-\bar{p})} \Big) \left( \int_{\Omega} \pd_{1}v^S\sdot |p-\bar{p}|^2 \,dx  \right)^\frac{1}{2} \left( \int_{\Omega} \pd_{1}v^S \,dx  \right)^\frac{1}{2} \\
    \leq&\; C\chi \norm{\nabla(p-\bar{p})} \sqrt{G^S(t)} \leq C\chi\D(t) + C\chi G^S(t).
  \end{aligned}
\end{equation}


$\bullet$ {\bf Estimates on $\Y_i(t)\; (i=3,4,5,6)$ }:
First, by \eqref{eq:apriori-assumption} and the wave interaction estimates in Lemma \ref{lem:interaction}, we have
\begin{align*}
  &|\Y_3(t)| = \left| \int_{\Omega} a\rho \pd_{1}v^S\sdot [p'(v^S)-p'(\bar{v})](v-\bar{v}) \,dx \right|  \\
  \leq&\; C \int_{\Omega} \pd_{1}v^S\sdot |\bar{v}-v^S| |\phi| \,dx \leq C \norm{\pd_{1}v^S(\bar{v}-v^S)} \norm{\phi} \leq C\chi \delta_R \delta_S^{3/2} e^{-C\delta_S t},
\end{align*}
and so
\begin{equation*}
  \frac{M}{\delta_S} |\Y_3(t)|^2 \leq C\chi^2 \delta_R^2 \delta_S^2\sdot  e^{-C\delta_S t} .
\end{equation*}
For $\Y_4(t)$, by Holder's inequality, we have
\begin{align*}
  |\Y_4(t)| \leq C\nu\int_{\Omega} \pd_{1}a \sdot \left| h_1-\bar{h}_1 - \frac{p-\bar{p}}{\s_*} \right|  \,dx \leq C \nu^{3/2} \sqrt{\G_2(t)},
\end{align*}
and thus
\begin{equation*}
  \frac{M}{\delta_S} |\Y_4(t)|^2 \leq C \nu \G_2(t).
\end{equation*}
Recall from \eqref{eq:h-def} and \eqref{eq:hbar-def} that
\begin{equation*}
  h-\bar{h} = \psi - (2\mu+\lambda)\nabla\phi - (2\mu+\lambda)\nabla v^R,
\end{equation*}
thus, by \eqref{eq:apriori-assumption} and Lemma \ref{lem:appro-rare-property}, it follows
\begin{equation*}
  \norm{h_1-\bar{h}_1 - \frac{p-\bar{p}}{\s_*}} \lesssim \norm{\phi} +  \norm{\nabla \phi} + \norm{\psi} + \norm{\pd_{1}v^R} \lesssim \chi +\delta_R .
\end{equation*}
Hence, we yield
\begin{equation*}
  |\Y_5(t)| \leq C \norm{\sqrt{\pd_{1}a}}_{L^\infty} \sqrt{\G_2(t)} (\chi +\delta_R) \leq C \sqrt{\nu \delta_S}(\chi+\delta_R)  \sqrt{\G_2(t)},
\end{equation*}
and so
\begin{equation*}
  \frac{M}{\delta_S} |\Y_5(t)|^2 \leq C \nu (\chi+\delta_R)^2 \G_2(t).
\end{equation*}
For $\Y_6(t)$, by \eqref{eq:relative-4} and \eqref{eq:apriori-assumption}, we get
\begin{align*}
  \frac{M}{\delta_S} |\Y_6(t)|^2 &\leq \frac{C}{\delta_S} \left( \int_{\Omega} \pd_{1}a \sdot |p-\bar{p}|^2 \,dx \right)^2 + \frac{C}{\delta_S} \left( \int_{\Omega} \pd_{1}a \sdot \frac{|h_2^2+h_3^2|}{2} \,dx  \right)^2  \\
  &\leq \frac{C}{\delta_S^2} \left( \int_{\Omega} \pd_{1}v^S \sdot |p-\bar{p}|^2 \,dx \right)^2 + C \nu \norm{h-\bar{h}}^2 \G_3(t)  \\
  &\leq C \norm{p-\bar{p}}^2 G^S(t) + C \nu \norm{h-\bar{h}}^2 \G_3(t)  \\
  &\leq C \chi^2 G^S(t) + C \nu \chi^2 \G_3(t).
\end{align*}


$\bullet$ {\bf Estimates on $\B_i(t) \; (i=3,...,8)$ }:
First, it follows from \eqref{eq:F-def} that
\begin{equation*}
  F = \rho \s_*(v-v^S) + \rho(u_1-u^S_1) = \rho(\s_*\phi + \psi_1) + [\rho \s_*(\bar{v}-v^S) + \rho(\bar{u}_1-u^S_1)] .
\end{equation*}
Thus, by assumption \eqref{eq:apriori-assumption}, we get
\begin{equation}\label{eq:F leq chi}
  \begin{aligned}
    |F| \lesssim |\phi| + |\psi_1| + |\bar{v}-v^S| + |\bar{u}_1-u^S_1| \lesssim \chi +\delta_R .
  \end{aligned}
\end{equation}
For $\B_3(t)$, using \eqref{eq:relative-4} and \eqref{eq:F leq chi}, we yield
\begin{align*}
  \B_3(t) &= \int_{\Omega} F \sdot \pd_{1}a \left( Q(v|\bar{v}) + \frac{|h_1-\bar{h}_1|^2}{2} +\frac{h_2^2+h_3^2}{2} \right) dx \\
  &\lesssim \int_{\Omega} |F| \sdot \pd_{1}a \left( |p-\bar{p}|^2 + \left| h_1-\bar{h}_1 - \frac{p-\bar{p}}{\s_*} \right|^2 + \frac{h_2^2+h_3^2}{2} \right) dx \\
  &\lesssim \int_{\Omega} |F| \sdot \pd_{1}a \sdot |p-\bar{p}|^2 dx + (\chi+\delta_R)(\G_2(t)+\G_3(t)).
\end{align*}
Recall that
\begin{equation*}
  \psi_1 = h_1-\bar{h}_1 + (2\mu+\lambda)\pd_{1}\phi + (2\mu+\lambda)\pd_{1}v^R,
\end{equation*}
which together with \eqref{eq:F leq chi} leads to
\begin{align*}
  \int_{\Omega} |F| \sdot \pd_{1}a \sdot |p-\bar{p}|^2 dx \lesssim \sum_{i=1}^4 \B_{3,i}(t),
\end{align*}
where
\begin{align*}
  \B_{3,1}(t) &:= \int_{\Omega} \pd_{1}a\sdot |p-\bar{p}|^3 \,dx ,  \\
  \B_{3,2}(t) &:= \int_{\Omega} \pd_{1}a\sdot |h_1-\bar{h}_1| \sdot |p-\bar{p}|^2 \,dx ,  \\
  \B_{3,3}(t) &:= \int_{\Omega} \pd_{1}a\sdot |\pd_{1}\phi| \sdot |p-\bar{p}|^2 \,dx ,  \\
  \B_{3,4}(t) &:= \int_{\Omega} \pd_{1}a\sdot \left( |\bar{v}-v^S| + |\bar{u}_1-u^S_1| + \pd_{1}v^R \right) \sdot |p-\bar{p}|^2 \,dx .
\end{align*}
Essentially the same with \eqref{eq:4.41}, it holds
\begin{equation*}
  \B_{3,1}(t) \leq C\chi\D(t) + C\chi G^S(t).
\end{equation*}
Using Lemma \ref{lem:viscous shock} and the  Gagliardo-Nirenberg inequality \eqref{eq:G-N ineq}, we yield
\begin{align*}
  \B_{3,2}(t) &\leq C \int_{\Omega} \pd_{1}a\sdot \left| h_1-\bar{h}_1 - \frac{p-\bar{p}}{\s_*} \right| \sdot |p-\bar{p}|^2 \,dx + C\int_{\Omega} \pd_{1}a\sdot |p-\bar{p}|^3 \,dx \\
  &\leq C \norm{p-\bar{p}}_{L^\infty}^2 \sqrt{\G_2(t)} \left( \int_{\Omega} \pd_{1}a \,dx \right)^\frac{1}{2}  + C\chi\D(t) + C\chi G^S(t)  \\
  &\leq C\sqrt{\nu} \sqrt{\G_2(t)} \Big( \norm{p-\bar{p}} + \norm{\nabla^2(p-\bar{p})} \Big) \norm{\nabla(p-\bar{p})}  + C\chi\D(t) + C\chi G^S(t) \\
  &\leq C\chi \sqrt{\nu} \sqrt{\G_2(t)} \sqrt{\D(t)} + C\chi\D(t) + C\chi G^S(t) \\
  &\leq C\chi^2 \nu\G_2(t) + C\chi\D(t) + C\chi G^S(t).
\end{align*}
Noticing the fact $|\phi|\sim |p-\bar{p}|$ and that
\begin{equation*}
  \nabla(p-\bar{p}) = p'(v)\nabla(v-\bar{v}) + (p'(v) - p'(\bar{v}))\nabla\bar{v},
\end{equation*}
we have
\begin{equation}\label{eq:nab-phi-bianxing}
  |\nabla\phi| \lesssim |\nabla(p-\bar{p})| + |\pd_{1}\bar{v} |\sdot |p-\bar{p}| \lesssim |\nabla(p-\bar{p})| + (\pd_{1}v^R + \pd_{1}v^S) \sdot |p-\bar{p}|.
\end{equation}
Using \eqref{eq:nab-phi-bianxing}, we have
\begin{align*}
  \B_{3,3}(t) &\leq C\int_{\Omega} \pd_{1}a\sdot |\nabla(p-\bar{p})| \sdot |p-\bar{p}|^2 \,dx + C\int_{\Omega} \pd_{1}a\sdot (\pd_{1}v^R + \pd_{1}v^S) \sdot |p-\bar{p}|^3 \,dx \\
  &\leq C\chi^2 \delta_S^2 \sdot \D(t) + C \delta \int_{\Omega} \pd_{1}a \sdot |p-\bar{p}|^3 \,dx \leq C\delta\chi\D(t) + C\delta\chi G^S(t).
\end{align*}
Using Lemma \ref{lem:interaction} and assumption \eqref{eq:apriori-assumption}, we get
\begin{align*}
  \B_{3,4}(t) &\leq \frac{C}{\nu} \chi \norm{p-\bar{p}} \Big( \norm{\pd_{1}v^S(\bar{v}-v^S)} + \norm{\pd_{1}v^S(\bar{u}_1-u^S_1)} + \norm{\pd_{1}v^R\sdot \pd_{1}v^S}  \Big)  \\
  &\leq  C\chi^2 \delta_R \delta_S e^{-C\delta_S t}.
\end{align*}
Combining the above estimates, we obtain
\begin{equation}\label{eq:b_3-est}
  \B_3(t) \leq C\chi\D(t) + C\chi G^S(t) + C\chi \G_2(t) + C\chi^2 \delta_R \delta_S e^{-C\delta_S t}.
\end{equation}
For $\B_4(t)$, we have
\begin{equation}\label{eq:b_4-est}
  \begin{aligned}
    \B_4(t) &\leq \frac{C}{\nu} \int_{\Omega} \pd_{1}v^S \sdot |p-\bar{p}| \sdot |\nabla(p-\bar{p})|\,dx \\
    &\leq \frac{C}{\nu} \sqrt{\norm{\pd_{1}v^S}_{L^\infty}} \sqrt{G^S(t)} \sqrt{\D(t)} \leq C\nu G^S(t) + C\nu \D(t)
  \end{aligned}
\end{equation}
For $\B_5(t)$, by Lemma \ref{lem:interaction} and assumption \eqref{eq:apriori-assumption}, we yield
\begin{equation}\label{eq:b_5-est}
  \begin{aligned}
    \B_5(t) &\leq C \int_{\Omega} \pd_{1}a\sdot |p-\bar{p}|^2\sdot (\pd_{1}v^R + \pd_{1}v^S) \,dx \\
    &\leq C\nu\delta_S G^S(t) + \frac{C\chi}{\nu} \norm{p-\bar{p}} \norm{\pd_{1}v^R \sdot \pd_{1}v^S} \\
    &\leq C\nu\delta_S G^S(t) + C\chi^2 \delta_R \delta_S e^{-C\delta_S t}.
  \end{aligned}
\end{equation}
For $\B_6(t)$, using Lemma \ref{lem:appro-rare-property}, assumption \eqref{eq:apriori-assumption} and the fact $p(v|\bar{v}) \sim |p-\bar{p}|^2$, we obtain
\begin{equation}\label{eq:b_6-est}
  \begin{aligned}
    \B_6(t) &\leq C \int_{\Omega} |\nabla(p-\bar{p})| \sdot |p-\bar{p}| \sdot (\pd_{1}v^R + \pd_{1}v^S) \,dx \\
    &\leq C\delta_S \sqrt{G^S(t)} \sqrt{\D(t)} + C\norm{\nabla(p-\bar{p})} \left( \int_{\Omega} |\pd_{1}v^R|^2 |p-\bar{p}|^2 \,dx  \right)^\frac{1}{2}  \\
    &\leq C\delta_S \sqrt{G^S(t)} \sqrt{\D(t)} + C\sqrt{\delta_R}  \sqrt{\G^R(t)} \sqrt{\D(t)} \\
    &\leq C\delta_S G^S(t) + C \left( \delta_S + \sqrt{\delta_R} \right) \D(t) + C\sqrt{\delta_R}\G^R(t).
  \end{aligned}
\end{equation}
For $\B_7(t)$, by assumption \eqref{eq:apriori-assumption}, the Gagliardo-Nirenberg inequality \eqref{eq:G-N ineq}, \eqref{eq:nab-phi-bianxing} and Young's inequality, we have
\begin{equation}\label{eq:b_7-est}
  \begin{aligned}
    \B_7(t) &\leq C \norm{\pd_{11}v^R}_{L^1} \norm{\phi}_{L^\infty} \\
    &\leq C\norm{\pd_{11}v^R}_{L^1} \sdot \Big( \norm{\phi}^{\frac{1}{2}} + \norm{\nabla^2\phi}^{\frac{1}{2}} \Big) \norm{\nabla\phi}^{\frac{1}{2}} \\
    &\leq C\chi \norm{\nabla\phi}^2 + C\chi^{\frac{1}{3}} \norm{\pd_{11}v^R}_{L^1}^{\frac{4}{3}} \\
    &\leq C\chi\D(t) + C\chi\delta_R \G^R(t) + C\chi\delta_S^2 G^S(t) + C\chi^{\frac{1}{3}} \norm{\pd_{11}v^R}_{L^1}^{\frac{4}{3}}.
  \end{aligned}
\end{equation}

To estimate $\B_8(t)$, we first notice that, by the definition of $R$ in \eqref{eq:R-def}, it holds
\begin{align*}
  R &= \rho(2\mu+\lambda) (\nabla u^t - \div u \,\mathbb{I}) \nabla v - \mu\curl\curl u  \\
  &= \rho(2\mu+\lambda) (\nabla \psi^t - \div\psi \,\mathbb{I})\nabla v  - \rho(2\mu+\lambda) \pd_{1}\bar{u}_1 
  \begin{bmatrix} 
    0  \qquad  \qquad   \\
    \qquad  1  \qquad \\
    \qquad  \qquad  1
  \end{bmatrix} 
  \nabla v - \mu\curl\curl u ,
\end{align*}
where $h':=(h_2,h_3)$ and $\mathbb{I}$ denotes the $3 \times 3$ identity matrix.
Meanwhile, it follows from \eqref{eq:apriori-assumption} and Sobolev inequality that
\begin{equation}\label{eq:L3-norm}
  \begin{aligned}
    \norm{\phi}_{L^3} + \norm{\psi}_{L^3} &\lesssim \norm{\phi}^{\frac{1}{2}}\norm{\phi}_{L^6}^{\frac{1}{2}} + \norm{\psi}^{\frac{1}{2}}\norm{\psi}_{L^6}^{\frac{1}{2}} \lesssim \norm{\phi}_{H^1} + \norm{\psi}_{H^1} \lesssim \chi \,, \\
    \norm{\nabla\phi}_{L^3} + \norm{\nabla\psi}_{L^3} &\lesssim \norm{\nabla\phi}^{\frac{1}{2}}\norm{\nabla\phi}_{L^6}^{\frac{1}{2}} + \norm{\nabla\psi}^{\frac{1}{2}}\norm{\nabla\psi}_{L^6}^{\frac{1}{2}} \lesssim \norm{\nabla(\phi,\psi)}_{H^1} \lesssim \chi \,,
  \end{aligned}
\end{equation}
which together with Lemma \ref{lem:appro-rare-property} and the fact that $|h-\bar{h}| \lesssim |\psi| + |\nabla\phi| + |\pd_{1}v^R| $ leads to
\begin{equation}\label{eq:h-L^3-est}
  \norm{h-\bar{h}}_{L^3} \lesssim \chi+\delta_R .
\end{equation}
For simplicity of presentation, we divide $\B_8(t)$ as follows:
\begin{equation*}
  \B_8(t) =: \sum_{i=1}^3 \B_{8,i}(t),
\end{equation*}
where
\begin{align*}
  \B_{8,1}(t) &:= \int_{\Omega} a\rho(2\mu+\lambda) (h-\bar{h})^t (\nabla \psi^t - \div\psi \,\mathbb{I}) \nabla v \,dx,  \\
  \B_{8,2}(t) &:= -\int_{\Omega} a\rho(2\mu+\lambda) \pd_{1}\bar{u}_1\sdot h'\sdot \nabla_{x'} v \,dx = -\int_{\Omega} a\rho(2\mu+\lambda) \pd_{1}\bar{u}_1\sdot h' \sdot \frac{\nabla_{x'}p(v)}{p'(v)}  \,dx  \\
  &= -\int_{\Omega} a\rho(2\mu+\lambda) \pd_{1}\bar{u}_1\sdot h' \sdot \frac{\nabla_{x'}(p-\bar{p})}{p'(v)}  \,dx, \\
  \B_{8,3}(t) &:= -\int_{\Omega} a\mu(h-\bar{h})\sdot \curl\curl u \,dx .
\end{align*}
For $\B_{8,1}(t)$, by \eqref{eq:apriori-assumption}, \eqref{eq:nab-phi-bianxing}, \eqref{eq:h-L^3-est} and Sobolev inequality, we have
\begin{align*}
  \B_{8,1}(t) &\leq C\int_{\Omega} |h-\bar{h}| \sdot |\nabla\psi| \sdot (|\nabla\phi|+|\nabla\bar{v}|) \,dx \\
  &\leq C\int_{\Omega} |h-\bar{h}| \sdot |\nabla\psi| \sdot \big( |\nabla(p-\bar{p})| + 2|\pd_{1}\bar{v}| \big) \,dx \\
  &\leq C\int_{\Omega} |h-\bar{h}| \sdot |\nabla\psi| \sdot |\nabla(p-\bar{p})| \,dx  + C\int_{\Omega} \pd_{1}v^R \sdot |h-\bar{h}| \sdot |\nabla\psi| \,dx \\
  &\quad + C\int_{\Omega} \pd_{1}v^S \sdot |h-\bar{h}| \sdot |\nabla\psi| \,dx \\
  &\leq C \norm{h-\bar{h}}_{L^3} \norm{\nabla\psi}_{L^6} \norm{\nabla(p-\bar{p})} + C \norm{h-\bar{h}}_{L^3} \norm{\pd_{1}v^R}_{L^6} \norm{\nabla\psi}  \\
  &\quad + C\delta_S \int_{\Omega} |\nabla\psi|^2 \,dx + C\delta_S \int_{\Omega} \pd_{1}v^S\sdot |h-\bar{h}|^2 \,dx \\
  &\leq C(\chi+\delta) \left( \D(t) + \norm{\nabla\psi}_{H^1}^2 + \norm{\pd_{1}v^R}_{L^6}^2 \right) + C\delta_S \int_{\Omega} \pd_{1}v^S \sdot (h_2^2+h_3^2) \,dx \\
  &\quad + C\delta_S \int_{\Omega} \pd_{1}v^S\sdot \left| h_1-\bar{h}_1 - \frac{p-\bar{p}}{\s_*} \right|^2 \,dx + C\delta_S \int_{\Omega} \pd_{1}v^S\sdot |p-\bar{p}|^2 \,dx \\
  &\leq C(\chi+\delta) \left( \D(t) + \norm{\nabla\psi}_{H^1}^2 + \norm{\pd_{1}v^R}_{L^6}^2 \right) + C\delta_S \big[ \G_2(t) + \G_3(t) + G^S(t) \big] .
\end{align*}
Also, it follows from Lemma \ref{lem:appro-rare-property} that
\begin{equation}\label{eq:v^R-L^6-est}
  \norm{\pd_{1}v^R}_{L^6}^2 \lesssim  \delta_R^{\frac{1}{3}}\sdot (1+t)^{-2+\frac{1}{3}} = \delta_R^{\frac{1}{3}}\sdot (1+t)^{-\frac{5}{3}}.
\end{equation}
Hence, we yield
\begin{align*}
  \B_{8,1}(t) \leq\; &C(\chi+\delta) \left( \D(t) + \norm{\nabla\psi}_{H^1}^2 \right) + C\delta_S \big[ \G_2(t) + \G_3(t) + G^S(t) \big] \\ 
  &+ C(\chi+\delta) \delta_R^{\frac{1}{3}}\sdot  (1+t)^{-\frac{5}{3}}.
\end{align*}
By \eqref{eq:h-L^3-est}, it holds
\begin{align*}
  \B_{8,2}(t) &\leq C\int_{\Omega} \pd_{1}v^S \sdot |h'| \sdot |\nabla_{x'}(p-\bar{p})| \,dx + C\int_{\Omega} \pd_{1}v^R \sdot |h'| \sdot |\nabla_{x'}(p-\bar{p})| \,dx \\
  &\leq C\delta_S \int_{\Omega} \pd_{1}v^S\sdot (h_2^2+h_3^2) \,dx + C\delta_S \D(t) +  C \norm{h-\bar{h}}_{L^3} \norm{\pd_{1}v^R}_{L^6} \norm{\nabla(p-\bar{p})} \\
  &\leq C\delta_S \G_3(t) + C(\chi+\delta_S) \D(t) + C(\chi+\delta_R) \delta_R^{\frac{1}{3}}\sdot  (1+t)^{-\frac{5}{3}}.
\end{align*}
In view of the property of curl operator, we have
\begin{equation*}
  \curl u = \curl \big[ u-(2\mu+\lambda) \nabla v \big] = \curl h,
\end{equation*}
and so
\begin{equation*}
  \curl\curl u = \curl\curl h = \curl\curl(h-\bar{h}).
\end{equation*}
Thus, it holds
\begin{align*}
  \B_{8,3}(t) &= -\mu\int_{\Omega} a(h-\bar{h}) \cdot \curl\curl(h-\bar{h}) \,dx  \\
  &= -\mu\int_{\Omega} \curl\big[ a(h-\bar{h}) \big] \cdot \curl(h-\bar{h}) \,dx  \\
  &= -\mu \int_{\Omega} a |\curl(h-\bar{h})|^2 \,dx -\mu \int_{\Omega} \big[ (h-\bar{h}) \times \nabla a \big] \cdot \curl(h-\bar{h}) \,dx \\
  &= -\mu \int_{\Omega} a |\curl(h-\bar{h})|^2 \,dx -\mu \int_{\Omega} \pd_{1}a (0,-h_3,h_2) \cdot \curl(h-\bar{h}) \,dx \\
  &\leq -\frac{\mu}{2} \int_{\Omega} a |\curl(h-\bar{h})|^2 \,dx + C \int_{\Omega} |\pd_{1}a|^2 \sdot (h_2^2+h_3^2) \,dx \\
  &\leq -\frac{\mu}{2} \int_{\Omega} a |\curl(h-\bar{h})|^2 \,dx + C\nu\delta_S \G_3(t).
\end{align*}

To sum up, we obtain
\begin{equation}\label{eq:b_8-est}
  \begin{aligned}
    \B_8(t) \leq&\;  C(\chi+\delta) \left( \D(t) + \norm{\nabla\psi}_{H^1}^2 \right) + C\delta_S \big[ \G_2(t) + \G_3(t) + G^S(t) \big] \\ 
    &+ C(\chi+\delta) \delta_R^{\frac{1}{3}}\sdot  (1+t)^{-\frac{5}{3}}.
  \end{aligned}
\end{equation}


$\bullet$ {\bf Estimates on $\S(t)$ }:
By definition,
\begin{equation*}
  \S(t) = \int_{\Omega} a\rho \left[ p'(\bar{v})(v-\bar{v})\sdot vQ_1 - (h_1-\bar{h}_1)\sdot vQ_2 \right] .
\end{equation*}
Recall that
\begin{equation*}
  \begin{cases}
    Q_1 = \rho(u_1-u^R_1)\sdot \pd_{1}v^R -\rho(v-v^R) \sdot \pd_{1}u^R_1 +F\sdot \pd_{1}v^S \,, \\
    Q_2 = \rho(u_1-u^R_1)\sdot \pd_{1}u^R_1 + \rho(v-v^R)\sdot \pd_{1}p^R + F\sdot \pd_{1}h^S_1 + \left[ \pd_{1}\bar{p}-\pd_{1}p^R-\pd_{1}p^S \right],  \\
    F = \rho \s_*(v-v^S) + \rho(u_1-u^S_1) \,.
  \end{cases}
\end{equation*}
So we can devide $vQ_1$ as follows:
\begin{equation}\label{eq:Q_1 chaifen}
  \begin{aligned}
    vQ_1 &= (u_1-u^R_1)\sdot \pd_{1}v^R -(v-v^R)\sdot \pd_{1}u^R_1 + \big[ \s_*(v-v^S) + (u_1-u^S_1) \big] \sdot \pd_{1}v^S  \\
    &=: Q_1^R +Q_1^S +Q_1^I \,,
  \end{aligned}
\end{equation}
where
\begin{align*}
  Q_1^R &:= \psi_1 \pd_{1}v^R - \phi \pd_{1}u^R_1 \,, \qquad
  Q_1^S := (\s_*\phi + \psi_1) \sdot \pd_{1}v^S \,, \\
  Q_1^I &:= (\bar{u}_1-u^R_1) \sdot \pd_{1}v^R -(\bar{v}-v^R)\sdot  \pd_{1}u^R_1 + \big[ \s_*(\bar{v}-v^S) + (\bar{u}_1-u^S_1) \big] \sdot \pd_{1}v^S \,.
\end{align*}
Similarly, we let 
\begin{equation}\label{eq:Q_2 chaifen}
  \begin{aligned}
    vQ_2 &= (u_1-u^R_1)\sdot \pd_{1}u^R_1 + (v-v^R)\sdot \pd_{1}p^R + \big[ \s_*(v-v^S) + (u_1-u^S_1) \big] \sdot \pd_{1}h^S_1 \\
    &\quad + v\left[ \pd_{1}\bar{p}-\pd_{1}p^R-\pd_{1}p^S \right] 
    =: Q_2^R +Q_2^S +Q_2^I \,,
  \end{aligned}
\end{equation}
where
\begin{align*}
  Q_2^R &:= \psi_1 \pd_{1}u^R_1 + \phi \pd_{1}p^R \,, \qquad
  Q_2^S := (\s_*\phi + \psi_1) \sdot \pd_{1}h^S_1 \,, \\
  Q_2^I &:= (\bar{u}_1-u^R_1)\sdot \pd_{1}u^R_1 -(\bar{v}-v^R)\pd_{1}p^R + \big[ \s_*(\bar{v}-v^S) + (\bar{u}_1-u^S_1) \big] \pd_{1}h^S_1 \\
  &\quad + v\left[ \pd_{1}\bar{p}-\pd_{1}p^R-\pd_{1}p^S \right].
\end{align*}
With the aid of the above notations, we divide $S(t)$ as follows:
\begin{equation*}
  \S(t) = \int_{\Omega} a\rho \left[ p'(\bar{v})(v-\bar{v})\sdot vQ_1 - (h_1-\bar{h}_1)\sdot vQ_2 \right] dx =: \S_1(t) + \S_2(t) + \S_3(t),
\end{equation*}
where
\begin{align*}
  \S_1(t) &:= \int_{\Omega} a\rho \left[ p'(\bar{v})(v-\bar{v})\sdot Q_1^S - (h_1-\bar{h}_1)\sdot Q_2^S \right] dx,  \\
  \S_2(t) &:= \int_{\Omega} a\rho \left[ p'(\bar{v})(v-\bar{v})\sdot Q_1^R - (h_1-\bar{h}_1)\sdot Q_2^R \right] dx,  \\
  \S_3(t) &:= \int_{\Omega} a\rho \left[ p'(\bar{v})(v-\bar{v})\sdot Q_1^I - (h_1-\bar{h}_1)\sdot Q_2^I \right] dx.  \\
\end{align*}

For $\S_1(t)$, we have
\begin{align*}
  \S_1(t) &= \int_{\Omega} a\rho(\s_*\phi + \psi_1) \left[ p'(\bar{v})\phi\sdot \pd_{1}v^S - (h_1-\bar{h}_1)\sdot \pd_{1}h^S_1 \right] dx \\
  &= \int_{\Omega} a\rho\left[ \s_* p'(\bar{v}) \pd_{1}v^S\sdot \phi^2 - \pd_{1}h^S_1\sdot (h_1-\bar{h}_1)\psi_1 \right] dx \\
  &\quad + \int_{\Omega} a\rho\phi \left[ p'(\bar{v}) \pd_{1}v^S \sdot \psi_1 - \s_* \pd_{1}h^S_1\sdot (h_1-\bar{h}_1) \right] dx =: \S_{1,1}(t) + \S_{1,2}(t).
\end{align*}
Using the fact that $\s_* \pd_{1}h^S_1 = \pd_{1}p^S$, we get 
\begin{align*}
  \S_{1,1}(t) &= \int_{\Omega} a\rho\left[ \s_* p'(\bar{v}) \pd_{1}v^S\sdot \phi^2 - \pd_{1}h^S_1\sdot (h_1-\bar{h}_1)\psi_1 \right] dx \\
  &= \int_{\Omega} a\rho \left[ \s_* \pd_{1}p^S \phi^2 - \pd_{1}h^S_1\sdot |h_1-\bar{h}_1|^2 \right] dx  + \int_{\Omega} a\rho \s_* [p'(\bar{v})-p'(v^S)] \pd_{1}v^S\phi^2 \,dx \\
  &\quad + \int_{\Omega} a\rho \pd_{1}h^S_1\sdot (h_1-\bar{h}_1)(h_1-\bar{h}_1 - \psi_1) \,dx \\
  &\leq \int_{\Omega} a\rho \frac{\pd_{1}p^S}{\s_*} \left[ \s_*^2 \phi^2 - |h_1-\bar{h}_1|^2 \right] dx + C\delta_R G^S(t) \\ 
  &\quad + C\int_{\Omega} \pd_{1}v^S\sdot |h_1-\bar{h}_1|\sdot |\pd_{1}\phi| \,dx + C\int_{\Omega} \pd_{1}v^S\sdot |h_1-\bar{h}_1|\sdot \pd_{1}v^R \,dx .
\end{align*} 
In addition, by \eqref{eq:equi-phi}, it holds
\begin{align*}
  \left| (\s_*\phi)^2 - \frac{|p-\bar{p}|^2}{\s_*^2} \right| \lesssim |p-\bar{p}| \sdot \left| \s_*\phi + \frac{p-\bar{p}}{\s_*} \right| \lesssim |p-\bar{p}|^3,
\end{align*}
Thus, using \eqref{eq:apriori-assumption}, \eqref{eq:nab-phi-bianxing} and Lemma \ref{lem:interaction}, we obtain
\begin{align*}
  \S_{1,1}(t) &\leq C\int_{\Omega} \pd_{1}v^S \sdot \left| |h_1-\bar{h}_1|^2 - \frac{|p-\bar{p}|^2}{\s_*^2} \right| dx + C\int_{\Omega} \pd_{1}v^S\sdot |p-\bar{p}|^3 \,dx + C\delta_R G^S(t) \\ 
  &\quad + C\!\int_{\Omega} \pd_{1}v^S |h_1-\bar{h}_1| \sdot (|\nabla(p-\bar{p})|+\delta|p-\bar{p}|) \,dx + C\!\int_{\Omega} \pd_{1}v^S |h_1-\bar{h}_1|\sdot \pd_{1}v^R \,dx \\
  &\leq C\int_{\Omega} \pd_{1}v^S \left| h_1-\bar{h}_1 - \frac{p-\bar{p}}{\s_*} \right|^2 \!dx + C\int_{\Omega} \pd_{1}v^S |p-\bar{p}|\sdot \left| h_1-\bar{h}_1 - \frac{p-\bar{p}}{\s_*} \right| dx \\
  &\quad + C(\chi+\delta_R) \left( G^S(t) + \norm{\pd_{1}v^S\sdot \pd_{1}v^R} \right) + C \int_{\Omega} \pd_{1}v^S |p-\bar{p}| \sdot |\nabla(p-\bar{p})| \,dx \\
  &\quad + C \int_{\Omega} \pd_{1}v^S \left| h_1-\bar{h}_1 - \frac{p-\bar{p}}{\s_*} \right| |\nabla(p-\bar{p})| \,dx  + C\delta\int_{\Omega} \pd_{1}v^S\sdot |p-\bar{p}|^2 \,dx \\
  &\leq C\nu\G_2(t) + C\sqrt{\nu} \big( \G_2(t)+G^S(t) \big) + C(\chi+\delta_R) \left( G^S(t) + \norm{\pd_{1}v^S\sdot \pd_{1}v^R} \right) \\
  &\quad + C\delta_S \G_2(t) + C\delta_S \D(t) +C\delta G^S(t)\\
  &\leq C(\chi+\delta)G^S(t) + C(\sqrt{\nu}+\delta_S)\G_2(t) + C\delta_S \D(t) + C(\chi+\delta_R)\delta_R \delta_S^{3/2} e^{-C\delta_S t}.
\end{align*}
Similarly, by \eqref{eq:apriori-assumption}, \eqref{eq:nab-phi-bianxing} and Lemma \ref{lem:interaction}, we yield
\begin{align*}
  \S_{1,2}(t) &\leq C \int_{\Omega} \pd_{1}v^S\sdot |\bar{v}-v^S| |\phi| |\psi_1| \,dx + C \int_{\Omega} |\pd_{1}p^S|\sdot |\phi| |h_1-\bar{h}_1-\psi_1| \,dx \\
  &\leq C\chi \norm{\phi} \norm{\pd_{1}v^S (\bar{v}-v^S)} + C \int_{\Omega} \pd_{1}v^S\sdot |p-\bar{p}| (|\nabla\phi| + \pd_{1}v^R) \,dx \\
  &\leq C\chi \delta_R \delta_S^{3/2} e^{-C\delta_S t} + C\delta G^S(t) + C\delta_S \D(t).
\end{align*}
Hence, it holds
\begin{equation}\label{eq:S_1-est}
  \S_1(t) \leq C(\chi+\delta)G^S(t) + C(\sqrt{\nu}+\delta_S)\G_2(t) + C\delta_S \D(t) + C(\chi+\delta_R)\delta_R \delta_S^{3/2} e^{-C\delta_S t}.
\end{equation}
For $\S_2(t)$, using \eqref{eq:nab-phi-bianxing}, \eqref{eq:L3-norm}, \eqref{eq:v^R-L^6-est} and assumption \eqref{eq:apriori-assumption}, we have
\begin{align*}
  \S_2(t) &= \int_{\Omega} a\rho \left[ p'(\bar{v})\phi \sdot Q_1^R - (h_1-\bar{h}_1)\sdot Q_2^R \right] dx  \\
  &= - \int_{\Omega} a\pd_{1}u^R_1 \sdot \rho p'(\bar{v}) \phi^2 \,dx + \int_{\Omega} a\rho \pd_{1}v^R \sdot p'(\bar{v})\phi\psi_1 \,dx - \int_{\Omega} a\rho \pd_{1}p^R \phi\psi_1 \,dx \\
  &\quad - \int_{\Omega} a\rho \pd_{1}u^R_1 \sdot \psi_1^2 \,dx + \int_{\Omega} a\rho \left[ (2\mu+\lambda)\pd_{1}\phi + \pd_{1}v^R \right] (\psi_1 \pd_{1}u^R_1 + \phi \pd_{1}p^R) \,dx \\
  &\leq - \int_{\Omega} a\pd_{1}u^R_1 \sdot \bar{\rho} p'(\bar{v}) \phi^2 \,dx -\int_{\Omega} a\rho \pd_{1}u^R_1 \sdot \psi^2 \,dx  + C\chi \int_{\Omega} \pd_{1}u^R_1 \sdot |\phi\psi_1| \,dx \\
  &\quad + C\int_{\Omega} \left| \pd_{1}v^R \right| ^2 (|\psi|+|\phi|) \,dx + C\int_{\Omega} \pd_{1}v^R |\pd_{1}\bar{v}| \sdot |p-\bar{p}| (|\psi|+|\phi|) \,dx \\
  &\quad + C\chi\G^R(t) + C\int_{\Omega} \pd_{1}v^R \sdot |\nabla(p-\bar{p})| (|\psi|+|\phi|) \,dx  \\
  &\leq - \int_{\Omega} a\pd_{1}u^R_1 \sdot \bar{\rho} p'(\bar{v}) \phi^2 \,dx -\int_{\Omega} a\rho \pd_{1}u^R_1 \sdot \psi^2 \,dx + C\chi \int_{\Omega} \pd_{1}u^R_1 \sdot \psi_1^2 \,dx \\
  &\quad + C\chi \norm{\pd_{1}v^R}_{L^4}^2 + C\chi \norm{\pd_{1}v^R\sdot \pd_{1}v^S} + C\chi\G^R(t) + C\chi \norm{\pd_{1}v^R}_{L^6} \norm{\nabla(p-\bar{p})} \\
  &\leq - \int_{\Omega} a\pd_{1}u^R_1 \sdot \bar{\rho} p'(\bar{v}) \phi^2 \,dx -\int_{\Omega} a\rho \pd_{1}u^R_1 \sdot \psi^2 \,dx + C\chi \int_{\Omega} \pd_{1}u^R_1 \sdot \psi_1^2 \,dx \\
  &\quad + C\chi \norm{\pd_{1}v^R}_{L^4}^2 + C\chi\G^R(t) + C\chi\D(t) + C\chi \delta_R \delta_S^{2/3} e^{-C\delta_S t} + C\chi \delta_R^{\frac{1}{3}}\sdot (1+t)^{-\frac{5}{3}} .
\end{align*}
Notice that
\begin{equation*}
  \G^R(t) = \int_{\Omega} a \pd_{1}u^R_1 \sdot p(v|\bar{v}) \,dx \geq \int_{\Omega} a \pd_{1}u^R_1 \sdot \frac{p''(\bar{v})}{2} \phi^2 \,dx + C\chi\G^R(t),
\end{equation*}
and that
\begin{equation*}
  -\bar{\rho} p'(\bar{v}) - \frac{p''(\bar{v})}{2} = \gamma \bar{v}^{\gamma-2} - \frac{\gamma(\gamma+1)}{2} \bar{v}^{\gamma-2} = - \frac{\gamma-1}{\gamma+1} \sdot \frac{p''(\bar{v})}{2},
\end{equation*}
so it follows
\begin{equation*}
  - \int_{\Omega} a\pd_{1}u^R_1 \sdot \bar{\rho} p'(\bar{v}) \phi^2 \,dx - \G^R(t) \leq -\frac{\gamma-1}{\gamma+1} \G^R(t) + C\chi\G^R(t) .
\end{equation*}
Therefore, choosing $\chi$ sufficiently small, we acquire
\begin{align*}
  \S_2(t) - \G^R(t) &\leq -\frac{\gamma-1}{\gamma+1} \G^R(t) +  C\chi\G^R(t) + C\chi\D(t) + C\chi \delta_R \delta_S^{2/3} e^{-C\delta_S t} \\
  &\quad + C\chi \delta_R^{\frac{1}{3}}\sdot (1+t)^{-\frac{5}{3}} + C\chi \norm{\pd_{1}v^R}_{L^4}^2 .
\end{align*}
For $\S_3(t)$, $\pd_{1}\bar{p}-\pd_{1}p^R-\pd_{1}p^S $ is the only nontrivial term in $Q_2^I$. Notice that
\begin{align*}
  |\pd_{1}\bar{p}-\pd_{1}p^R-\pd_{1}p^S| &= |p'(\bar{v})(\pd_{1}v^R + \pd_{1}v^S) - p'(v^R)\pd_{1}v^R - p'(v^S)\pd_{1}v^S| \\
  &\lesssim \pd_{1}v^R\sdot |\bar{v}-v^R| + \pd_{1}v^S\sdot |\bar{v}-v^S|.
\end{align*}
According to Lemma \ref{lem:interaction}, it holds
\begin{equation*}
  \norm{Q^I_i} \lesssim \delta_R \delta_S \sdot e^{-C\delta_S t}, \quad i=1,2.
\end{equation*}
Also, by assumption \eqref{eq:apriori-assumption}, we have
\begin{equation*}
  \norm{h-\bar{h}} \lesssim \norm{\psi} + \norm{\nabla\phi} + \norm{\pd_{1}v^R} \lesssim \chi+\delta_R.
\end{equation*}
Thus, using Lemma \ref{lem:interaction}, we obtain
\begin{equation}\label{eq:S_3-est}
  \begin{aligned}
    \S_3(t) &\leq C \norm{\phi} \norm{Q^I_1} + C \norm{h-\bar{h}} \norm{Q^I_2}  \\
    &\leq C(\chi+\delta_R) \delta_R \delta_S \sdot e^{-C\delta_S t}.
  \end{aligned}
\end{equation}
Combining the above estimates, we arrive at
\begin{equation}\label{eq:S-est}
  \begin{aligned}
    &\quad \S(t) - \G^R(t) \\
    &\leq -\frac{\gamma-1}{\gamma+1} \G^R(t) + C(\sqrt{\nu}+\delta_S)\G_2(t) + C\chi\G^R(t) + C(\chi+\delta_S)\D(t) + C(\chi+\delta)G^S(t)  \\
    &\quad + C(\chi+\delta_R) \delta_R \delta_S \sdot e^{-C\delta_S t} + C(\chi+\delta)\delta_R^{\frac{1}{3}}\sdot (1+t)^{-\frac{5}{3}} + C\chi \norm{\pd_{1}v^R}_{L^4}^2 .
  \end{aligned}
\end{equation}

Finally, using Lemma \ref{lem:4.1} and Lemma \ref{lem:key lemma} and choosing $\nu,\delta_R, \delta_S, \chi$ sufficiently small, we yield
\begin{equation}\label{eq:zhongjian}
  \begin{aligned}
    &\frac{d}{dt} \int_{\Omega} a\rho \left[ Q(v|\bar{v}) + \frac{|h-\bar{h}|^2}{2} \right] dx + \delta_S |\dot{\X}(t)|^2 \\
    &\quad + \G_2(t) + \G_3(t) + \G^R(t) + \D(t) + G^S(t) \\
    \lesssim &\; (\chi+\delta) \norm{\nabla\psi}_{H^1}^2 + \delta_R \delta_S \sdot e^{-C\delta_S t} + \delta_R^{\frac{1}{3}}\sdot (1+t)^{-\frac{5}{3}} + \norm{\pd_{1}v^R}_{L^4}^2 + \norm{\pd_{11}v^R}_{L^1}^{\frac{4}{3}}.
  \end{aligned}
\end{equation}
According to Lemma \ref{lem:appro-rare-property}, it holds
\begin{equation*}
  \norm{\pd_{11}v^R}_{L^1} \lesssim 
  \begin{cases}
    \delta_R,        &1+t \leq \frac{1}{\delta_R}, \\
    \frac{1}{1+t},   &1+t \geq \frac{1}{\delta_R},
  \end{cases}
\end{equation*}
and
\begin{equation*}
  \norm{\pd_{1}v^R}_{L^4} \lesssim  
  \begin{cases}
    \delta_R,        &1+t \leq \frac{1}{\delta_R}, \\
    \displaystyle \frac{\delta_R^{\frac{1}{4}}}{1+t},   &1+t \geq \frac{1}{\delta_R} .
  \end{cases}
\end{equation*}
So direct computations give 
\begin{equation}\label{eq:rare-time jifen}
  \int_0^\infty \norm{\pd_{11}v^R}_{L^1}^{\frac{4}{3}} d\tau \lesssim  \delta_R^{\frac{1}{3}} \,, \qquad \int_0^\infty \norm{\pd_{1}v^R}_{L^4}^2 d\tau \lesssim  \delta_R.
\end{equation}
Notice that
\begin{equation*}
  Q(v|\bar{v}) \sim |v-\bar{v}|^2,
\end{equation*}
and that
\begin{equation*}
  \G_2(t) \sim G_2(t), \quad \G_3(t) \sim G_3(t), \quad \G^R(t) \sim G^R(t), \quad \D(T) \sim D(t).
\end{equation*}
Hence, integrating \eqref{eq:zhongjian} over $[0,t]$, we acquire the desired inequality \eqref{eq:phi,h-L^2-est}, which completes the proof of Lemma \ref{lem:phi,h-estimate}.
\hfill  $\blacksquare$


\subsection{Wave interaction estimates}\label{sec:interaction}

According to \eqref{eq:Q_1-def} and \eqref{eq:F-def}, it holds
\begin{align*}
    vQ_1 &= (u_1-u^R_1)\pd_{1}v^R -(v-v^R)\pd_{1}u^R_1 + \big[ \s_*(v-v^S) + (u_1-u^S_1) \big] \pd_{1}v^S  \\
    &=: Q_1^R +Q_1^S +Q_1^I \,,
\end{align*}
where
\begin{align*}
  Q_1^R &:= \psi_1 \pd_{1}v^R - \phi \pd_{1}u^R_1 \,,  \qquad  Q_1^S := (\s_*\phi + \psi_1) \pd_{1}v^S \,, \\
  Q_1^I &:= (\bar{u}_1-u^R_1)\pd_{1}v^R -(\bar{v}-v^R)\pd_{1}u^R_1 + \big[ \s_*(\bar{v}-v^S) + (\bar{u}_1-u^S_1) \big] \pd_{1}v^S \,.
\end{align*}
Note that the above $Q^I_1$ is a nonlinear term reflecting the interaction of the approximate rarefaction and viscous shock, which can be estimated as follows.

\begin{lemma}\label{lem:interaction}
  Under the assumptions of Proposition \ref{prop:apriori-estimate}, there exists constant $C>0$ (independent of $\nu,\delta,\chi$ and $T$) such that for all $t \in [0,T]$, it holds
  \begin{align}
    &\norm{\pd_{1}v^R(\bar{v}-v^R)} + \norm{\pd_{1}v^R(\bar{u}_1-u^R_1)} \leq C\delta_R \delta_S e^{-C\delta_S t}, \label{eq:inter_1} \\
    &\norm{\pd_{1}v^S(\bar{v}-v^S)} + \norm{\pd_{1}v^S(\bar{u}_1-u^S_1)} + \norm{\pd_{1}v^R\sdot \pd_{1}v^S} \leq C\delta_R \delta_S^{3/2} e^{-C\delta_S t}, \label{eq:inter_2}
  \end{align}
\end{lemma}

\begin{proof}
  Notice that the Navier-Stokes system and Euler system are of Galilean invariance, so we can assume $u_{1m}=0$ without loss of generality. It follows from \eqref{eq:Lax-E condition} that $\lambda_{1m} < \s < u_{1m}=0$. 
  Since $\bar{v}=v^R+v^S-v_m$, it holds
  \begin{equation*}
    \norm{\pd_{1}v^R(\bar{v}-v^R)}^2 = \int_{\Omega} |\pd_{1}v^R|^2 |v^S-v_m|^2 \,dx = \int_{\R} |\pd_{1}v^R|^2 |v^S-v_m|^2 \,dx_1 .
  \end{equation*}
  To estimate the above integral, we introduce two time-dependent sets:
  \begin{align*}
    E_1(t)&:= \left\{ x_1 \in \R: x_1 > \frac{\lambda_{1m} +\s}{2}(1+t) \right\}, \\
    E_2(t)&:= \left\{ x_1 \in \R: x_1 \leq \frac{\lambda_{1m} +\s}{2}(1+t) \right\}.
  \end{align*}
  First, by assumption \eqref{eq:apriori-assumption} and Sobolev inequality, it follows
  \begin{equation}\label{eq:4.9}
    \norm{p-\bar{p}}_{L^\infty} \leq C \norm{v-\bar{v}}_{L^\infty} \leq C\chi,
  \end{equation}
  In addition, by the definition of $\X(t)$ \eqref{eq:X(t)-def}, we have
  \begin{equation}\label{eq:4.10}
    |\dot{\X}(t)| \leq \frac{C}{\delta_S} \norm{v-\bar{v}}_{L^\infty} \int_{\Omega} \pd_{1}v^S \,dx \leq C \norm{v-\bar{v}}_{L^\infty} \leq C\chi, \quad t \in [0,T],
  \end{equation}
  which together with $\X(0)=0$ leads to
  \begin{equation*}
    |\X(t)| \leq C\chi \,t, \quad t \in [0,T] .
  \end{equation*}
  According to the Lax entropy conditions \eqref{eq:Lax-E condition}, it holds $\s-\lambda_{1m}>0$. Thus choosing $\chi$ sufficiently small, we have
  \begin{equation*}
    |\X(t)| \leq C\chi \,t < \frac{\s-\lambda_{1m}}{4}t , \quad t \in [0,T] .
  \end{equation*}
  Hence, for each fixed $t\in[0,T]$, if $x_1 \in E_2$, then
  \begin{align*}
    &x_1 -\s t - \X(t) \leq \frac{\lambda_{1m}+\s}{2}(1+t) -\s t - \X(t) \\
    \leq\; &\frac{\lambda_{1m}+\s}{2} + \frac{1}{2}(\lambda_{1m}-\s)t -\X(t) < -\frac{1}{4}(\s-\lambda_{1m})t <0,
  \end{align*}
  Also, according to Lemma \ref{lem:viscous shock} and the definition of $v^S$ \eqref{eq:v^S-def}, it holds
  \begin{equation}
    \begin{aligned}
      |v^S-v_m| &\leq C\delta_S  e^{-C\delta_S |x_1 -\s t - \X(t)|}  \\
      &\leq C\delta_S  e^\frac{-C\delta_S |x_1 -\s t - \X(t)|}{2} \sdot e^{-C\delta_S t}, \quad x_1 \in E_2,
    \end{aligned}
  \end{equation}
  Therefore, using Lemma \ref{lem:appro-rare-property}, we obtain
  \begin{equation}
    \int_{E_2} |\pd_{1}v^R|^2 |v^S-v_m|^2 \,dx_1 \leq C\delta_S^2 e^{-2C\delta_S t} \int_{\R} |\pd_{1}v^R|^2 \,dx_1 \leq C\delta_R^2\delta_S^2 e^{-2C\delta_S t} .
  \end{equation}

  Secondly, for each fixed $t\in[0,T]$, if $x_1 \in E_1$, then
  \begin{equation*}
    x_1-\lambda_{1m}(1+t) = x_1 - \frac{\lambda_{1m} +\s}{2}(1+t) + \frac{\s-\lambda_{1m}}{2}(1+t) >0 .
  \end{equation*} 
  Notice that
  \begin{equation*}
    e^{-2|x_1-\lambda_{1m}(1+t)|} = e^{-2[x_1 - \frac{\lambda_{1m} +\s}{2}(1+t) + \frac{\s-\lambda_{1m}}{2}(1+t)]} = e^{[x_1 - \frac{\lambda_{1m} +\s}{2}(1+t)]} \sdot e^{[\frac{\s-\lambda_{1m}}{2}(1+t)]}
  \end{equation*}
  Thus, by Lemma \ref{lem:appro-rare-property}, it holds
  \begin{equation}
    \begin{aligned}
      \int_{E_1} |\pd_{1}v^R|^2 |v^S-v_m|^2 \,dx_1 
      &\leq C\delta_R^2\delta_S^2 \sdot e^{[\frac{\s-\lambda_{1m}}{2}(1+t)]}  \int_{E_1} e^{[x_1 - \frac{\lambda_{1m} +\s}{2}(1+t)]} \,dx_1,  \\
      &\leq C\delta_R^2\delta_S^2 \sdot e^{-Ct}.
    \end{aligned}
  \end{equation}
  To sum up, we have
  \begin{equation*}
    \norm{\pd_{1}v^R(\bar{v}-v^R)} \leq C\delta_R \delta_S \sdot e^{-C\delta_S t} .
  \end{equation*}
  Similarly, we can acquire the same estimate for $\norm{\pd_{1}v^R(\bar{u}_1-u^R_1)}$ and then verify \eqref{eq:inter_1}.

  Similar to the proof of \eqref{eq:inter_1}, in order to prove  \eqref{eq:inter_2}, we only need to verify
  \begin{equation*}
    \norm{\pd_{1}v^S(\bar{v}-v^S)} \leq C\delta_R \delta_S^{3/2} e^{-C\delta_S t}.
  \end{equation*}
  Notice that
  \begin{align*}
    \norm{\pd_{1}v^S(\bar{v}-v^S)}^2 &= \int_{\Omega} |\pd_{1}v^R|^2 |v^R-v_m|^2 \,dx = \int_{\R} |\pd_{1}v^S|^2 |v^R-v_m|^2 \,dx_1  \\
    &= \int_{E_1} |\pd_{1}v^S|^2 |v^R-v_m|^2 \,dx_1 + \int_{E_2} |\pd_{1}v^S|^2 |v^R-v_m|^2 \,dx_1 .
  \end{align*}
  Thus for each fixed $t\in[0,T]$, if $x_1 \in E_2$, then by Lemma \ref{lem:viscous shock} and the definition of $v^S$ \eqref{eq:v^S-def}, it follows
  \begin{equation}
    \begin{aligned}
      |\pd_{1}v^S| &\leq C\delta_S^2 \sdot e^{-C\delta_S |x_1 -\s t - \X(t)|}  \\
      &\leq C\delta_S^2 \sdot e^\frac{-C\delta_S |x_1 -\s t - \X(t)|}{2} \sdot e^{-C\delta_S t}, \quad x_1 \in E_2 .
    \end{aligned}
  \end{equation}
  Hence, 
  \begin{align*}
    \int_{E_2} |\pd_{1}v^S|^2 |v^R-v_m|^2 \,dx_1 &\leq C\delta_R^2 \delta_S^2 \sdot e^{-C\delta_S t} \int_{E_2} \pd_{1}v^S \,dx_1 \\
    &\leq C\delta_R^2 \delta_S^3 \sdot e^{-C\delta_S t} .
  \end{align*}
  Meanwhile, 
  \begin{align*}
    \int_{E_1} |\pd_{1}v^S|^2 |v^R-v_m|^2 \,dx_1 
      &\leq C\delta_R^2\delta_S^4 \sdot e^{[\frac{\s-\lambda_{1m}}{2}(1+t)]}  \int_{E_1} e^{[x_1 - \frac{\lambda_{1m} +\s}{2}(1+t)]} \,dx_1,  \\
      &\leq C\delta_R^2\delta_S^4 \sdot e^{-Ct} .
  \end{align*}
  Therefore, we have
  \begin{equation*}
    \norm{\pd_{1}v^S(\bar{v}-v^S)} \leq C\delta_R \delta_S^{3/2} e^{-C\delta_S t}.
  \end{equation*}
  This completes the proof of Lemma \ref{lem:interaction}.
\end{proof}


\subsection{Estimates on $\norm{u-\bar{u}}$ and $\norm{v-\bar{v}}_{H^1}$ }  
From now on, we will consider the energy estimates with respect to $(v-\bar{v}, u-\bar{u})$ .
In this subsection, we first give the zero-th order estimates.

\begin{lemma}\label{lem:psi-0-estimate}
  Under the assumptions of Proposition \ref{prop:apriori-estimate}, there exists constant $C>0$ (independent of $\nu,\delta,\chi$ and $T$)  such that for all $t \in [0,T]$, it holds
  \begin{equation}\label{eq:psi-L^2-est}
    \begin{aligned}
      &\norm{(v-\bar{v})(t)}_{H^1}^2  + \norm{(u-\bar{u})(t)}^2 + \delta_S \int_0^t |\dot{\X}(\tau)|^2 d\tau  \\
      &+ \int_0^t G_2(\tau) + G_3(\tau) + \mathcal{G}^R(\tau) + G^S(\tau) + D(\tau) + \norm{\nabla(u-\bar{u})}^2 \,d\tau \\
      \leq&\; C \left( \norm{v_0-\bar{v}(0,\cdot)}_{H^1}^2 + \norm{u_0-\bar{u}(0,\cdot)}^2 \right)  +C(\delta+\chi)\int_0^t \norm{\nabla^2(u-\bar{u})}^2 d\tau + C\delta_R^{\frac{1}{3}} ,
    \end{aligned}
  \end{equation}
  where $G^S(t), G^R(t), D(t)$ are as in \eqref{eq:new good-def} and
  \begin{equation*}
    \mathcal{G}^R(t) := \frac{\gamma-1}{\gamma+1}G^R(t) + \int_{\Omega} \rho \pd_{1}u^R_1\sdot \psi_1^2 \,dx = \int_{\Omega} \frac{\gamma-1}{\gamma+1} \pd_{1}u^R_1 \sdot p(v|\bar{v}) + \rho \pd_{1}u^R_1\sdot \psi_1^2 \,dx .
  \end{equation*}
\end{lemma}

\begin{proof}
  Since $\bar{v}=v^S+v^R-v_m$, $\bar{u}=u^S-u^R-u_m$, direct computations give
  \begin{equation}\label{eq:vbar,ubar-eqn}
    \begin{cases}
      \rmd{\bar{v}} -\div\bar{u} + \rho\dot{\X}(t) \pd_{1}v^S = Q_1, \\
      \rmd{\bar{u}} + \nabla\bar{p} + \rho\dot{\X}(t) \pd_{1}u^S = (2\mu+\lambda) \Delta(\bar{u}-u^R) + Q_3\sdot e_1.
    \end{cases}
  \end{equation}
  where $Q_1$ is as in \eqref{eq:Q_1-def} and $e_1=(1,0,0)^t$. $Q_3$ is defined as follows:
  \begin{equation}\label{eq:Q_3-def}
    Q_3 := \rho(v-v^R)\pd_{1}p^R + \rho(u_1-u^R_1)\pd_{1}u^R_1 + F\sdot \pd_{1}u^S_1 + \pd_{1}(\bar{p}-p^S-p^R) ,
  \end{equation}
  where $F$ is as defined in \eqref{eq:F-def} .
  Subtracting \eqref{eq:vbar,ubar-eqn} from \eqref{eq:3D-NS-v,u}, we acquire the perturbed system for $(\phi,\psi)= (v-\bar{v},u-\bar{u})$ which read as
  \begin{equation}\label{eq:psi,phi-eqn}
    \left\{
    \begin{aligned}
      &\rmd{\phi} -\div\psi - \rho\dot{\X}(t) \pd_{1}v^S = -Q_1, \\
      &\rmd{\psi} + \nabla(p-\bar{p}) - \rho\dot{\X}(t) \pd_{1}u^S = \mu\Delta\psi + (\mu+\lambda)\nabla\div\psi \\ 
      & \hspace{18em} + (2\mu+\lambda) \Delta u^R - Q_3\sdot e_1.
    \end{aligned} 
    \right. 
  \end{equation}
  Multiplying \eqref{eq:psi,phi-eqn}$_1$ by $-(p-\bar{p})$, we have
  \begin{equation}\label{eq:0jie-1}
    \begin{aligned}
      &\pd_{t}\big[\rho Q(v|\bar{v})\big] + \div \big[ \rho Q(v|\bar{v}) u \big] + (p-\bar{p})\div\psi  \\
      =&\; -\rho \dot{\X}(t) \pd_{1}v^S\sdot (p-\bar{p}) 
      - p(v|\bar{v}) \big[ \rmd{\bar{v}} \big] + Q_1(p-\bar{p}) \\
      =&\; -\rho \dot{\X}(t) \pd_{1}v^S \sdot p'(\bar{v})(v-\bar{v}) - p(v|\bar{v})\sdot \div\bar{u} + p'(\bar{v})(v-\bar{v})Q_1.
    \end{aligned}
  \end{equation}
  Multiplying \eqref{eq:psi,phi-eqn}$_2$ by $\psi=u-\bar{u}$, we have
  \begin{equation}\label{eq:0jie-2}
    \begin{aligned}
      &\pd_{t}\left[ \rho \frac{|\psi|^2}{2} \right] + \div \left[ \rho \frac{|\psi|^2}{2} u \right] + \div \big[ (p-\bar{p})\psi \big] - (p-\bar{p})\div\psi \\
      =&\; \rho\dot{\X}(t) \pd_{1}u^S_1\sdot \psi_1 + (2\mu+\lambda) \pd_{11}u^R_1\sdot \psi_1 -Q_3\sdot \psi_1 + \mu\psi\Delta\psi + (\mu+\lambda)\nabla\div\psi .
    \end{aligned}
  \end{equation}
  Adding \eqref{eq:0jie-1} and \eqref{eq:0jie-2} together and integrating the resultant equation over $\Omega$, we yield
  \begin{equation}\label{eq:0-weifen}
    \frac{d}{dt} \int_{\Omega} \rho Q(v|\bar{v}) + \rho \frac{|\psi|^2}{2} \,dx + G^R(t) + D_1(t) = \dot{\X}(t) \mathcal{Y}(t) + \sum_{i=1}^3 I_i(t),
  \end{equation}
  where
  \begin{align*}
    D_1(t) &:= \int_{\Omega} \mu|\nabla\psi|^2 + (\mu+\lambda) |\div \psi|^2 \,dx , \\
    \mathcal{Y}(t) &:= -\int_{\Omega} \rho \pd_{1}v^S\sdot p'(\bar{v})\phi \,dx  + \int_{\Omega} \rho \pd_{1}u^S_1\sdot \psi_1 \,dx =: \mathcal{Y}_1(t) + \mathcal{Y}_2(t), \\
    I_1(t) &:= \int_{\Omega} -\pd_{1}u^S_1 \sdot p(v|\bar{v})  \,dx , \\
    I_2(t) &:= \int_{\Omega} p'(\bar{v})\phi \sdot Q_1 - Q_3\sdot \psi_1 \,dx  , \\
    I_3(t) &:= \int_{\Omega} (2\mu+\lambda) \pd_{11}u^R_1\sdot \psi_1 \,dx .
  \end{align*}
  For $I_1(t)$, it is easy to see that
  \begin{equation*}
    I_1(t) \leq  C G^S(t).
  \end{equation*}
  The treatment for $I_2(t)$ is similar to that for $\S(t)$. We first divide $vQ_3$ as follows, where $Q_3$ is defined in \eqref{eq:Q_3-def}.
  \begin{equation}\label{eq:Q_3 chaifen}
    \begin{aligned}
      vQ_3 &= (u_1-u^R_1)\sdot \pd_{1}u^R_1 + (v-v^R)\sdot \pd_{1}p^R + \big[ \s_*(v-v^S) + (u_1-u^S_1) \big] \sdot \pd_{1}u^S_1 \\
    &\quad + v\left[ \pd_{1}\bar{p}-\pd_{1}p^R-\pd_{1}p^S \right] 
    =: Q_3^R +Q_3^S +Q_3^I \,,
    \end{aligned}
  \end{equation}
  where
  \begin{align*}
    Q_3^R &:= \psi_1 \pd_{1}u^R_1 + \phi \pd_{1}p^R \,, \qquad
    Q_3^S := (\s_*\phi + \psi_1) \sdot \pd_{1}u^S_1 \,, \\
    Q_3^I &:= (\bar{u}_1-u^R_1)\sdot \pd_{1}u^R_1 -(\bar{v}-v^R)\pd_{1}p^R + \big[ \s_*(\bar{v}-v^S) + (\bar{u}_1-u^S_1) \big] \pd_{1}u^S_1 \\
    &\quad + v\left[ \pd_{1}\bar{p}-\pd_{1}p^R-\pd_{1}p^S \right].
  \end{align*}
  Thus,
  \begin{align*}
    I_2(t) =: \sum_{i=1}^3 I_{2,i}(t) 
    =& \int_{\Omega} \rho \left[ p'(\bar{v})\phi\sdot Q_1^S - \psi_1\sdot Q_3^S \right] dx  \\
    & + \int_{\Omega} \rho \left[ p'(\bar{v})\phi\sdot Q_1^R - \psi_1\sdot Q_3^R \right] dx  \\ 
    & + \int_{\Omega} \rho \left[ p'(\bar{v})\phi\sdot Q_1^I - \psi_1\sdot Q_3^I \right] dx. 
  \end{align*}
  For $I_{2,1}(t)$, by Cauchy's inequality, it holds
  \begin{align*}
    I_{2,1}(t) \lesssim \int_{\Omega} \pd_{1}v^S (\phi^2+\psi_1^2) \,dx \lesssim G^S(t) + \int_{\Omega} \pd_{1}v^S\sdot \psi_1^2 \,dx .
  \end{align*}
  Notice that
  \begin{align*}
    |\psi_1| &\lesssim \left| h_1-\bar{h}_1 - \frac{p-\bar{p}}{\s_*} \right|  + |p-\bar{p}| + |\nabla\phi| + \pd_{1}v^R  \\
    &\lesssim \left| h_1-\bar{h}_1 - \frac{p-\bar{p}}{\s_*} \right|  + |p-\bar{p}| + |\nabla(p-\bar{p})| + \pd_{1}v^R ,
  \end{align*}
  which together with Lemma \ref{lem:interaction} leads to
  \begin{align*}
    I_{2,1}(t) \leq CG^S(t) + C\nu G_2(t) + C\delta_S D(t) + C\chi^2 \delta_R \delta_S^{\frac{3}{2}} e^{-C\delta_S t}.
  \end{align*}
  Similar to the estimate for $\S_2(t)$, we have
  \begin{align*}
    &I_{2,2}(t) - G^R(t)  = \int_{\Omega} \rho \left[ p'(\bar{v})\phi\sdot Q_1^R - \psi_1\sdot Q_3^R \right] - \pd_{1}u^R_1 \sdot p(v|\bar{v}) \,dx \\
    =& \int_{\Omega} \rho \left[ -p'(\bar{v})\pd_{1}u^R_1 \sdot \phi^2 + p'(\bar{v})\pd_{1}v^R\sdot \phi\psi_1 -\pd_{1}p^R \sdot \phi\psi_1 - \pd_{1}u^R_1 \sdot \psi_1^2  \right] - \pd_{1}u^R_1 \sdot p(v|\bar{v})  \,dx \\
    \leq& -\frac{\gamma-1}{\gamma+1}G^R(t) - \int_{\Omega} \rho \pd_{1}u^R_1\sdot \psi_1^2 \,dx + C(\chi+\delta)G^R(t) + C \chi^2 \norm{\pd_{1}v^R (\bar{v}-v^R)} \\
    \leq& -\mathcal{G}^R(t) + C(\chi+\delta)G^R(t) + C \chi^2 \delta_R \delta_R \sdot e^{-C\delta_S t} .
  \end{align*}
  The estimate for $I_{2,3}(t)$ is exactly the same with that for $\S_3(t)$, which reads as 
  \begin{equation*}
    I_{2,3}(t) \leq C(\chi+\delta_R) \delta_R \delta_S \sdot e^{-C\delta_S t}.
  \end{equation*}
  Hence, we obtain 
  \begin{equation}
    \begin{aligned}
      I_2(t) - G^R(t) &\leq -\mathcal{G}^R(t) + C(\chi+\delta)G^R(t) + C \delta_R \delta_R \sdot e^{-C\delta_S t} \\
      &\quad +CG^S(t) + C\nu G_2(t) + C\delta_S D(t).
    \end{aligned}
  \end{equation}
  For $I_3(t)$, by assumption \eqref{eq:apriori-assumption}, the Gagliardo-Nirenberg inequality \eqref{eq:G-N ineq} and \eqref{eq:nab-phi-bianxing}, we have
  \begin{equation}
    \begin{aligned}
      I_3(t) &\leq C \norm{\pd_{11}u^R_1}_{L^1} \norm{\psi_1}_{L^\infty} \\
      &\leq C\norm{\pd_{11}u^R_1}_{L^1} \sdot \Big( \norm{\psi_1}^{\frac{1}{2}} + \norm{\nabla^2\psi_1}^{\frac{1}{2}} \Big) \norm{\nabla\psi_1}^{\frac{1}{2}} \\
      &\leq C\chi D_1(t) + C\chi^{\frac{1}{3}} \norm{\pd_{11}v^R}_{L^1}^{\frac{4}{3}} .
    \end{aligned}
  \end{equation}
  Since $|p-\bar{p}|\sim |\phi|$, it follows from Holder inequality that
  \begin{equation*}
    |\mathcal{Y}_1(t)| \leq \left( \int_{\Omega} \pd_{1}v^S \,dx \right)^\frac{1}{2} \left( \int_{\Omega} \pd_{1}v^S\sdot \phi^2 \,dx  \right) \leq C\nu \sqrt{G^S(t)}.
  \end{equation*}
  Recall that
  \begin{equation*}
    \phi_1 = h_1-\bar{h}_1 + (2\mu+\lambda) \pd_{1}\phi + (2\mu+\lambda) \pd_{1}v^R,
  \end{equation*}
  and that
  \begin{equation*}
    |\nabla\phi| \lesssim |\nabla(p-\bar{p})| + |\pd_{1}\bar{v}| |p-\bar{p}| .
  \end{equation*}
  Thus, using Lemma \ref{lem:interaction}, we have
  \begin{align*}
    |\mathcal{Y}_2(t)| &\leq C \int_{\Omega} \pd_{1}v^S \left( \left| h_1-\bar{h}_1 - \frac{p-\bar{p}}{\s_*} \right| + |p-\bar{p}| + |\nabla(p-\bar{p})| + \pd_{1}v^R \right) dx \\
    &\leq C\nu^{\frac{3}{4}} \sqrt{G_2(t)} + C\nu \sqrt{G^S(t)} + \delta_S \sqrt{D(t)} + C\delta_R \delta_S \sdot e^{-C\delta_S t} .
  \end{align*}
  Therefore, we obtain
  \begin{align*}
    \dot{\X}(t) \mathcal{Y}(t) &\leq  |\dot{\X}(t)| \left( |\mathcal{Y}_1(t)| + |\mathcal{Y}_1(t)| \right)  \\
    &\leq \frac{\delta_S}{4} |\dot{\X}(t)|^2 + C\nu G_2(t) + C_1 G^S(t) + C\delta_S D(t) +  C\delta_R^2 \delta_S \sdot e^{-C\delta_S t}.
  \end{align*}
  Integrating \eqref{eq:0-weifen} over $[0,t]$ for any $t \in [0,T]$, using the above estimates and \eqref{eq:rare-time jifen}, we can find some constant $C_2>0$ independent of $t$ such that
  \begin{align*}
    &\int_{\Omega} \rho Q(v|\bar{v}) +  \frac{\rho|\psi|^2}{2} \,dx + \int_0^t \frac{\gamma-1}{\gamma+1}G^R(\tau) + D_1(\tau)  d\tau \\
    \leq&\; C \int_{\Omega} |v_0-\bar{v}(0,x)|^2 + |u_0-\bar{u}(0,x)|^2 \,dx + \int_0^t \frac{\delta_S}{4} |\dot{\X}(\tau)|^2 + C_1 G^S(\tau) d\tau   \\
    &\; + \int_0^t C\nu G_2(\tau) + C\delta_S D(\tau) + C(\chi+\delta)G^R(\tau) + C\chi D_1(\tau) d\tau + C\delta_R^{\frac{1}{3}} .
  \end{align*}
  Multiplying the above inequality by $\frac{1}{2\max \{ 1,C_2 \}}$, adding the resultant expression to \eqref{eq:phi,h-L^2-est}, and choosing $\nu, \chi, \delta$ sufficiently small, we yield
  \begin{align*}
    &\norm{(v-\bar{v})(t)}^2 + \norm{(h-\bar{h})(t)}^2 + \norm{(u-\bar{u})(t)}^2 + \delta_S \int_0^t \dot{\X}(\tau) d\tau \\ 
    &+ \int_0^t G_2(\tau) + G_3(\tau) + G^R(\tau) + G^S(\tau) + D(\tau) + D_1(t) \,d\tau \\
    \leq&\; C \left( \norm{v_0-\bar{v}(0,\cdot)}^2 + \norm{h_0-\bar{h}(0,\cdot)}^2 + \norm{u_0-\bar{u}(0,\cdot)}^2 \right)  \\
    &\quad + C(\chi+\delta) \int_0^t \norm{\nabla^2(u-\bar{u})} d\tau + C \delta_R^{\frac{1}{3}} \,.
  \end{align*}
  Finally, since
  \begin{align*}
    &|\nabla(v-\bar{v})| \lesssim |h-\bar{h}| + |u-\bar{u}| + \pd_{1}v^R, \\
    &|h_0-\bar{h}_0(0,\cdot)| \lesssim |u_0-\bar{u}_0(0,\cdot)| + |v_0-\bar{v}_0(0,\cdot)| + \pd_{1}v^R ,
  \end{align*}
  and
  \begin{equation*}
    D_1(t) \sim \norm{\nabla(u-\bar{u})}^2,
  \end{equation*}
  we can acquire the desired inequality \eqref{eq:psi-L^2-est}, which completes the proof of Lemma \ref{lem:psi-0-estimate}.
\end{proof}


\subsection{Estimates on $\norm{\nabla(u-\bar{u})}$ }

\begin{lemma}\label{lem:psi-1-estimate}
  Under the assumptions of Proposition \ref{prop:apriori-estimate}, there exists constant $C>0$ (independent of $\nu,\delta,\chi$ and $T$)  such that for all $t \in [0,T]$, it holds
  \begin{equation}\label{eq:psi-1-est}
    \begin{aligned}
      &\norm{\nabla(u-\bar{u})(t)}^2 + \int_0^t \norm{\nabla^2(u-\bar{u})}^2 d\tau \\
      \leq&\; C \left( \norm{v_0-\bar{v}(0,\cdot)}_{H^1}^2 + \norm{u_0-\bar{u}(0,\cdot)}_{H^1}^2 \right) + C\delta_R^{\frac{1}{3}}.
    \end{aligned}
  \end{equation}
\end{lemma}

\begin{proof}
  Multiplying \eqref{eq:psi,phi-eqn}$_2$ by $-v\Delta\psi$, and integrating the resultant equation over $\Omega$, it follows
  \begin{equation*}
    \frac{d}{dt} \int_{\Omega} \frac{|\nabla\psi|^2}{2} \,dx + D_2(t) = \sum_{i=1}^5 J_i(t) ,
  \end{equation*}
  where
  \begin{align*}
    D_2(t) &:= \int_{\Omega} \mu v |\Delta\psi|^2 + (\mu+\lambda)v |\nabla\div\psi|^2 \,dx,  \\
    J_1(t) &:= \int_{\Omega} \div u \frac{|\nabla\psi|^2}{2} - \left( \nabla\psi \cdot \nabla u \right) : \nabla\psi \,dx, \\
    J_2(t) &:= - \dot{\X}(t) \int_{\Omega} \pd_{1}u^S_1\sdot \Delta\psi \,dx, \\
    J_3(t) &:= (\mu+\lambda) \int_{\Omega} \div\psi \left( \Delta\psi - \nabla\div\psi \right) \sdot \nabla v \,dx, \\
    J_4(t) &:= \int_{\Omega} v\Delta\psi \sdot \nabla(p-\bar{p}) \,dx,  \\
    J_5(t) &:= \int_{\Omega} vQ_3\sdot \Delta\psi_1 \,dx .
  \end{align*}
  Using assumption \eqref{eq:apriori-assumption}, Sobolev inequality and the  $L^2$ theory for elliptic equations, we have
  \begin{align*}
    J_1(t) &\leq C \norm{\nabla\psi}_{L^3} \norm{\nabla\psi}_{L^6} \norm{\nabla\psi} + C\delta \norm{\nabla\psi}^2 \\
    &\leq C\chi \norm{\nabla\psi}_{H^1} \norm{\nabla\psi} + C\delta \norm{\nabla\psi}^2 \\
    &\leq C\chi D_2(t) + C(\chi+\delta)\norm{\nabla\psi}^2 .
  \end{align*}
  By Holder inequality and Lemma \ref{lem:viscous shock}, we yield 
  \begin{equation*}
    J_2(t) \leq C |\dot{\X}(t)| \sdot \norm{\pd_{1}u^S_1} \sqrt{D_2(t)} \leq C \delta_S^{\frac{3}{2}} |\dot{\X}(t)| \sqrt{D_2(t)} \leq C\delta_S^2 |\dot{\X}(t)|^2 + C\delta_S D_2(t).
  \end{equation*}
  Using assumption \eqref{eq:apriori-assumption} and the fact that $\norm{\nabla\phi}_{L^3} \lesssim  \norm{\nabla\phi}_{H^1} \lesssim  \chi$, we obtain
  \begin{align*}
    J_3(t) &\leq C\delta \norm{\Delta\psi} \norm{\div\psi} + C \norm{\Delta\psi} \norm{\nabla\phi}_{L^3} \norm{\div\psi}_{L^6} \\
    &\leq C\delta \sqrt{D_2(t)} \norm{\nabla\psi} + C\chi \sqrt{D_2(t)} \norm{\div\psi}_{H^1} \\
    &\leq C(\chi+\delta)D_2(t) + C(\chi+\delta) \norm{\nabla\psi}^2. 
  \end{align*}
  It follows from Cauchy's inequality that
  \begin{equation*}
    J_4(t) \leq \frac{1}{4} D_2(t) + C D(t).
  \end{equation*}
  Using \eqref{eq:Q_3 chaifen}, we decompose $J_5(t)$ as follows:
  \begin{equation*}
    J_5(t) = \int_{\Omega} Q_3^S \sdot \Delta\psi_1 \,dx + \int_{\Omega} Q_3^R \sdot \Delta\psi_1 \,dx + \int_{\Omega} Q_3^I \sdot \Delta\psi_1 \,dx =: \sum_{ii=1}^3 J_{5,i}(t) .
  \end{equation*}
  According to Lemma \ref{lem:viscous shock}, it holds
  \begin{align*}
    J_{5,1}(t) &\leq C \int_{\Omega} \pd_{1}v^S\sdot \Delta\psi_1 \left( \left| h_1-\bar{h}_1 - \frac{p-\bar{p}}{\s_*} \right| + |p-\bar{p}| + |\nabla(p-\bar{p})| \right)  \,dx \\
    &\leq C \left( \delta_S \sqrt{G_3(t)} + \delta_S \sqrt{G^S(t)} + \delta_S^2 \sqrt{D(t)} \right) \sqrt{D_2(t)} \\
    &\leq C\delta_S G_3(t) + C\delta_S G^S(t) + C\delta_S D(t) + C\delta_S D_2(t) ,
  \end{align*}
  Using Lemma \ref{lem:appro-rare-property} and Cauchy's inequality, we have
  \begin{align*}
    J_{5,2}(t) &\leq C \sqrt{\delta_R} \mathcal{G}^R(t) + C \sqrt{\delta_R} D_2(t).
  \end{align*}
  By assumption \eqref{eq:apriori-assumption} and Lemma \ref{lem:interaction}, we get
  \begin{equation*}
    J_{5,3}(t) \leq C \norm{Q_3^I} \norm{\Delta\psi_1} \leq C\chi \delta_R \delta_S \sdot e^{-C\delta_S t},
  \end{equation*}
  Combining the above estimates, we have
  \begin{align*}
    \frac{d}{dt} \norm{\nabla\psi}^2 + D_2(t) \leq &\;C\delta_S^2 |\dot{\X}(t)|^2 + C\delta_S G_3(t) + C\delta_S G^S(t) + C D(t) \\
    & + C(\chi+\delta) \norm{\nabla\psi}^2 + C \sqrt{\delta_R} \mathcal{G}^R(t) + C\chi \delta_R \delta_S \sdot e^{-C\delta_S t} .
  \end{align*}
  Finally, integrating the above inequality over $[0,t]$, using Lemma \ref{lem:psi-0-estimate} and taking $\chi, \delta, \delta_R, \delta_S $ sufficiently small, we can obtain the desired inequality \eqref{eq:psi-1-est}, where we use the fact
  \begin{equation*}
    D_2(t) \sim \norm{\nabla^2\psi}^2 .
  \end{equation*}
  The proof of Lemma \ref{lem:psi-1-estimate} is completed.
\end{proof}


\subsection{Estimates on $\norm{\nabla^2(v-\bar{v})}$ }

\begin{lemma}\label{lem:phi-2-estimate}
  Under the assumptions of Proposition \ref{prop:apriori-estimate}, there exists constant $C>0$ (independent of $\nu,\delta,\chi$ and $T$)  such that for all $t \in [0,T]$, it holds
  \begin{equation}\label{eq:phi-2-est}
    \begin{aligned}
      &\norm{\nabla^2(v-\bar{v})(t)}^2 + \int_0^t \norm{\nabla^2(v-\bar{v})}^2 d\tau \\
      \leq&\; C \left( \norm{v_0-\bar{v}(0,\cdot)}_{H^2}^2 + \norm{u_0-\bar{u}(0,\cdot)}_{H^1}^2 \right) + C(\chi+\delta) \int_0^t \norm{\nabla^3(u-\bar{u})}^2 d\tau   + C\delta_R^{\frac{1}{3}}.
    \end{aligned}
  \end{equation}
\end{lemma}

\begin{proof}
  We first deform the perturbed system \eqref{eq:psi,phi-eqn} as follows:
  \begin{equation}\label{eq:psi,phi-bianxing}
    \left\{
    \begin{aligned}
      &\md{\phi} - v\div\psi - \dot{\X}(t) \pd_{1}v^S = -vQ_1, \\
      &\md{\psi} + vp'(v)\nabla\phi + v \big( p'(v) - p'(\bar{v}) \big) \nabla\bar{v} - \dot{\X}(t) \pd_{1}u^S \\
      & \hspace{6em}= \mu v\Delta\psi + (\mu+\lambda)v\nabla\div\psi + (2\mu+\lambda) v\Delta u^R - vQ_3\sdot e_1.
    \end{aligned} 
    \right. 
  \end{equation}
  Applying $\pd_{ij}\; (i,j=1,2,3)$ to \eqref{eq:psi,phi-bianxing}$_1$ leads to
  \begin{equation}\label{eq:pd_ij-phi}
    \begin{aligned}
      &\md{\pd_{ij}\phi} + \pd_{i}u \sdot \nabla\pd_{j}u + \pd_{j}u \sdot \nabla\pd_{i}\phi + \pd_{ij}u \sdot \nabla\phi  - \dot{\X}(t) \pd_{1ij}v^S \\
      =&\; \pd_{ij}v \sdot \div\psi + \pd_{i}v \sdot \pd_{j}\div\psi + \pd_{j}v \sdot 
      \pd_{i}\div\psi + v\pd_{ij}\div\psi - \pd_{ij}(vQ_1). 
    \end{aligned} 
  \end{equation}
  Multiplying \eqref{eq:pd_ij-phi} by $\rho (2\mu+\lambda) \pd_{ij}\phi$ and integrating the resultant equation over $\Omega$, we have
  \begin{equation}\label{eq:4.61}
    \begin{aligned}
      &(2\mu+\lambda) \frac{d}{dt} \int_{\Omega} \rho \frac{|\nabla^2\phi|^2}{2} \,dx - (2\mu+\lambda) \int_{\Omega} \nabla\pd_{i}\phi \sdot \nabla\pd_{i}\div\psi \,dx \\
      =& -(2\mu+\lambda) \int_{\Omega} 2\rho \pd_{ij}\phi \sdot \pd_{i}u \sdot \nabla\pd_{j}\phi \,dx - (2\mu+\lambda) \int_{\Omega} \rho\pd_{ij}\phi \sdot \pd_{ij}u \sdot \nabla\phi \,dx \\
      & + (2\mu+\lambda) \int_{\Omega} \rho\pd_{ij}\phi \sdot \pd_{ij}v \sdot \div\psi \,dx + (2\mu+\lambda) \int_{\Omega} 2\rho \pd_{ij}\phi \sdot \pd_{i}v \sdot \pd_{j}\div\psi \,dx \\
      & + (2\mu+\lambda) \dot{\X}(t) \int_{\Omega} \rho \pd_{111}v^S \sdot \pd_{11}\phi \,dx - (2\mu+\lambda) \int_{\Omega} \rho \pd_{ij}\phi \sdot \pd_{ij}(vQ_1) \,dx .
    \end{aligned}
  \end{equation}
  Note that throughout this paper we obey the Einstein summation convention, i.e. repeated indices are implicitly summed from 1 to 3.
  
  Applying $\pd_{i} \; (i=1,2,3)$ to \eqref{eq:psi,phi-bianxing}$_2$, it follows
  \begin{equation}\label{eq:pd_i-psi}
    \begin{aligned}
      &\md{\pd_{i}\psi} + \pd_{i}u \sdot \nabla\psi + vp'(v)\nabla\pd_{i}\phi + \pd_{i}(vp'(v))\nabla\phi  - \dot{\X}(t) \pd_{1i}u^S_1 \\
      =& -v(p'(v) - p'(\bar{v}))\nabla\pd_{i}\bar{v} - \pd_{i}\big[ v(p'(v) - p'(\bar{v})) \big]\nabla\bar{v} +\mu v \Delta\pd_{i}\psi + (\mu+\lambda)v\nabla\pd_{i}\div\psi \\
      & +(2\mu+\lambda)v\Delta\pd_{i}u^R - \pd_{i}v \left( \mu\Delta\psi + (\mu+\lambda)\nabla\div\psi + (2\mu+\lambda) \Delta u^R  \right) - \pd_{i}(vQ_3)\sdot e_1 .
    \end{aligned} 
  \end{equation}
  Multiplying \eqref{eq:pd_i-psi} by $ -\rho \nabla\pd_{i}\phi$ and integrating the resultant equation over $\Omega$, we have
  \begin{equation}\label{eq:4.62}
    \begin{aligned}
      &\int_{\Omega} -p'(v)|\nabla^2\phi|^2 \,dx + (2\mu+\lambda) \int_{\Omega} \nabla\pd_{i}\phi \sdot \nabla\pd_{i}\div\psi \,dx \\
      =&\; \frac{d}{dt} \int_{\Omega} \rho\pd_{i}\psi \sdot \nabla\pd_{i}\phi \,dx - \int_{\Omega} \rho\pd_{i}\psi_j [\md{\pd_{ij}\phi}] \,dx + \int_{\Omega} \rho(\pd_{i}u \sdot \nabla\psi) \sdot \nabla\pd_{i}\phi \,dx \\
      & + \int_{\Omega} \rho \pd_{i}(vp'(v))\nabla\phi \sdot \nabla\pd_{i}\phi \,dx + \int_{\Omega} \rho \pd_{i}\big[ v(p'(v) - p'(\bar{v})) \big] \pd_{1}\bar{v} \sdot \pd_{1i}\phi \,dx \\
      & + \int_{\Omega} (p'(v) - p'(\bar{v})) \pd_{11}\bar{v} \sdot \pd_{11}\phi \,dx - \dot{\X}(t)\! \int_{\Omega} \rho \pd_{11}u^S_1 \sdot \pd_{11}\phi \,dx + \int_{\Omega} \rho \pd_{i}(vQ_3) \pd_{1i}\phi \,dx  \\
      & - (2\mu+\lambda) \int_{\Omega} \pd_{111}u^R_1 \sdot \pd_{11}\phi \,dx - (2\mu+\lambda) \int_{\Omega} \rho \pd_{i}v \sdot \pd_{11}u^R_1 \sdot \pd_{11}\phi \,dx \\
      & - (2\mu+\lambda) \int_{\Omega} \rho\pd_{i}v (\mu\Delta\psi + (\mu+\lambda)\nabla\div\psi ) \sdot \nabla \pd_{i}\phi \,dx \,.
    \end{aligned}
  \end{equation}
  Adding \eqref{eq:4.61} and \eqref{eq:4.62} together, we arrive at
  \begin{equation}\label{eq:phi-2-weifen}
    (2\mu+\lambda) \frac{d}{dt} \int_{\Omega} \rho \frac{|\nabla^2\phi|^2}{2} \,dx + D_\phi(t) = \sum_{i=1}^9 K_i(t),
  \end{equation}
  where
  \begin{align*}
      D_\phi(t) &= \int_{\Omega} |p'(v)|\sdot |\nabla^2\phi|^2 \,dx \,,  \\
      K_1(t) &= \frac{d}{dt} \int_{\Omega} \rho\pd_{i} \psi \sdot \nabla\pd_{i}\phi \,dx \,, \\
      K_2(t) &= \dot{\X}(t) \int_{\Omega} \rho \pd_{11}\phi \big[ (2\mu+\lambda) \pd_{111}v^S - \pd_{11}u^S_1 \big]  \,dx \,, \\
      K_3(t) &= -(2\mu+\lambda) \int_{\Omega} \rho\pd_{ij}\phi \, \pd_{ij}(vQ_1) \,dx \,,  \\
      K_4(t) &= \int_{\Omega} \rho \pd_{1i}\phi \, \pd_{i}(vQ_3) \,dx \,,  \\
      K_5(t) &= -\int_{\Omega} \rho\pd_{i}\psi_j \big[ \md{\pd_{ij}\phi} \big]    \,dx  \,,  \\
      K_6(t) &= -(2\mu+\lambda) \int_{\Omega} \rho\pd_{ij}\phi \, (2\pd_{i}u \sdot \nabla\pd_{j}\phi + \pd_{ij}u \sdot \nabla\phi) \,dx  \\ & + \int_{\Omega} \rho\nabla\pd_{i}\phi \sdot \big[ \pd_{i}u \sdot \nabla\psi + \pd_{i}(vp'(v))\nabla\phi \big] \,dx  \,,  \\
      K_7(t) &= \int_{\Omega} \pd_{11}\bar{v} \big[p'(v)-p'(\bar{v})\big] \pd_{11}\phi \,dx + \int_{\Omega} \rho\pd_{1}\bar{v}\, \pd_{i}\big[v\big(p'(v)-p'(\bar{v})\big)\big] \pd_{1i}\phi \,dx \,,  \\
      K_8(t) &= (2\mu+\lambda)\int_{\Omega} \rho\pd_{ij}\phi \big[ \pd_{ij}v\,\div\psi + 2\pd_{i}v\, \pd_{j}\div\psi \big] \,dx \\ &- \int_{\Omega} \rho\pd_{i}v \big[ \mu\Delta\psi + (\mu+\lambda)\nabla\div\psi \big] \nabla\pd_{i}\phi \,dx \,,  \\
      K_9(t) &= -(2\mu+\lambda)\int_{\Omega} \pd_{111}u^R_1\, \pd_{11}\phi + \rho\pd_{11}u^R_1\, \pd_{i}v \pd_{1i}\phi \,dx  \,.
  \end{align*}
  Using Cauchy's inequality, we have
  \begin{align*}
      &\int_0^t K_1(\tau) d\tau = \int_{\Omega} \rho\pd_{i} \psi \sdot \nabla\pd_{i}\phi \,dx \Big|_{\tau=0}^{\tau=t}  \\
      \leq& \frac{2\mu+\lambda}{8} \norm{\sqrt{\rho}\nabla\phi(t)}^2 + C\norm{\nabla\psi(t)}^2 + C\big( \norm{\nabla\psi(0,\cdot)}^2 + \norm{\nabla^2\phi(0,\cdot)}^2 \big).
  \end{align*}
  By Lemma \ref{lem:viscous shock} and Holder inequality, it holds
  \begin{equation*}
    K_2(t) \lesssim \delta_S |\dot{\X}(t)| \int_{\Omega} \pd_{1}v^S \sqrt{|p'(v)|}|\nabla^2\phi| \,dx \lesssim \delta_S^2 |\dot{\X}(t)|^2 + \delta_S^2 D_\phi(t).
  \end{equation*}

  Using the decomposition of $vQ_1$ in \eqref{eq:Q_1 chaifen}, it follows from Cauchy's inequality that
  \begin{align*}
    K_3(t) \leq C\int_{\Omega} |\nabla^2\phi|\sdot |\nabla^2(vQ_1)| \,dx \leq C\sum_{i=1}^3 K_{3,i}(t) ,
  \end{align*}
  where
  \begin{align*}
      K_{3,1}(t) &= \int_{\Omega} |\nabla^2\phi|\sdot |\nabla^2(Q_1^R)| \,dx \,, \\
      K_{3,2}(t) &= \int_{\Omega} |\nabla^2\phi|\sdot |\nabla^2(Q_1^S)| \,dx \,, \\
      K_{3,3}(t) &= \int_{\Omega} |\nabla^2\phi|\sdot |\nabla^2(Q_1^I)| \,dx \,. 
  \end{align*}
  Using \eqref{eq:L3-norm}, Lemma \ref{lem:appro-rare-property} and Holder inequality, we have
  \begin{align*}
      K_{3,1}(t) &\leq C\int_{\Omega} |\nabla^2\phi|\sdot \pd_{1}v^R \big( |\psi_1|+|\phi|+|\nabla\psi_1|+|\nabla\phi|+ |\nabla^2\psi_1|+|\nabla^2\phi|\big) \,dx  \\
      &\leq C \norm{\pd_{1}u^R_1}_{L^6} \norm{(\phi,\psi)}_{L^3} \norm{\nabla^2\phi} + C\delta \big( \norm{\nabla(\phi,\psi)}^2 + \norm{\nabla^2(\phi,\psi)}^2 \big) \\
      &\leq C\chi\norm{\pd_{1}u^R_1}_{L^6}^2 + C(\chi+\delta) \big( \norm{\nabla(\phi,\psi)}^2 + \norm{\nabla^2(\phi,\psi)}^2 \big)  \\
      &\leq C\chi\, \delta_R^{\frac{1}{3}}\sdot (1+t)^{-\frac{5}{3}} + C(\chi+\delta) \big( \norm{\nabla(\phi,\psi)}^2 + \norm{\nabla^2(\phi,\psi)}^2 \big).
  \end{align*}
  Notice that
  \begin{equation}\label{eq:nabla-phi-est}
      \norm{\nabla\phi}^2 \lesssim D(t) + \delta_S^2 G^S(t) + \delta_R G^R(t),
  \end{equation}
  so
  \begin{align*}
      K_{3,1}(t) &\leq C\chi\, \delta_R^{\frac{1}{3}}\sdot (1+t)^{-\frac{5}{3}} + C(\chi+\delta) \big( \norm{\nabla^2\phi}^2 + \norm{\nabla\psi}_{H^1}^2 \big) \\
      &\quad + C(\chi+\delta) \big(D(t) + G^S(t) + G^R(t)\big).
  \end{align*}
  By the definition of $h$, it holds
  $$ 
  \psi_1 = (h_1-\bar{h}_1) + (2\mu+\lambda)\pd_{1}\phi - (2\mu+\lambda)\pd_{1}v^R \,,
  $$
  and so
  \begin{align*}
    |\psi_1| &\lesssim |h_1-\bar{h}_1| + |\nabla(p-\bar{p})| +\pd_{1}\bar{v}\sdot |p-\bar{p}| + \pd_{1}v^R \\
      &\lesssim \left| h_1-\bar{h}_1 - \frac{p-\bar{p}}{\s_*} \right| + |\nabla(p-\bar{p})| + |p-\bar{p}| + \pd_{1}u^R_1 \,,
  \end{align*}
  which together with $ |\phi| \sim |p-\bar{p}| $ lead to
  \begin{align*}
      K_{3,2}(t) &\lesssim \int_{\Omega} |\nabla^2\phi|\sdot \pd_{1}v^S \big( |\psi_1|+|\phi|+|\nabla\psi_1|+|\nabla\phi|+ |\nabla^2\psi_1|+|\nabla^2\phi|\big) \,dx  \\
      &\lesssim \int_{\Omega} |\nabla^2\phi|\sdot \pd_{1}v^S \left( \left| h_1-\bar{h}_1 - \frac{p-\bar{p}}{\s_*} \right| + |\nabla(p-\bar{p})| + |p-\bar{p}| + \pd_{1}u^R_1 \right) dx \\
      &\quad + \delta \big( \norm{\nabla(\phi,\psi)}^2 + \norm{\nabla^2(\phi,\psi)}^2 \big) \\
      &\lesssim \chi \norm{\pd_{1}u^R_1\sdot \pd_{1}v^S} + \delta \big( G_2(t) + D(t) + G^S(t) \big) + \delta \big( \norm{\nabla(\phi,\psi)}^2 + \norm{\nabla^2(\phi,\psi)}^2 \big) \,.
  \end{align*}
  According to \ref{lem:interaction}, it holds
  \begin{equation*}
      \norm{\pd_{1}u^R_1\sdot \pd_{1}v^S} \leq C\delta_R \delta_S e^{-C\delta_S t} .
  \end{equation*}
  In addition,
  \begin{equation*}
      \norm{\nabla\phi}^2 \lesssim D(t) + \delta_S^2 G^S(t) + \delta_R^2 G^R(t) .
  \end{equation*}
  Hence, we have
  \begin{align*}
      K_{3,2}(t) &\leq C\delta_R \delta_S e^{-C\delta_S t} + C\delta \big( G_2(t) + D(t) + G^R(t) + G^S(t) \big)  \\
      &\quad + C\delta \norm{\nabla^2\phi}^2 + C\delta \norm{\nabla\psi}_{H^1}^2 .
  \end{align*}
  It follows from Lemma \ref{lem:interaction} that
  \begin{equation*}
    \norm{\nabla^2 Q_1^I} \leq C\delta_R \delta_S e^{-C\delta_S t},
  \end{equation*}
  so 
  \begin{equation*}
      K_{3,3}(t) \leq \norm{\nabla^2\phi}\sdot \norm{\nabla^2 Q_1^I} \leq C\chi\delta_R \delta_S e^{-C\delta_S t} .
  \end{equation*}
  Combining the above estimates, we yield
  \begin{align*}
    K_3(t) &\leq C\chi\, \delta_R^{\frac{1}{3}}\sdot (1+t)^{-\frac{5}{3}} + C\delta_R \delta_S e^{-C\delta_S t} + C\delta D_\phi(t) + C\delta \norm{\nabla\psi}_{H^1}^2 \\
    &\quad + C(\chi+\delta) \big( G_2(t) + D(t) + G^R(t) + G^S(t) \big).
  \end{align*}

  Using the decomposition of $vQ_3$ \eqref{eq:Q_3 chaifen}, we have
  \begin{align*}
      K_4(t) \leq C\int_{\Omega} |\nabla^2\phi|\sdot |\nabla(vQ_3)| \,dx \leq C\sum_{i=1}^3 K_{4,i}(t),
  \end{align*}
  where
  \begin{align*}
    K_{4,1}(t) &= \int_{\Omega} |\nabla^2\phi|\sdot |\nabla(Q_3^R)| \,dx \,, \\
    K_{4,2}(t) &= \int_{\Omega} |\nabla^2\phi|\sdot |\nabla(Q_3^S)| \,dx \,, \\
    K_{4,3}(t) &= \int_{\Omega} |\nabla^2\phi|\sdot |\nabla(Q_3^I)| \,dx \,. 
  \end{align*}
  Similarly, we can obtain
  \begin{align*}
      K_{4,1}(t) &\leq C\int_{\Omega} |\nabla^2\phi|\sdot \pd_{1}u^R_1 \big( |\psi_1|+|\phi|+|\nabla\psi_1|+|\nabla\phi| \big) \,dx  \\
      &\leq C \norm{\pd_{1}u^R_1}_{L^6} \norm{(\phi,\psi)}_{L^3} \norm{\nabla^2\phi} + C\delta_R \norm{\nabla(\phi,\psi)}^2 + C\delta_R \norm{\nabla^2\phi}^2 \\
      &\leq C\chi\norm{\pd_{1}u^R_1}_{L^6}^2 + C(\chi+\delta) \norm{\nabla(\phi,\psi)}^2 + C(\chi+\delta)\norm{\nabla^2\phi}^2 \\
      &\leq C\chi\, \delta_R^{\frac{1}{3}}\sdot (1+t)^{-\frac{5}{3}} + C(\chi+\delta) \big( \norm{\nabla(\phi,\psi)}^2 + \norm{\nabla^2\phi}^2 \big) \\
      &\leq C\chi\, \delta_R^{\frac{1}{3}}\sdot (1+t)^{-\frac{5}{3}} + C(\chi+\delta) \big( \norm{\nabla^2\phi}^2 + \norm{\nabla\psi}^2 \big) \\
      &\quad + C(\chi+\delta) \big(D(t) + G^S(t) + G^R(t)\big) ,
  \end{align*}
  and
  \begin{align*}
      K_{4,2}(t) &\leq C\int_{\Omega} |\nabla^2\phi|\sdot \pd_{1}v^S \big( |\psi_1|+|\phi|+|\nabla\psi_1|+|\nabla\phi| \big) \,dx  \\
      &\leq C\int_{\Omega} |\nabla^2\phi|\sdot \pd_{1}v^S \left( \left| h_1-\bar{h}_1 - \frac{p-\bar{p}}{\s_*} \right| + |\nabla(p-\bar{p})| + |p-\bar{p}| + \pd_{1}u^R_1 \right) dx \\
      &\quad + C\delta \norm{\nabla(\phi,\psi)}^2 + C\delta \norm{\nabla^2\phi}^2 \\
      &\leq C\delta_R \delta_S e^{-C\delta_S t} + C\delta \big( G_2(t) + D(t) + G^R(t) + G^S(t) + \norm{\nabla^2\phi}^2 + \norm{\nabla\psi}^2 \big).
  \end{align*}
  For $ K_{4,3}(t) $, the only nontrivial thing is to estimate
  \begin{equation*}
      \int_{\Omega} |\nabla^2\phi|\sdot \left| \nabla \big(v\pd_{1}(\bar{p}-p^S-p^R)\big) \right| dx .
  \end{equation*}
  It follows from direct computations that
  \begin{align*}
    &|\nabla\pd_{1}(\bar{p}-p^S-p^R)| \\
    \lesssim \;& \left| \nabla\big[ (p'(\bar{v})-p'(v^R))\pd_{1}v^R \big] \right| + \left| \nabla\big[ (p'(\bar{v})-p'(v^S))\pd_{1}v^S \big] \right|  \\
    \lesssim \;& |\bar{v}-v^R|\pd_{1}v^R + |\bar{v}-v^S|\pd_{1}v^S + \pd_{1}v^R \sdot \pd_{1}v^S \,,
  \end{align*}
  thus by Lemma \ref{lem:interaction}, we have
  \begin{equation*}
      \norm{\nabla\pd_{1}(\bar{p}-p^S-p^R)} \leq C\delta_R \delta_S e^{-C\delta_S t} .
  \end{equation*}
  Notice that
  \begin{equation*}
    \norm{\pd_{1}(\bar{p}-p^S-p^R)}_{L^\infty} \lesssim \norm{\pd_{1}v^R}_{L^\infty} + \norm{\pd_{1}v^S}_{L^\infty} \lesssim \delta_R +\delta_S \lesssim \delta \,.
  \end{equation*}
  Hence, it holds
  \begin{align*}
      &\int_{\Omega} |\nabla^2\phi|\sdot \left|\nabla \big(v\pd_{1}(\bar{p}-p^S-p^R)\big)\right| dx \\
      \leq\;& C \int_{\Omega} |\nabla^2\phi|\sdot |\nabla\pd_{1}(\bar{p}-p^S-p^R)| dx + C \int_{\Omega} |\nabla^2\phi|\sdot |\nabla v|\sdot |\pd_{1}(\bar{p}-p^S-p^R)| \,dx \\
      \leq\;& C\chi\,\delta_R \delta_S e^{-C\delta_S t} + C\delta \int_{\Omega} |\nabla^2\phi|\sdot |\nabla \phi| \,dx  \\
      \leq\;& C\chi\,\delta_R \delta_S e^{-C\delta_S t} + C\delta\big(D(t) + G^S(t) + G^R(t)\big) +C\delta \norm{\nabla^2\phi}^2 \,.
  \end{align*}
  Therefore, we have
  \begin{align*}
      K_4(t) &\leq C\chi\, \delta_R^{\frac{1}{3}}\sdot (1+t)^{-\frac{5}{3}} + C\delta_R \delta_S e^{-C\delta_S t} + C\delta D_\phi(t) + C\delta \norm{\nabla\psi}^2 \\
      &\quad + C(\chi+\delta) \big( G_2(t) + D(t) + G^R(t) + G^S(t) \big).
    \end{align*}
  For $ K_5(t) $, according to \eqref{eq:pd_ij-phi}, it holds
  \begin{equation*}
      K_5(t) = -\int_{\Omega} \rho\pd_{i}\psi_j \big[ \md{\pd_{ij}\phi} \big] \,dx =: \sum_{i=1}^5 K_{5,i}(t) ,
  \end{equation*}
  where
  \begin{align*}
      K_{5,1}(t) &= -\dot{\X}(t) \int_{\Omega} \rho \pd_{1}\psi_1 \pd_{111}v^S \,dx\,,  \\
      K_{5,2}(t) &= \int_{\Omega} \rho\pd_{i}\psi_j \big[ \pd_{i}u\sdot \nabla\pd_{j}\phi +
      \pd_{j}u\sdot \nabla\pd_{i}\phi + \pd_{ij}u\sdot \nabla\phi \big]  \,dx \,, \\
      K_{5,3}(t) &= -\int_{\Omega} \rho\pd_{i}\psi_j \big[ \pd_{j}\psi \sdot \pd_{i}\div\psi + \pd_{i}\psi \sdot \pd_{j}\div\psi + \pd_{ij}v \sdot \div\psi \big]  \,dx \,, \\
      K_{5,4}(t) &= -\int_{\Omega} \pd_{i}\psi_j \sdot \pd_{ij}\div\psi \,dx \,, \\
      K_{5,5}(t) &= \int_{\Omega} \rho\pd_{i}\psi_j \sdot \pd_{ij}(vQ_1) \,dx .
  \end{align*}
  Using Lemma \ref{lem:viscous shock} and Cauchy's inequality, we have
  \begin{equation*}
      K_{5,1}(t) \leq C\delta_S^2 |\dot{\X}(t)|^2 + C\delta_S^2 \norm{\nabla\psi}^2 .
  \end{equation*}
  By \eqref{eq:L3-norm} and \eqref{eq:nabla-phi-est}, it holds
  \begin{align*}
      K_{5,2}(t) &\leq C \int_{\Omega} |\nabla\psi| \big(|\nabla^2\psi||\nabla\phi| + |\nabla\psi||\nabla^2\phi| \big) + |\nabla\psi| \big(|\pd_{11}\bar{u}_1||\nabla\phi| + |\pd_{1}\bar{u}_1||\nabla^2\phi| \big) \,dx  \\
      &\leq C \norm{\nabla\psi}_{L^6} \norm{\nabla(\phi,\psi)}_{L^3} \norm{\nabla^2(\phi,\psi)} + C\delta \norm{\nabla\psi} \left( \norm{\nabla\phi} + \norm{\nabla^2\phi} \right) \\
      &\leq C(\chi+\delta)\left( \norm{\nabla^2\phi}^2 + \norm{\nabla\psi}_{H^1}^2 \right) + C\delta \left( G^R(t) + G^S(t) +D(t) \right) .
  \end{align*}
  Similar to $K_{5,2}(t)$, we have
  \begin{align*}
      K_{5,3}(t) &\leq C \int_{\Omega} |\nabla\psi| \big(|\nabla\psi||\nabla^2\phi| + |\nabla^2\psi||\nabla\phi| \big) + |\nabla\psi| \big(|\pd_{11}\bar{u}_1||\nabla\psi| + |\pd_{1}\bar{u}_1||\nabla^2\psi| \big) \,dx  \\
      &\leq C \norm{\nabla\psi}_{L^6} \norm{\nabla(\phi,\psi)}_{L^3} \norm{\nabla^2(\phi,\psi)} + C\delta \norm{\nabla\psi} \left( \norm{\nabla\psi} + \norm{\nabla^2\psi} \right) \\
      &\leq C(\chi+\delta)\left( \norm{\nabla^2\phi}^2 + \norm{\nabla\psi}_{H^1}^2 \right).
  \end{align*}
  Integration by parts over $\Omega$ leads to
  \begin{equation*}
      K_{5,4}(t) = -\int_{\Omega} \pd_{i}\psi_j \sdot \pd_{ij}\div\psi \,dx = \int_{\Omega} \pd_{i}\div\psi \sdot \pd_{i}\div\psi \,dx = \norm{\nabla\div\psi}^2 \lesssim \norm{\nabla^2\psi}^2.
  \end{equation*}
  Essentially the same with $K_3(t)$, we have
  \begin{align*}
      K_{5,5}(t) &\leq C\chi\, \delta_R^{\frac{1}{3}}\sdot (1+t)^{-\frac{5}{3}} + C\delta_R \delta_S e^{-C\delta_S t} + C\delta D_\phi(t) + C\delta \norm{\nabla\psi}_{H^1}^2 \\
      &\quad + C(\chi+\delta) \big( G_2(t) + D(t) + G^R(t) + G^S(t) \big).
  \end{align*}
  Combining the above estimates, we yield
  \begin{align*}
    K_5(t) &\leq C\chi\, \delta_R^{\frac{1}{3}}\sdot (1+t)^{-\frac{5}{3}} + C\delta_R \delta_S e^{-C\delta_S t} + C\delta D_\phi(t) +C\delta_S^2 |\dot{\X}(t)|^2  \\
    &\quad + C\norm{\nabla\psi}_{H^1}^2 + C(\chi+\delta) \big( G_2(t) + D(t) + G^R(t) + G^S(t) \big).
  \end{align*}
  Using \eqref{eq:L3-norm}, Cauchy's inequality and the fact that
  \begin{equation*}
    \norm{\nabla^2\phi}^2 \norm{\nabla\psi}_{L^\infty} \lesssim \norm{\nabla^2\phi}^2 \norm{\nabla\psi}_{H^2} \lesssim \chi \norm{\nabla^2\phi} \norm{\nabla\psi}_{H^2} \lesssim \chi D_\phi(t) + \chi \norm{\nabla\psi}_{H^2}^2 ,
  \end{equation*}
  we have
  \begin{align*}
      K_6(t) &\leq C \int_{\Omega} |\nabla^2\phi| \big(|\nabla^2\psi||\nabla\phi| + |\nabla\psi||\nabla^2\phi| + |\nabla\psi|^2 +|\nabla\phi|^2 \big) \\ &\quad + |\nabla^2\phi| \big(|\pd_{11}\bar{u}_1||\nabla\phi| + |\pd_{1}\bar{u}_1||\nabla\psi| + |\pd_{1}\bar{u}_1||\nabla^2\phi| + |\pd_{1}\bar{v}|\nabla\phi|| \big) \,dx  \\
      &\leq C\norm{\nabla^2\phi} \left( \norm{\nabla^2\psi}_{L^6} \norm{\nabla\phi}_{L^3} +  \norm{\nabla\psi}_{L^6} \norm{\nabla\psi}_{L^3} + \norm{\nabla^2\phi} \norm{\nabla\psi}_{L^\infty}\right)  \\  &\quad + C \norm{\nabla\phi}_{L^6} \norm{\nabla\phi}_{L^3} \norm{\nabla^2\phi} + C\delta \norm{\nabla^2\phi} \left( \norm{\nabla(\phi,\psi)} + \norm{\nabla^2\phi} \right) \\
      &\leq C(\chi+\delta) \left( D_\phi(t) + D(t) + G^R(t) + G^S(t) + \norm{\nabla\psi}_{H^2}^2 \right) .
  \end{align*}

  Notice that $ vp'(v) = -\gamma p(v) $, so it holds
  \begin{align*}
      &\nabla \big[v\big(p'(v)-p'(\bar{v})\big)\big] = \nabla \big[ vp'(v)- \bar{v}p'(\bar{v}) - \phi p'(\bar{v}) \big] 
      \lesssim |\nabla(p-\bar{p})| + |\nabla\phi| + |\phi| |\pd_{1}\bar{v}|.
  \end{align*}
  Hence, we have
  \begin{align*}
    K_7(t) &\leq C\int_{\Omega} |\pd_{11}\bar{v}| |\phi| |\nabla^2\phi| dx + C\int_{\Omega} |\pd_{1}\bar{v}| |\nabla^2\phi| \big( |\nabla(p-\bar{p})| + |\nabla\phi| + |\phi| |\pd_{1}\bar{v}| \big) dx \\
    &\leq C \int_{\Omega} \left( \pd_{11}v^R + |\pd_{1}v^R|^2 \right) |\phi| |\nabla^2\phi| \,dx + C \int_{\Omega} \pd_{1}v^S |p-\bar{p}| |\nabla^2\phi| \,dx  \\ &\quad + C\delta \int_{\Omega} |\nabla^2\phi|^2 + |\nabla(p-\bar{p})|^2 + |\nabla\phi|^2 \,dx \\
    &\leq C\sqrt{\delta} \left( G^R(t) + D_\phi(t) \right) + C\delta \left( D(t)+G^S(t) \right) .
  \end{align*}
  Similar to $K_6(t)$, it holds
  \begin{align*}
    K_8(t) &\leq C \int_{\Omega} |\nabla^2\phi| \big( |\nabla^2\psi||\nabla\phi| + |\nabla\psi||\nabla^2\phi| + |\pd_{1}\bar{v}| |\nabla^2\psi| + |\pd_{11}\bar{v}| |\nabla\psi| \big) dx  \\
    &\leq C \norm{\nabla^2\phi} \left( \norm{\nabla^2\psi}_{L^6} \norm{\nabla\phi}_{L^3} + \norm{\nabla\psi}_{L^\infty} \norm{\nabla^2\phi} + \delta \norm{\nabla\psi}_{H^1} \right) \\
    &\leq C \chi D_\phi(t) +C(\chi+\delta) \norm{\nabla\psi}_{H^2}^2.
  \end{align*}
  Using Lemma \ref{lem:appro-rare-property}, we get
  \begin{align*}
    K_9(t) &\leq C \int_{\Omega} |\pd_{11}v^R|\sdot |\nabla^2\phi| \,dx + C \int_{\Omega} |\pd_{1}v^R|\sdot |\nabla\phi|\sdot |\nabla^2\phi| \,dx \\
    &\leq C \norm{\pd_{11}v^R} \norm{\nabla^2\phi} + \norm{\pd_{1}v^R}_{L^6} \norm{\nabla\phi}_{L^3} \norm{\nabla^2\phi} \\
    &\leq C(\delta_R)^{-\frac{1}{2}}\norm{\pd_{11}v^R}^2 + C\chi \norm{\pd_{1}v^R}_{L^6}^2 +C(\chi + \sqrt{\delta_R}) \norm{\nabla^2\phi}^2 \\
    &\leq C(\delta_R)^{-\frac{1}{2}}\norm{\pd_{11}v^R}^2 + C\chi\,\delta_R^{\frac{1}{3}}\sdot (1+t)^{-\frac{5}{3}} + C(\chi + \sqrt{\delta_R}) D_\phi(t).
  \end{align*}
  In addition, it follows from Lemma \ref{lem:appro-rare-property} that
  \begin{equation*}
      \norm{\pd_{11}v^R}^2 \lesssim 
    \begin{cases}
      \delta_R\;, \quad       &1+t \leq \delta_R ^{-1} \\
      \frac{1}{(1+t)^2} \;, \quad  &1+t \geq \delta_R ^{-1}
    \end{cases},
  \end{equation*}
  which yields
  \begin{equation*}
    \int_0^t (\delta_R)^{-\frac{1}{2}}\norm{\pd_{11}v^R}^2 d\tau  \leq \int_0^\infty (\delta_R)^{-\frac{1}{2}}\norm{\pd_{11}v^R}^2 d\tau \lesssim \sqrt{\delta_R} \,.
  \end{equation*}
  Finally, notice the fact that
  \begin{equation*}
    D_\phi(t) \sim \norm{\nabla^2\phi}^2 .
  \end{equation*}
  Thus, integrating \eqref{eq:phi-2-weifen} over $[0,t]$ and using the above estimates, we can deduce the desired inequality \eqref{eq:phi-2-est} from Lemma \ref{lem:psi-0-estimate} and Lemma \ref{lem:psi-1-estimate}. The proof of Lemma \ref{lem:phi-2-estimate} is completed.
\end{proof}
\vspace{1em}


\subsection{Estimates on $\norm{\nabla^2(u-\bar{u})}$ }

\begin{lemma}\label{lem:psi-2-estimate}
  Under the assumptions of Proposition \ref{prop:apriori-estimate}, there exists constant $C>0$ (independent of $\nu,\delta,\chi$ and $T$)  such that for all $t \in [0,T]$, it holds
  \begin{equation}\label{eq:psi-2-est}
    \begin{aligned}
      &\norm{\nabla^2(u-\bar{u})(t)}^2 + \int_0^t \norm{\nabla^3(u-\bar{u})}^2 d\tau \\
      \leq&\; C \left( \norm{v_0-\bar{v}(0,\cdot)}_{H^2}^2 + \norm{u_0-\bar{u}(0,\cdot)}_{H^2}^2 \right) + C\delta_R^{\frac{1}{3}}.
    \end{aligned}
  \end{equation}
\end{lemma}

\begin{proof}
  Multiplying \eqref{eq:pd_i-psi} by $-\Delta\pd_{i}\psi$ and summing $i$ from 1 to 3, then integrating the resultant equation over $\Omega$, we have
  \begin{equation}\label{eq:psi-2-weifen}
    \frac{d}{dt} \int_{\Omega} \frac{|\nabla^2\psi|^2}{2} \,dx + D_3(t) =: \sum_{i=1}^8 L_i(t),
  \end{equation}
  where
  \begin{align*}
    D_3(t) &:= \mu\int_{\Omega} v|\nabla\Delta\psi|^2 \,dx + (\mu+\lambda) \int_{\Omega} v|\nabla^2\div\psi|^2 \,dx , \\
    L_1(t) &:= \int_{\Omega} \frac{|\nabla^2\psi|^2}{2} \sdot \div u \,dx - \int_{\Omega} (\pd_{j}u \sdot \nabla)\pd_{i}\psi \cdot \pd_{ij}\psi \,dx ,    \\
    L_2(t) &:= \int_{\Omega} (\pd_{i}u \sdot \nabla)\psi \cdot \Delta\pd_{i}\psi + \pd_{i}(vp'(v))\nabla\phi \cdot \Delta\pd_{i}\psi \,dx ,       \\
    L_3(t) &:= \int_{\Omega} vp'(v)\nabla\pd_{i}\phi \cdot \Delta\pd_{i}\psi \,dx ,     \\
    L_4(t) &:= \int_{\Omega} v(p'(v) - p'(\bar{v}))\pd_{11}\bar{v} \cdot \Delta\pd_{1}\psi_1 + \pd_{1}\bar{v}\, \nabla\big[ v(p'(v) - p'(\bar{v})) \big] \cdot \nabla\Delta\psi_1 \,dx ,      \\
    L_5(t) &:= -\dot{\X}(t) \int_{\Omega} \pd_{11}u^S_1 \cdot \Delta\pd_{1}\psi_1 \,dx ,       \\
    L_6(t) &:= \int_{\Omega} \nabla(vQ_3) \cdot \nabla\Delta\psi_1 \,dx ,      \\
    L_7(t) &:= -(2\mu+\lambda)\int_{\Omega} \nabla(v \sdot \pd_{11}u^R_1) \cdot \nabla\Delta\psi_1 \,dx ,       \\
    L_8(t) &:= (\mu+\lambda)\int_{\Omega} \pd_{i}\div\psi (\Delta\pd_{i}\psi - \nabla\pd_{i}\div\psi) \sdot \nabla v \,dx \\ 
    &\qquad - \int_{\Omega} \pd_{i}v \left( \mu\Delta\psi + (\mu+\lambda)\nabla\div\psi \right) \cdot \Delta\pd_{i}\psi \,dx .  
  \end{align*}
  It follows from Holder inequality and \eqref{eq:L3-norm} that
  \begin{align*}
    L_1(t) &\leq C\delta \norm{\nabla^2\psi}^2 + C \norm{\nabla\psi}_{L^3} \norm{\nabla^2\psi}_{L^6} \norm{\nabla^2\psi} \\
    &\leq C(\delta+\chi) \norm{\nabla^2\psi}^2 + C\chi \norm{\nabla^3\psi}^2 \\
    &\leq C(\delta+\chi) \norm{\nabla^2\psi}^2 + C\chi D_3(t).
  \end{align*}
  Notice that $vp'(v) = -\gamma p(v)$ and that
  \begin{equation*}
    |\nabla(p-\bar{p})| \lesssim |\nabla\phi| + |\pd_{1}\bar{v}| |\phi| \lesssim |\nabla\phi| + \delta |\phi| .
  \end{equation*}
  Thus, using Holder inequality, Sobolev inequality and \eqref{eq:L3-norm}, we have
  \begin{align*}
    L_2(t) &\leq C\delta \norm{\nabla\psi} \norm{\nabla\Delta\psi} + C \norm{\nabla\psi}_{L^3} \norm{\nabla\psi}_{L^6} \norm{\nabla\Delta\psi} \\
    &\quad + C\delta \norm{\nabla\phi} \norm{\nabla\Delta\psi} + \int_{\Omega} |\nabla(p-\bar{p})| \sdot |\nabla\phi| \sdot |\nabla\Delta\psi| \,dx \\
    &\leq C(\delta+\chi) \Big( D_3(t) + \norm{\nabla\psi}_{H^1}^2 + \norm{\nabla\phi}^2 + \norm{\nabla^2\phi}^2 \Big) \\
    &\leq C(\delta+\chi) \Big( D_3(t) + G^R(t) + G^S(t) + D(t) + \norm{\nabla\psi}_{H^1}^2 + \norm{\nabla^2\phi}^2 \Big).
  \end{align*}
  By Cauchy's inequality, it holds
  \begin{equation*}
    L_3(t) \leq \frac{1}{8}D_3(t) + C \norm{\nabla^2\phi}^2.
  \end{equation*}
  Notice that
  \begin{align*}
    &\nabla\big[ v(p'(v) - p'(\bar{v})) \big] = -\gamma \nabla(p-\bar{p}) - \nabla[p'(\bar{v})\phi] \\
    \lesssim&\; |\nabla(p-\bar{p})| + (\pd_{1}v^S + \pd_{1}v^R)|\phi| + |\nabla\phi| \lesssim |\nabla(p-\bar{p})| + (\pd_{1}v^S + \pd_{1}v^R)|\phi|.
  \end{align*}
  Thus, using Lemma \ref{lem:appro-rare-property}, Lemma \ref{lem:viscous shock} and \eqref{eq:L3-norm}, we have
  \begin{align*}
    L_4(t) &\leq C\int_{\Omega} |\pd_{11}\bar{v}| |\phi| |\nabla\Delta\psi| \,dx + C\delta \Big( D_3(t) + G^R(t) + G^S(t) + D(t) \Big) \\
    &\leq C(\delta+\sqrt{\delta_R})\Big(D_3(t) + G^R(t) \Big) + C\delta \Big( G^S(t) + D(t) \Big).
  \end{align*}
  By Lemma \ref{lem:viscous shock}, we have
  \begin{align*}
    L_5(t) &\leq C\delta_S |\dot{\X}(t)| \sqrt{D_3(t)} \norm{\pd_{1}v^S} \leq C\delta_S^{\frac{5}{2}} |\dot{\X}(t)| \sqrt{D_3(t)} \leq C\delta_S^{\frac{5}{2}} |\dot{\X}(t)|^2 + C\delta_S^{\frac{5}{2}} D_3(t).
  \end{align*}
  Using the decomposition of $vQ_3$ in \eqref{eq:Q_3 chaifen}, we divide $L_6(t)$ as follows:
  \begin{align*}
    L_6(t) = &\int_{\Omega} \nabla(Q_3^R) \cdot \nabla\Delta\psi_1 \,dx + \int_{\Omega} \nabla(Q_3^S) \cdot \nabla\Delta\psi_1 \,dx \\
    & + \int_{\Omega} \nabla(Q_3^I) \cdot \nabla\Delta\psi_1 \,dx =: \sum_{i=1}^3 L_{6,i}(t).
  \end{align*}
  Using Lemma \ref{lem:appro-rare-property}, assumption \eqref{eq:apriori-assumption} and \eqref{eq:L3-norm}, we yield
  \begin{align*}
    L_{6,1}(t) &\leq C\int_{\Omega} |\nabla\Delta\psi|\sdot \pd_{1}u^R_1 \big( |\psi_1|+|\phi|+|\nabla\psi_1|+|\nabla\phi| \big) \,dx  \\
    &\leq C \norm{\pd_{1}u^R_1}_{L^6} \norm{(\phi,\psi)}_{L^3} \norm{\nabla\Delta\psi} + C\delta_R \norm{\nabla(\phi,\psi)}^2 + C\delta_R \norm{\nabla\Delta\psi}^2 \\
    &\leq C\chi\norm{\pd_{1}u^R_1}_{L^6}^2 + C(\chi+\delta) \norm{\nabla(\phi,\psi)}^2 + C(\chi+\delta)D_3(t) \\
    &\leq C\chi\, \delta_R^{\frac{1}{3}}\sdot (1+t)^{-\frac{5}{3}} + C(\chi+\delta) \big( \norm{\nabla(\phi,\psi)}^2 + D_3(t) \big) \\
    &\leq C\chi\, \delta_R^{\frac{1}{3}}\sdot (1+t)^{-\frac{5}{3}} + C(\chi+\delta) \big( \norm{\nabla\psi}^2 + D_3(t) + D(t) + G^S(t) + G^R(t) \big) .
  \end{align*}
  It follows from Lemma \ref{lem:viscous shock} and Lemma \ref{lem:interaction} that
  \begin{align*}
    L_{6,2}(t) &\leq C\int_{\Omega} |\nabla\Delta\psi|\sdot \pd_{1}v^S \big( |\psi_1|+|\phi|+|\nabla\psi_1|+|\nabla\phi| \big) \,dx  \\
    &\leq C\int_{\Omega} |\nabla\Delta\psi|\sdot \pd_{1}v^S \left( \left| h_1-\bar{h}_1 - \frac{p-\bar{p}}{\s_*} \right| + |\nabla(p-\bar{p})| + |p-\bar{p}| + \pd_{1}u^R_1 \right) dx \\
    &\quad + C\delta \norm{\nabla(\phi,\psi)}^2 + C\delta \norm{\nabla\Delta\psi}^2 \\
    &\leq C\delta_R \delta_S e^{-C\delta_S t} + C\delta \big( G_2(t) + D(t) + G^R(t) + G^S(t) + D_3(t) + \norm{\nabla\psi}^2 \big).
  \end{align*}
  Similar to $K_{4,3}(t)$, we have
  \begin{equation*}
    L_{6,3}(t) \leq C\delta_R \delta_S e^{-C\delta_S t} + C\delta\big( D(t) + G^S(t) + G^R(t) + D_3(t) \big).
  \end{equation*}
  Using Lemma \ref{lem:appro-rare-property}, assumption \eqref{eq:apriori-assumption} and \eqref{eq:L3-norm}, we get
  \begin{align*}
    L_7(t) &\leq C\int_{\Omega} |\pd_{11}u^R_1| \sdot |\nabla\Delta\psi| \,dx + C\int_{\Omega} |\pd_{1}u^R_1| \sdot |\pd_{1}\bar{v} + \nabla\phi| \sdot |\nabla\Delta\psi| \,dx \\
    &\leq C \norm{\pd_{11}u^R_1} \norm{\nabla\Delta\psi} +  C \norm{\pd_{1}u^R_1}_{L^6} \norm{\nabla\phi}_{L^3} \norm{\nabla\Delta\psi}  \\
    &\leq C (\delta_R)^{-\frac{1}{2}} \norm{\pd_{11}u^R_1}^2 + C(\sqrt{\delta_R} + \chi) D_3(t) + C\chi \norm{\pd_{1}u^R_1}_{L^6}^2 .
  \end{align*}
  It follows Sobolev inequality and \eqref{eq:L3-norm} that
  \begin{align*}
    L_8(t) &\leq C \int_{\Omega} |\pd_{1}\bar{v} + \nabla\phi| \sdot |\nabla\div\psi| \sdot (|\nabla\Delta\psi| + |\nabla^2\div\psi|) \,dx  \\ 
    &\quad + C \int_{\Omega} |\pd_{1}\bar{v} + \nabla\phi| \sdot |\mu\Delta\psi + (\mu+\lambda)\nabla\div\psi| \sdot |\nabla\Delta\psi| \,dx   \\
    &\leq C\delta \left( \norm{\nabla^2\psi}^2 +D_3(t) \right) +
    C \norm{\nabla^2\div\psi} \norm{\nabla\phi}_{L^3} \norm{\nabla\div\psi}_{L^6} \\
    &\quad + C \norm{\nabla\Delta\psi} \norm{\nabla\phi}_{L^3} \left( \norm{\Delta\psi}_{L^6} + \norm{\nabla\div\psi}_{L^6} \right) \\
    &\leq C\delta \left( \norm{\nabla^2\psi}^2 +D_3(t) \right) + C\chi \left( \norm{\nabla^3\psi}^2 +D_3(t) \right) .
  \end{align*}
  Finally, notice that
  \begin{align*}
    D_3(t) \sim \norm{\nabla^3\psi}^2 .
  \end{align*}
  Hence, integrating \eqref{eq:psi-2-weifen} over $[0,t]$ and using the above estimates, we can deduce the desired inequality \eqref{eq:psi-2-est} from Lemma \ref{lem:psi-0-estimate}, Lemma \ref{lem:psi-1-estimate} and Lemma \ref{lem:phi-2-estimate}.
  The proof of Lemma \ref{lem:psi-2-estimate} is completed.
\end{proof}

\subsection{Proof of Proposition \ref{prop:apriori-estimate} }

Notice that
\begin{equation*}
  \norm{\nabla(v-\bar{v})}^2 \lesssim D(t) + \delta_S^2 G^S(t) + \delta_R G^R(t) .
\end{equation*}
Thus, using Lemma \ref{lem:psi-0-estimate}-\ref{lem:psi-2-estimate}, we immediately obtain \eqref{eq:apriori-estimate}.
In addition, by the definition of $\dot{\X}(t)$ \eqref{eq:X(t)-def}, it follows from Lemma \ref{lem:viscous shock} and assumption \eqref{eq:apriori-assumption} that
\begin{align*}
  |\dot{\X}(t)| \leq \frac{C}{\delta_S} \norm{v-\bar{v}}_{L^\infty} \int_{\Omega} \pd_{1}v^S \,dx \leq C \norm{v-\bar{v}}_{L^\infty},
\end{align*}
which implies \eqref{eq:xdot-property}.
The proof of Proposition \ref{prop:apriori-estimate} is completed.


\addtocontents{toc}{\setcounter{tocdepth}{2}}

\bibliographystyle{plain}
\bibliography{reference}

\begin{thebibliography}{10}

\bibitem{BD-2006}
Didier Bresch and Beno\^{i}t Desjardins.
\newblock On the construction of approximate solutions for the 2{D} viscous
  shallow water model and for compressible {N}avier-{S}tokes models.
\newblock {\em J. Math. Pures Appl. (9)}, 86(4):362--368, 2006.

\bibitem{Kang-ODE-2020}
Kyudong Choi, Moon-Jin Kang, Young-Sam Kwon, and Alexis~F. Vasseur.
\newblock Contraction for large perturbations of traveling waves in a
  hyperbolic-parabolic system arising from a chemotaxis model.
\newblock {\em Math. Models Methods Appl. Sci.}, 30(2):387--437, 2020.

\bibitem{Dafermos1996}
C.~M. Dafermos.
\newblock Entropy and the stability of classical solutions of hyperbolic
  systems of conservation laws.
\newblock In {\em Recent mathematical methods in nonlinear wave propagation
  ({M}ontecatini {T}erme, 1994)}, volume 1640 of {\em Lecture Notes in Math.},
  pages 48--69. Springer, Berlin, 1996.

\bibitem{DiPerna1979}
Ronald~J. DiPerna.
\newblock Uniqueness of solutions to hyperbolic conservation laws.
\newblock {\em Indiana Univ. Math. J.}, 28(1):137--188, 1979.

\bibitem{Goodman1986}
Jonathan Goodman.
\newblock Nonlinear asymptotic stability of viscous shock profiles for
  conservation laws.
\newblock {\em Arch. Ration. Mech. Anal.}, 95(4):325--344, 1986.

\bibitem{Goodman-1989}
Jonathan Goodman.
\newblock Stability of viscous scalar shock fronts in several dimensions.
\newblock {\em Trans. Amer. Math. Soc.}, 311(2):683--695, 1989.

\bibitem{Hoff-Zumbrun-2000}
David Hoff and Kevin Zumbrun.
\newblock Asymptotic behavior of multidimensional scalar viscous shock fronts.
\newblock {\em Indiana Univ. Math. J.}, 49(2):427--474, 2000.

\bibitem{Hoff-Zumbrun-2002}
David Hoff and Kevin Zumbrun.
\newblock Pointwise {G}reen's function bounds for multidimensional scalar
  viscous shock fronts.
\newblock {\em J. Differential Equations}, 183(2):368--408, 2002.

\bibitem{huang-li-Matsumura-2010}
Feimin Huang, Jing Li, and Akitaka Matsumura.
\newblock Asymptotic stability of combination of viscous contact wave with
  rarefaction waves for one-dimensional compressible {N}avier-{S}tokes system.
\newblock {\em Arch. Ration. Mech. Anal.}, 197(1):89--116, 2010.

\bibitem{Huang-Matsumura-2009}
Feimin Huang and Akitaka Matsumura.
\newblock Stability of a composite wave of two viscous shock waves for the full
  compressible {N}avier-{S}tokes equation.
\newblock {\em Comm. Math. Phys.}, 289(3):841--861, 2009.

\bibitem{Zumbrun-2017}
Jeffrey Humpherys, Gregory Lyng, and Kevin Zumbrun.
\newblock Multidimensional stability of large-amplitude {N}avier-{S}tokes
  shocks.
\newblock {\em Arch. Ration. Mech. Anal.}, 226(3):923--973, 2017.

\bibitem{Oleinik1960}
A.~M. Il'in and O.~A. Ole\u{\i}nik.
\newblock Asymptotic behavior of solutions of the {C}auchy problem for some
  quasi-linear equations for large values of the time.
\newblock {\em Mat. Sb. (N.S.)}, 51(93):191--216, 1960.

\bibitem{Ito-1996}
Kazuo Ito.
\newblock Asymptotic decay toward the planar rarefaction waves of solutions for
  viscous conservation laws in several space dimensions.
\newblock Number 966, pages 116--135. 1996.
\newblock Nonlinear evolution equations and their applications (Japanese)
  (Kyoto, 1995).

\bibitem{K-V-W-NSF}
Moon-Jin Kang, Alexis F., and Yi~Wang.
\newblock Time-asymptotic stability of generic {R}iemann solutions for
  compressible {N}avier-{S}tokes-{F}ourier equations, 2023.
\newblock arXiv: 2306.05604.

\bibitem{kang-lee-2024}
Moon-Jin Kang and Hobin Lee.
\newblock Long-time behavior toward composite wave of shocks for 3{D}
  barotropic {N}avier-{S}tokes system, 2024.
\newblock arXiv: 2406.11215.

\bibitem{KV-2017}
Moon-Jin Kang and Alexis~F. Vasseur.
\newblock {$L^2$}-contraction for shock waves of scalar viscous conservation
  laws.
\newblock {\em Ann. Inst. H. Poincar\'{e} C Anal. Non Lin\'{e}aire},
  34(1):139--156, 2017.

\bibitem{weighted-poincare}
Moon-Jin Kang and Alexis~F. Vasseur.
\newblock Contraction property for large perturbations of shocks of the
  barotropic {N}avier-{S}tokes system.
\newblock {\em J. Eur. Math. Soc. (JEMS)}, 23(2):585--638, 2021.

\bibitem{KVW-2019}
Moon-Jin Kang, Alexis~F. Vasseur, and Yi~Wang.
\newblock {$L^2$}-contraction of large planar shock waves for multi-dimensional
  scalar viscous conservation laws.
\newblock {\em J. Differential Equations}, 267(5):2737--2791, 2019.

\bibitem{Kang-Vasseur-Wang-2023}
Moon-Jin Kang, Alexis~F. Vasseur, and Yi~Wang.
\newblock Time-asymptotic stability of composite waves of viscous shock and
  rarefaction for barotropic {N}avier-{S}tokes equations.
\newblock {\em Adv. Math.}, 419:Paper No. 108963, 66, 2023.

\bibitem{Leger-2011}
Nicholas Leger and Alexis~F. Vasseur.
\newblock Relative entropy and the stability of shocks and contact
  discontinuities for systems of conservation laws with non-{BV} perturbations.
\newblock {\em Arch. Ration. Mech. Anal.}, 201(1):271--302, 2011.

\bibitem{3D-Gargliardo2022}
Hai-Liang Li, Teng Wang, and Yi~Wang.
\newblock Wave phenomena to the three-dimensional fluid-particle model.
\newblock {\em Arch. Ration. Mech. Anal.}, 243(2):1019--1089, 2022.

\bibitem{Li-teng-yi-2018}
Lin-an Li, Teng Wang, and Yi~Wang.
\newblock Stability of planar rarefaction wave to 3{D} full compressible
  {N}avier-{S}tokes equations.
\newblock {\em Arch. Ration. Mech. Anal.}, 230(3):911--937, 2018.

\bibitem{Liu-1985}
Tai-Ping Liu.
\newblock Nonlinear stability of shock waves for viscous conservation laws.
\newblock {\em Mem. Amer. Math. Soc.}, 56(328):v+108, 1985.

\bibitem{liu-xin-1997}
Tai-Ping Liu and Zhouping Xin.
\newblock Pointwise decay to contact discontinuities for systems of viscous
  conservation laws.
\newblock {\em Asian J. Math.}, 1(1):34--84, 1997.

\bibitem{Liu-Zeng-2015}
Tai-Ping Liu and Yanni Zeng.
\newblock Shock waves in conservation laws with physical viscosity.
\newblock {\em Mem. Amer. Math. Soc.}, 234(1105):vi+168, 2015.

\bibitem{Matsumura2018}
Akitaka Matsumura.
\newblock Waves in compressible fluids: viscous shock, rarefaction, and contact
  waves.
\newblock In {\em Handbook of mathematical analysis in mechanics of viscous
  fluids}, pages 2495--2548. Springer, Cham, 2018.

\bibitem{Matsumura-shock1985}
Akitaka Matsumura and Kenji Nishihara.
\newblock On the stability of travelling wave solutions of a one-dimensional
  model system for compressible viscous gas.
\newblock {\em Japan J. Appl. Math.}, 2(1):17--25, 1985.

\bibitem{Matsumura-rarefaction1986}
Akitaka Matsumura and Kenji Nishihara.
\newblock Asymptotics toward the rarefaction waves of the solutions of a
  one-dimensional model system for compressible viscous gas.
\newblock {\em Japan J. Appl. Math.}, 3(1):1--13, 1986.

\bibitem{Matsumura-rarefaction1992}
Akitaka Matsumura and Kenji Nishihara.
\newblock Global stability of the rarefaction wave of a one-dimensional model
  system for compressible viscous gas.
\newblock {\em Comm. Math. Phys.}, 144(2):325--335, 1992.

\bibitem{meng2025}
Jiayun Meng.
\newblock Time-asymptotic stability of composite weak planar waves for a
  general $n\times n$ multi-{D} viscous system, 2025.
\newblock arXiv: 2501.14188.

\bibitem{Nash1962}
John Nash.
\newblock Le probl\`eme de {C}auchy pour les \'{e}quations diff\'{e}rentielles
  d'un fluide g\'{e}n\'{e}ral.
\newblock {\em Bull. Soc. Math. France}, 90:487--497, 1962.

\bibitem{Yang-Zhao-2004}
Kenji Nishihara, Tong Yang, and Huijiang Zhao.
\newblock Nonlinear stability of strong rarefaction waves for compressible
  {N}avier-{S}tokes equations.
\newblock {\em SIAM J. Math. Anal.}, 35(6):1561--1597, 2004.

\bibitem{Nishikawa-2000}
Masataka Nishikawa and Kenji Nishihara.
\newblock Asymptotics toward the planar rarefaction wave for viscous
  conservation law in two space dimensions.
\newblock {\em Trans. Amer. Math. Soc.}, 352(3):1203--1215, 2000.

\bibitem{Solonnikov1976}
V.~A. Solonnikov.
\newblock The solvability of the initial-boundary value problem for the
  equations of motion of a viscous compressible fluid.
\newblock {\em Zap. Nau\v{c}n. Sem. Leningrad. Otdel. Mat. Inst. Steklov.
  (LOMI)}, 56:128--142, 197, 1976.
\newblock Investigations on linear operators and theory of functions, VI.

\bibitem{Szepessy-1993}
Anders Szepessy and Zhouping Xin.
\newblock Nonlinear stability of viscous shock waves.
\newblock {\em Arch. Ration. Mech. Anal.}, 122(1):53--103, 1993.

\bibitem{Teng-Yi-JEMS}
Teng Wang and Yi~Wang.
\newblock Nonlinear stability of planar viscous shock wave to three-dimensional
  compressible {N}avier-{S}tokes equations.
\newblock {\em J. Eur. Math. Soc. (JEMS)}, 2023.

\bibitem{Xin-1990}
Zhouping Xin.
\newblock Asymptotic stability of planar rarefaction waves for viscous
  conservation laws in several dimensions.
\newblock {\em Trans. Amer. Math. Soc.}, 319(2):805--820, 1990.

\end{thebibliography}

\end{document}